\documentclass[12pt,reqno]{amsart}
\usepackage{amssymb}
\usepackage[T1]{fontenc}
\usepackage[latin1]{inputenc} 
\usepackage{dsfont, mathrsfs}
\usepackage{a4wide}  
\usepackage{stmaryrd}
\usepackage{amsthm}
\usepackage{xspace}

\usepackage[colorlinks, hyperindex]{hyperref} 

\hypersetup{
    colorlinks=true,         
    linkcolor=blue,          
    citecolor=green,         
    filecolor=magenta,       
    urlcolor=cyan            
}

\newtheorem{theorem}{Theorem}[section]
\newtheorem{lemma}[theorem]{Lemma}

\newtheorem{corollary}[theorem]{Corollary}

\newtheorem{question}[theorem]{Question}

\theoremstyle{definition}
\newtheorem{definition}[theorem]{Definition}
\newtheorem{example}[theorem]{Example}

\theoremstyle{remark}
\newtheorem{remark}[theorem]{Remark}

\newcommand{\Cc}{\mathbb{C}}

\newcommand{\Ee}{\mathbb{E}}

\newcommand{\Pp}{\mathbb{P}}

\newcommand{\Rr}{\mathbb{R}}

\newcommand{\Uu}{\mathbb{U}}

\newcommand{\Un}{\mathds{1}}

\newcommand{\Ce}{\mathcal{C}}

\newcommand{\Fe}{\mathcal{F}}

\newcommand{\He}{\mathcal{H}}

\newcommand{\Le}{\mathcal{L}}

\newcommand{\Ue}{\mathcal{U}}

\newcommand{\Xe}{\mathcal{X}}
\newcommand{\Ze}{\mathcal{Z}}

\newcommand{\Hs}{\mathscr{H}}

\newcommand{\Ns}{\mathscr{N}}

\newcommand{\Us}{\mathscr{U}}

\newcommand{\Rg}{\mathfrak{R}}

\newcommand{\eg}{\mathfrak{e}}

\def\eqlaw{\stackrel{\Le}{=}}
\def\inv{^{-1}}
\def\longlongrightarrow{\hspace{+0.1ex} - \hspace{-1.1ex} - \hspace{-1.1ex} - \hspace{-1.1ex}\longrightarrow  }

\def\sign{\operatorname{sign}}

\def\Kol{\operatorname{Kol}}
\def\dkol{d_{\Kol}}

\def\zer{ ^{  _{(0)}  } }
\def\zzer{ ^{  _{(0)}  }_{\vphantom{ X }} }
\def\zerR{ ^{  _{(0, 0)}  } }

\def\zzerc{ ^{  _{(C)}  }_{\vphantom{ X }} }

\def\hn{ _{\vphantom{ \Le_{H_n}\inv } H_n} }

\def\hgam{ _{\vphantom{ \Le_{H_\gamma}\inv } H_\gamma} }

\def\geq{\geqslant}
\def\leq{\leqslant}

\let\oldforall\forall
\def\forall{\oldforall\,} 

\let\oldexists\exists
\def\exists{\oldexists\,}

\def\Re{\Rg \eg}

\newcommand{\ensemble}[1]{ \left\lbrace #1 \right\rbrace }
\newcommand{\prth}[1]{\!\left( #1 \right) } 
\newcommand{\Esp}[1]{ \Ee \prth{ #1 } }
\newcommand{\Prob}[1]{ \Pp \prth{ #1 } }
\newcommand{\crochet}[1]{\left[ #1 \right] }
\newcommand{\intcrochet}[1]{\llbracket #1 \rrbracket} 
\newcommand{\abs}[1]{\left| #1 \right|}
\newcommand{\norm}[1]{\left| \! \left| #1 \right| \! \right|}
\newcommand{\Unens}[1]{ \Un_{ \ensemble{#1} } }  
\newcommand{\tendvers}[2]{ \underset{#1 \rightarrow #2}{\longlongrightarrow} } 
\newcommand{\cvlaw}[2]{\stackrel{\Le}{\underset{#1 \, \rightarrow \, #2}{\longlongrightarrow}}} 
\newcommand{\cvmodg}[2]{   \   {\underset{#1 \rightarrow #2}{  \overset{ \operatorname{mod-G} }{\longlongrightarrow } } } \ }
\newcommand{\Bgammax}[1]{ \rho_\gamma''(#1) + 3 \rho_\gamma(#1) \rho_\gamma'(#1)  + \rho_\gamma(#1)^3 }

\newcommand{\emailhref}[1]{ \email{\href{mailto:#1}{#1}} }

\begin{document}


\title[On Stein's method and mod-* convergence]{On Stein's method and mod-* convergence} 


\author[Y. Barhoumi-Andr\'eani]{Yacine Barhoumi-Andr\'eani}
\address{Department of Statistics, University of Warwick, Coventry CV4 7AL, U.K.}
\emailhref{y.barhoumi-andreani@warwick.ac.uk}


\subjclass[2000]{60E10, 60E05, 60F05, 60G50, 60B10}


\date{\today}



\begin{abstract}
Stein's method allows to prove distributional convergence of a sequence of random variables and to quantify it with respect to a given metric such as Kolmogorov's (a Berry-Esséen type theorem). Mod-* convergence quantifies the convergence of a sequence of random variables to a given distribution in a sense unusual in probability theory, a priori unrelated to a metric on probability measures. 

This article gives a connection between these two notions. It shows that mod-* convergence can be understood as a higher order approximation in distribution when the limiting function is integrable and proves a refined Berry-Esséen type theorem for sequences converging in the mod-Gaussian sense.
\end{abstract}

\maketitle

\vspace{-0.9cm}

\tableofcontents

\section{Introduction}

Let  $ (X_n)_n $ be a sequence of random variables converging in law to $ Z \sim \Ns(0,1) $ ; for instance, take $ X_n $ to be the sum of $ n $ i.i.d. random variables of expectation $0$ and variance $ 1/ n $. The Central Limit Theorem asserts that 
\begin{align*}
\dkol(X_n, Z) := \sup_{x \in \Rr} \abs{ \Prob{X_n \leq x} - \Prob{ Z \leq x } } \tendvers{n}{ +\infty } 0
\end{align*}

The Berry-Esséen theorem \cite{Berry, Esseen} is a direct continuation of the Central Limit Theorem~: it gives the rate of convergence of this latest limit under the form
\begin{align*}
\dkol(X_n, Z)  \leq \frac{C}{ \sqrt{n} }
\end{align*}
with a constant $ C $ depending on the sequence $ (X_n)_n $.

To prove their bound, Berry and Esséen used a Fourier inversion. Such a method perfectly applies in the framework of a sum of independent variables but becomes less efficient in the context of a marked dependence.

Charles Stein introduced his eponymous method in \cite{SteinProceedings} as an alternative to the Fourier formalism to achieve a Berry-Esséen bound. The key point consisted in replacing the characteristic function by a characteristic operator easier to handle in situations of dependency. Many paradigm shifts were then observed in the theory ; initially designed for the Gaussian distribution, the method was extended to the Poisson setting in \cite{Chen} and the characterisation of the distribution via the operator was replaced by a fixed point equation in law using a probabilistic transformation such as the $0$-bias or the size-bias transform (see \cite{GoldsteinReinert, GoldsteinRinott}).

$ $

Mod-Gaussian convergence was introduced in \cite{JacodAl}. A sequence of random variables $ (X_n) _n $ is said to converge in the mod-Gaussian sense if there exists a sequence $ (\gamma_n)_n $ of strictly positive reals and a function $ \Phi : \Rr \to \Cc $ such that, locally uniformly in $ u \in \Rr $ 
\begin{align*}
\frac{ \Esp{ e^{i u (X_n - \Ee(X_n)) } } }{ e^{- u^2 \gamma_n^2 / 2} } \tendvers{n }{ + \infty } \Phi(u)
\end{align*}

Due to the type of convergence and the properties of the converging sequence, $ \Phi $ is a continuous function satisfying $ \Phi(0) = 1 $ and $ \overline{\Phi(u)} = \Phi(-u) $. Moreover, $ \Phi $ is not necessarily the Fourier transform of a probability distribution (see \cite{JacodAl}). This last fact impeds the naive probabilistic interpretation that would think $ X_n - \Esp{X_n} $ as the sum of a random variable and a Gaussian noise.

The same notion can be defined in the Poisson framework (see \cite{KowalskiNikeghbali1})~: a sequence $ (Z_n)_n $ is  is said to converge in the mod-Poisson sense at speed $ (\gamma_n)_n $ if, locally uniformly in $x \in \Uu $, 
\begin{align*}
\frac{ \Esp{ x^{  Z_n  } } }{ e^{ \gamma_n (x-1) } } \tendvers{n }{ + \infty } \Phi(x)
\end{align*}
for a continuous function $ \Phi : \Uu \to \Rr $ satisfying $ \Phi(1) = 1 $. Here, $ \Uu $ designates the unit circle. 

These two notions define mod-* convergence with * $ \in \ensemble{ \operatorname{Gaussian}, \operatorname{Poisson} }$, but the set of admissible distributions can also extend to the infinitely divisible case (see \cite{DelbaenAl}) or any distribution that arises as a limit in law whose Fourier transform does not vanish.

As we consider unormalised (hence diverging) random variables, renormalising directly their Fourier transform gives a non trivial limiting function. Note that a change of renormalisation (setting $ u = v/\gamma_n $) implies the convergence in law of $ (X_n - \Esp{X_n})/\gamma_n $ to the Gaussian distribution~:
\begin{align*}
\frac{ \Esp{ e^{i v (X_n - \Ee(X_n))/\gamma_n } } }{ e^{- v^2  / 2} } \tendvers{n }{ + \infty } 1
\end{align*}

This former type of convergence is thus more precise than the usual convergence in law. It is unusual in probability theory and was only investigated in a few occurrences, for instance in \cite{Hwang1}. But it is well exploited in other branches of mathematics such as number theory since, for example, Keating and Snaith's celebrated moments conjecture writes (see \cite{KeatingSnaith})
\begin{align*}
\frac{ \Esp{ e^{ \lambda \log\abs{ \zeta\prth{ \frac{1}{2} + i T U } } } } }{ e^{  (\frac{1}{2} \log\log T)  \lambda^2 / 2 } } \tendvers{ T }{ + \infty } \Phi_\zeta(\lambda)
\end{align*}

Here, $ \zeta $ denotes the Riemann Zeta function, $ \Phi_\zeta $ is a particular function described in \cite{KeatingSnaith} and $ U $ denotes a random variable uniformly distributed in $ \crochet{0, 1} $. The convergence holds locally uniformly in $ \lambda \in \ensemble{ \Re > -1 } $, hence for $ \lambda \in i \Rr $. The fact that $ \Phi_\zeta $ does not write as the Fourier-Laplace transform of a probability distribution asks the question of the existence of a random variable that would naturally converge in the mod-Gaussian sense to the same function, to be able to compare it to the original sequence of random variables.


$ $

The first result of this paper answers the following

\begin{question}\label{Question:Mod} Given a function $ \Phi : \Rr \to \Cc $ satisfying some admissibility assumptions (described in detail in theorem \ref{Thm:InterprModG}) and $ \gamma \to \infty $, can we construct a family of random variables $ (\He_\gamma(\Phi))_\gamma $ that converges in the mod-Gaussian sense to $ \Phi $ ?

\end{question}

The answer to this question is given in theorem \ref{Thm:InterprModG}. A slight shift of point of view, using the Laplace transform in place of the Fourier transform in the precedent definition allows to already give a flavour of the result~:

\begin{theorem}\label{Thm:Intro} Let $ \Phi : \Rr \to \Rr $ be a continuous positive integrable function on $ \Rr $ such that $ \Phi(0) = 1 $. Let $ \gamma > 0 $ and $ X_\gamma \sim \Ns(0, \gamma^2) $. Define the random variable $ \He_\gamma(\Phi) $ by the change of probability 
\begin{align}\label{Def:VariableHIntro}
\Esp{ f(\He_\gamma(\Phi) ) } := \frac{  \Esp{ f(X_\gamma) \Phi\prth{ \frac{ X_\gamma  }{\gamma^2} } }  }{  \Esp{  \Phi\prth{ \frac{ X_\gamma  }{\gamma^2} } } }
\end{align}
for all bounded measurable $ f : \Rr \to \Rr $. Then, 
\begin{align}\label{Eq:DualiteModGaussienne}
\frac{  \Esp{ e^{ u \He_\gamma(\Phi) } } }{ \Esp{ e^{ u X_\gamma } } } = \frac{ \Esp{ \Phi\prth{ \frac{X_\gamma}{\gamma^2} + u } } }{ \Esp{ \Phi\prth{ \frac{X_\gamma}{\gamma^2} } } }
\end{align}
and in particular, locally uniformly in $ u \in \Rr $
\begin{align*}
\frac{  \Esp{ e^{ u \He_\gamma(\Phi) } } }{ \Esp{ e^{ u X_\gamma } } } \tendvers{ \gamma }{ + \infty }  \Phi(u)
\end{align*}
\end{theorem}


The duality relation \eqref{Eq:DualiteModGaussienne} will appear to be a direct avatar of the Gaussian change of probability. More generally, it will hold in the context of any infinitely divisible distribution (see theorem \ref{Thm:InfDiv}).

Since there exists now a ``canonical'' random variable associated to $ \Phi $ (with some additional restrictive hypotheses, though), it is tempting to think of \textit{metrising} mod-Gaussian convergence in this restricted setting by performing a probabilistic approximation with this new distribution. This asks the

\begin{question}\label{Question:Bound} Given a sequence of random variables $ (X_n)_n $ that converges in the mod-Gaussian sense at speed $ (\gamma_n)_n $ to a given function $ \Phi $ satisfying the hypotheses of theorem \eqref{Thm:Intro}, can we find a bound for
\begin{align*}
\dkol\prth{ X_n, \He_{\gamma_n }(\Phi) } := \sup_{x \in \Rr} \abs{ \Prob{X_n \leq x} - \Prob{ \He_{\gamma_n }(\Phi) \leq x } } 
\end{align*}
and compare it to the classical Berry-Ess\'een bound obtained for the convergence in law of $ (X_n/\gamma_n)_n $ to the Gaussian distribution ?
\end{question}

The answer to question \ref{Question:Bound} will use Stein's method, by first describing a characteristic operator of the law of $ \He_{\gamma_n}(\Phi) $, and then following Stein's steps in \cite{SteinApprox}. We will address the problem with complex analytic methods in a subsequent publication, using the information on the values of $ \Phi $ and its derivatives in $ 0 $ in the same vein as Berry and Ess\'een.

The most famous sequence of random variables converging in law to the Gaussian distribution is the sum of i.i.d. random variables. Such a sequence is also convergent in the mod-Gaussian sense to an explicit function satisfying the hypotheses of theorem \ref{Thm:Intro}, namely $ \Phi_C : x \mapsto e^{-C x^4/4} $ for a certain constant $ C $ (see example \ref{Ex:SommeIID}). We will treat this example in section \ref{Chapitre:ExempleIID}. The distribution of $ \He_\gamma(\Phi_C) $ is the Gaussian with quartic interaction potential, a real-valued version of the celebrated $ \Phi^4 $ model given by $ \Prob{ \He_\gamma(\Phi_C) \leq x } = \int_{-\infty}^x \exp\big( - a \frac{t^2}{2} - b \frac{t^4}{4} \big) dt /\Ze_{a, b}  $ where $ a, b $ and $ \Ze_{a, b} $ are explicitely described in section \ref{Chapitre:ExempleIID}. The Stein estimates developed in this context in section \ref{Chapitre:SteinEstimates} can thus be of independent use when one deals with such a distribution.

$ $

This article is structured in the following way : section \ref{Chapitre:Mod} defines mod$\operatorname{-}$* convergence and constructs a canonical family of distributions associated to it, section \ref{Chapitre:Stein} explains the links between this canonical family and Edgeworth expansion or signed expansion and develops the fondamentals of Stein's method for a sequence of random variables converging in the mod-Gaussian sense to a general function $ \Phi $ ; last, section \ref{Chapitre:ExempleIID} treats the example of the sum of i.i.d. symmetric random variables and proves a Berry-Ess\'een theorem using Stein's method and the estimates of section \ref{Chapitre:SteinEstimates}.


\section*{Notations}

We gather here some notations used throughout the paper.

If $ \gamma $ is a positive real number, $ \Ns(0, \gamma^2) $ designates the Gaussian distribution of expectation 0 and variance $ \gamma^2 $ and $ \Us\prth{I} $ the uniform distribution in the set $ I $. 
If $n$ is an integer, $ \intcrochet{1, n} $ designates the set $ \ensemble{1, 2, \dots, n} $ ; $ z \mapsto \overline{z} $ denotes the complex conjugation.

All random variables will be considered on a probability space $ (\Omega, \Fe, \Pp) $. The distribution of the random variable $ X : \Omega \longrightarrow \Rr $ will be denoted by $ \Pp_X $~: if $ A \in \Fe $ is a measurable set, $ \Pp_X(A) := \Prob{X \in A } $. If $X$ and $Y$ are two random variables having the same distribution, that is $ \Pp_X = \Pp_Y $, we will note $ X \eqlaw Y $. The convergence in law/in distribution will be denoted by $ \stackrel{\Le}{\longrightarrow} $. 

For $ f \in L^1(\Pp_X) $, $ f \geq 0 $, the penalisation or bias of $ \Pp_X $ by $f$ is the probability measure $ \Pp_Y $ denoted by
\begin{align*}
\Pp_{ Y } := \frac{ f(X) }{ \Esp{ f(X) } } \bullet \Pp_{ X }
\end{align*}

This definition is equivalent to the following~: for all $ g \in L^\infty(\Pp_X) $, 
\begin{align*}
\Esp{ g(Y) } = \frac{ \Esp{ f(X) g(X) }  }{ \Esp{ f(X) } } 
\end{align*}


\section{Mod-* convergence}\label{Chapitre:Mod}

\subsection{Reminder of the main notions}

\begin{definition}[Mod-Gaussian convergence] Let $ (X_n)_n $ be a sequence of random variables of expectation $ 0 $ and $ (\gamma_n)_n $ be a sequence of strictly positive real numbers. Let $ G \sim \Ns(0,1) $. We say that $ (X_n)_n $ converges in the mod-Gaussian sense if 
\begin{align*}
\frac{ \Esp{ e^{i u X_n} } }{ \Esp{ e^{i u \gamma_n G } } } \tendvers{n }{ + \infty } \Phi(u)
\end{align*}
the convergence being locally uniform in $ u \in \Rr $ and $ \Phi : \Rr \to \Cc $ hence being a continuous function satisfying $ \Phi(0) = 1 $ and $ \overline{\Phi(u)} = \Phi(-u) $.

When such a convergence holds, we write it as
\begin{align*}
(X_n, \gamma_n ) \cvmodg{n }{ + \infty } \Phi
\end{align*}

\end{definition}


\begin{remark} 

One can always be reduced to the case of a sequence of random variables with zero expectation. Otherwise, we include additional renormalization in the Fourier transform of the Gaussian random variable, which corresponds to the original definition of \cite{JacodAl}.

\end{remark}


A trivial example but a useful insight for the intuition allows to illustrate the concept~:

\begin{example} Consider $ X_n := Y_n + \gamma_n G $ where $ (Y_n)_n $ is independent of $ G \sim \Ns(0,1) $, where $ Y_n \cvlaw{n}{ + \infty} Y_\infty $ and $ \gamma_n \to + \infty $. Then,
\begin{align*}
\frac{ \Esp{ e^{i u X_n} } }{ \Esp{ e^{i u \gamma_n G } } } = \Esp{e^{i u Y_n } } \tendvers{n }{ + \infty } \Phi(u) := \Esp{e^{i u Y_\infty } }
\end{align*}

Thus, in the case of an \textit{additive} independent Gaussian noise, such a renormalisation gives at the limit the Fourier transform of a probability measure.

\end{example}


An interesting question related to question \ref{Question:Mod} concerns the probabilistic meaning of this particular type convergence.  Since the limiting function is not always the Fourier transform of a probability measure (see e.g. example \ref{Ex:SommeIID}), the intuitive idea of an additive correlated noise that disappears with this particular type of renormalisation (a deconvolution) is not satisfactory if we escape from the domain of probability theory at the limit. 

As pointed out in the introduction, one solution to this problem is to change the probability using $ \Phi $ as a weight~:

\begin{theorem}[A probabilistic interpretation of mod-Gaussian convergence]\label{Thm:InterprModG}

Let $ (\gamma_n)_n $ be a sequence of strictly positive real numbers such that $ \gamma_n \to + \infty $ when $ n \to +\infty $ and let $ \Phi $ be an admissible function for the mod-Gaussian convergence, i.e. a continuous complex function satisfying $ \Phi(0) = 1 $ and $ \overline{\Phi(u)} = \Phi(-u) $.

Suppose moreover that 

\begin{enumerate}

\item $ \Phi $ can be analytically extended on the whole complex plane and satisfies, $ \forall \beta \in \Rr $,
\begin{align}\label{Ineq:SupBande}
&  \sup_{z \in  a + i \crochet{0, \beta} } \abs{ \Phi(z) } < \infty  \ \ \ \forall a \in \Rr  \notag \\
&  \sup_{z \in  a + i \crochet{0, \beta} } \abs{ \Phi(z) } \tendvers{ a }{  \pm \infty } 0  
\end{align}

\item $ \Phi(i x ) = \Phi(x) $ for all $ x \in \Rr $,

\item $ \Phi(x ) \geq 0 $ for all $ x \in \Rr $,

\end{enumerate}

Define the distribution $ \Pp_{\He_{\gamma_n}(\Phi)} $ of a random variable $ \He_{\gamma_n}(\Phi) $ by the following penalisation
\begin{align}\label{Def:DistribH}
\Pp_{ \He_{\gamma_n}(\Phi) } := \frac{ \Phi\prth{ \frac{G}{\gamma_n} } }{ \Esp{ \Phi\prth{ \frac{G}{\gamma_n} } } } \bullet \Pp_{ \gamma_n G }
\end{align}

Then, 
\begin{align*}
\He_{\gamma_n}(\Phi)/\gamma_n & \cvlaw{n}{\infty}   \Ns(0,1) \\
(\He_{\gamma_n}(\Phi), \gamma_n) & \cvmodg{n}{\infty}   \Phi 
\end{align*}
\end{theorem}

Note that \eqref{Def:DistribH} coincides with \eqref{Def:VariableHIntro}, but in the Fourier framework, which imposes $ \Phi $ to be positive and real on $ \Rr $. We will compare these two settings in section \ref{Section:LaplaceSetting}. Before proving the theorem, we give an example of such a function $ \Phi $ that will be our guiding example.

\begin{example}\label{Ex:SommeIID} For $ C > 0 $, set
\begin{align}\label{Def:FonctionPhiC}
\Phi_C(x) = e^{- C \frac{x^4}{4} }
\end{align}

This function is the mod-Gaussian limit of 
\begin{align*}
Z_n := \frac{1}{ n^{1/4}  } \sum_{k = 1}^n X_k
\end{align*}
where $ (X_k)_k $ is a sequence of i.i.d. symmetric random variables (that is $ X \eqlaw -X $) satisfying $ \Esp{X^2} = 1 $ and $ \kappa := \Esp{X^4}  < 3 $.

More precisely, for $ C = (3 - \Esp{ X^4 } )/6 $, we have mod-Gaussian convergence of $ ( Z_n )_n $ to $ \Phi_C $ at speed $  n^{1/4}  $ :
\begin{align*}
\Esp{e^{ix Z_n} } & =  \Esp{e^{ix \sum_{k = 1}^n X_k / n^{1/4} } } = \prth{\Esp{e^{ix X / n^{1/4} } }}^n \\
				  		  & =   e^{ n \log\,\prth{ \Ee\,\prth{ \exp\,\prth{ix X / n^{1/4}} } } } \\ 
				  		  & =  e^{ n \log\,\prth{ 1 +  \frac{x^2}{ 2 \sqrt{n} } + \kappa \frac{x^4}{ 24 n } + \frac{x^4}{ n } \varepsilon_1\prth{ \frac{x}{n^{1/4}} }  } } \\
				  		  & =  e^{ n \, \prth{ \frac{x^2}{ 2 \sqrt{n} } + \kappa \frac{x^4}{ 24 n } - \frac{1}{2} \prth{ \frac{x^2}{ 2 \sqrt{n}} }^2  + \frac{x^4}{ n } \varepsilon_2\prth{ \frac{x}{n^{1/4}} }  } } \\
				  		  & =  e^{  \sqrt{n} \frac{x^2}{ 2  } + (\kappa - 3) \frac{x^4}{ 24  }  +  \frac{x^4}{ n } \varepsilon_2\prth{ \frac{x}{n^{1/4}} } }  
\end{align*}

Here, $ \varepsilon_1 $ and $ \varepsilon_2 $ are functions that tend to $0$ in $0$ and are bounded on a compact neighborhood of $0$. We thus have the following convergence that holds locally uniformly in $x$ in a certain interval around $0$
\begin{align*}
\frac{  \Esp{e^{ix Z_n} }  }{ \Esp{e^{ix n^{1/4} G } } }   \tendvers{n}{ + \infty }  e^{ - \frac{ (3 - \kappa) }{ 24 } x^4}
\end{align*}

We can check moreover that the required assumptions of analyticity and boundedness in a horizontal strip are fullfilled~: 
\begin{align*}
\sup_{z \in  a + i \crochet{0, \beta} } \abs{ e^{-C z^4} } = \sup_{ y \in  \crochet{0, \beta} } \abs{ e^{-C (a + i y)^4} } =  \sup_{ y \in  \crochet{0, \beta} } e^{-C ( a^4 + y^4 - 6 a^2 y^2)}   \leq  C'_\beta e^{-C a^4 /2 } \tendvers{ a }{  \pm \infty } 0 
\end{align*}

In accordance to theorem \ref{Thm:InterprModG}, the random variable $ \He_{ n^{1/4} }(\Phi_C) $ of distribution given by \eqref{Def:DistribH} with $ \Phi_C(x) = e^{- C x^4/4 } $ satisfies the same type of convergence, and one can write
\begin{align*}
\frac{ \Esp{ e^{iu Z_n} }  }{ \Esp{ e^{iu \He_{ n^{1/4} }(\Phi_C) } } } \tendvers{ n }{  + \infty } 1 
\end{align*}
in the same vein as convergence in law writes
\begin{align*}
\frac{ \Esp{ e^{iu Z_n / n^{1/4} } }  }{ \Esp{ e^{iu G} } } \tendvers{ n }{  + \infty } 1 
\end{align*}
\end{example}
We now prove theorem \ref{Thm:InterprModG}.

\begin{proof}
For $ \theta \in \Rr $, write
\begin{align*}
\frac{ \Esp{e^{i \theta \He_{\gamma_n}(\Phi) } } }{\Esp{ e^{i \theta \gamma_n G} } } =  \frac{ \Esp{  \Phi\prth{ \frac{G}{\gamma_n} } e^{i \theta \gamma_n G } } }{  \Esp{e^{i \theta \gamma_n G } } \Esp{  \Phi\prth{ \frac{G}{\gamma_n} } } } =  \frac{   \Esp{  \frac{ \vphantom{ \big ( } e^{i \theta \gamma_n G } }{  \Esp{  \vphantom{ \big ( } e^{i \theta \gamma_n G } } }   \Phi\prth{ \frac{G}{\gamma_n} }    }  }{ \Esp{  \Phi\prth{ \frac{G}{\gamma_n} } } }  =: \frac{ \int_\Rr \Phi(x) \mu_n^{(\theta)}(dx) }{ \Esp{  \Phi\prth{ \frac{G}{\gamma_n} } }}
\end{align*}
where
\begin{align*}
\int_\Rr \Phi(x) \mu_n^{(\theta)}(dx)  & := \Esp{   \frac{ e^{i \theta \gamma_n G } }{ \Esp{e^{i \theta \gamma_n G } } }   \Phi\prth{ \frac{G}{\gamma_n} }    } \\
				& = e^{ \theta^2 \gamma_n^2 / 2 } \int_\Rr  e^{ i \theta \gamma_n x } \Phi\prth{ \frac{ x }{ \gamma_n } } e^{-x^2 / 2 } \frac{dx}{ \sqrt{2 \pi} } \\
				& =  \int_\Rr  \Phi\prth{ \frac{ x }{ \gamma_n } }  e^{ -\frac{1}{2} ( x - i \theta \gamma_n )^2 } \frac{dx}{ \sqrt{2 \pi} } \\
				& = \int_{\Rr - i \theta \gamma_n }  \Phi\prth{ \frac{ y }{ \gamma_n } + i \theta } e^{ -\frac{1}{2} y^2 } \frac{dy}{ \sqrt{2 \pi} } 
\end{align*}

Set
\begin{align*}
g(z) := \Phi(z/\gamma_n + i\theta) e^{-z^2 / 2} 
\end{align*}

If $g$ is analytic on the whole complex plane, the Cauchy formula gives
\begin{align*}
\int_{ \crochet{-a, a}  } g + \int_{ a + i \crochet{ 0, \beta }  } g - \int_{ \crochet{-a, a} + i \beta  } g - \int_{ -a + i \crochet{ 0, \beta }  } g = 0
\end{align*}

If moreover $g$ satisfies the hypothesis \eqref{Ineq:SupBande}, we can write
\begin{align*}
\abs{ \int_{ a + i \crochet{ 0, \beta }  } g(x) dx }  \leq \abs{ \beta } \sup_{z \in  a + i \crochet{ 0, \beta } } \abs{ g(z) }  \tendvers{ a }{  \pm \infty } 0 
\end{align*}

Hence, 
\begin{align*}
\int_{ \crochet{-a, a} + i \beta  } g = \int_{ \crochet{-a, a}  } g +  \prth{\int_{ a + i \crochet{ 0, \beta }  } g - \int_{ -a + i \crochet{ 0, \beta }  } g } =:  \int_{ \crochet{-a, a}  } g + R(a)
\end{align*}
with 
\begin{align*}
\abs{ R(a) }  \leq 2\abs{ \beta } \sup_{z \in  a + i \crochet{ 0, \beta } } \abs{ g(z) }  \tendvers{ a }{  \pm \infty } 0 
\end{align*}

Passing to the limit on $ a \to + \infty $, we get 
\begin{align*}
\int_{ \Rr - i \beta  } g = \int_{ \Rr  } g 
\end{align*}

Now, 
\begin{align*}
& \sup_{z \in  a + i \crochet{ 0, \beta } } \abs{ e^{ -z^2 / 2 } }   =  \sup_{u \in \crochet{ 0, \beta } } \abs{ e^{ -(a + iu)^2 / 2 } }   =  \sup_{u \in \crochet{ 0, \beta } }  e^{ -a^2 / 2 + u^2 / 2 }   = e^{ \beta^2 / 2 }  e^{ -a^2 / 2 }  \tendvers{ a }{  \pm \infty } 0  \\
& \sup_{z \in  a + i \crochet{ 0, \beta } } \abs{ \Phi(z) }   =  \sup_{u \in \crochet{ 0, \beta } } \abs{ \Phi(a + iu) } \tendvers{ a }{  \pm \infty } 0  \mbox{ $ $ by the hypothesis \eqref{Ineq:SupBande} } \\
& \sup_{z \in  a + i \crochet{ 0, \beta } } \abs{ e^{ -z^2 / 2 } \Phi\prth{ z } } \tendvers{ a }{  \pm \infty } 0 
\end{align*}

We can thus write 
\begin{align*}
\int_\Rr \Phi(x) \mu_n^{(\theta)}(dx) = \int_{\Rr }  \Phi \prth{ \frac{ y }{ \gamma_n } + i \theta } e^{ -\frac{1}{2} y^2 } \frac{dy}{ \sqrt{2 \pi} } = \Esp{ \Phi\prth{ \frac{G}{\gamma_n} + i \theta } }
\end{align*}

The condition \eqref{Ineq:SupBande} ensures that $ \Phi $ is bounded on a horizontal strip, hence, by the dominated convergence theorem, the continuity of $ \Phi $ on the complex plane and the hypothesis $ \Phi(i\theta) = \Phi(\theta) $ for all $ \theta \in \Rr $, we get 
\begin{align*}
\lim_{n \to + \infty} \int_\Rr \Phi(x) \mu_n^{(\theta)}(dx) = \Esp{  \lim_{n \to + \infty} \Phi\prth{ \frac{G}{\gamma_n} + i \theta } }  = \Phi(i \theta) = \Phi(\theta) 
\end{align*}

Finally, dominated convergence implies
\begin{align*}
\lim_{n \to + \infty} \Esp{ \Phi \prth{ \frac{ G }{\gamma_n} } } = \Phi(0) = 1
\end{align*}
which proves the theorem.
\end{proof}


\begin{remark}\label{Ex:CvFaibleInterpr} The fact that the signed (complex) measures $ \mu_n^{(\theta)} $ satisfy
\begin{align*}
\lim_{n \to + \infty} \int_\Rr \Phi(x) \mu_n^{(\theta)}(dx) = \Phi(\theta) =  \int_\Rr \Phi(x) \delta_\theta(dx)
\end{align*}
for all $ \Phi $ satisfying the assumptions of theorem \ref{Thm:InterprModG} can be rephrased into a weak convergence of the sequence $ ( \mu_n^{(\theta)} )_n $ to the measure $ \delta_\theta $. Note that the space of functions on which this convergence holds is restrictive and is a strict subset of the space of continuous bounded functions. On this last space, the weak convergence does not hold as one can check by considering the limit of the Fourier transform $ \int_\Rr e^{i \alpha x} \mu_n^{(\theta)}(dx) $.
\end{remark}

The last theorem motivates the following 

\begin{definition} Let $ G \sim \Ns(0, 1) $, $ \gamma > 0 $ and $ \Phi $ be a function satisfying the hypotheses of theorem \ref{Thm:InterprModG}. We define the distribution $ \Hs(\Phi, \gamma) $ by 
\begin{align}\label{Def:LoiHgamma}
H_\gamma \sim \Hs(\Phi, \gamma)   \ \ \ \ \Longleftrightarrow \ \ \ \  \Pp_{ H_\gamma } := \frac{ \Phi\prth{ \frac{G}{\gamma} } }{ \Esp{ \Phi\prth{ \frac{G}{\gamma} } } } \bullet \Pp_{ \gamma G }
\end{align}
\end{definition}

\begin{remark} Another way of writing \eqref{Def:LoiHgamma} is to say that $ H_\gamma $ has a Lebesgue-density given by
\begin{align}\label{Eq:HLebesgueDensite}
f_\gamma(x) = \frac{1}{c_\gamma } \Phi \prth{ \frac{x}{\gamma^2} } e^{ -\frac{1}{2}  \prth{ \frac{x}{\gamma} }^2 },  \quad\quad c_\gamma := \gamma \sqrt{2\pi} \Esp{ \Phi ( G/\gamma ) }
\end{align}
%
%
%
%
\end{remark}

\subsection{Mod-Gaussian convergence in the Laplace setting}\label{Section:LaplaceSetting}

As noticed in remark \ref{Ex:CvFaibleInterpr}, the key point in theorem \ref{Thm:InterprModG} is to show that $ (\mu_n^{(\theta)})_n $ converges weakly to  $ \delta_{ \theta } $ for a certain notion of weak convergence of measures. But the fact that $ \lim_{n \to + \infty} \int_\Rr \Phi(x) \mu_n^{(\theta)}(dx) = \Phi(i\theta) $ forces the function $ \Phi $ to have an additionnal symmetry and gives the hint that this is the variable $ i \theta $ that should be the relevant parameter. It thus becomes natural to consider the Laplace transform in place of the Fourier transform.

\begin{definition}\label{Def:ModGLaplace} Let $ (X_n)_n $ be a sequence of random variables of expectation $ 0 $ and  $ (\gamma_n)_n $ a sequence of strictly positive real numbers. Suppose moreover that $ \Esp{ e^{u X_n} } < \infty $ for all $ u \in A \subset \Rr $ where $ A $ is an open set containing $ 0 $ or $ A = \Rr_+ $

$ (X_n)_n $ is said to converge in the mod-Gaussian-Laplace sense at speed $ (\gamma_n)_n $ if 
\begin{align*}
\frac{ \Esp{ e^{  u X_n} } }{ \Esp{ e^{ u \gamma_n G } } } \tendvers{n }{ + \infty } \Phi\prth{ u }
\end{align*}
where $ \Phi : A \to \Rr_+ $ is a continuous function satisfying $ \Phi(0) = 1 $, the last convergence being locally uniform in $ u \in A $.
\end{definition}


\begin{remark}
Note that the function $ \Phi $ here defined must always be positive, as a limit of a sequence of positive functions. The advantage of choosing the Fourier transform in place of the Laplace transform is clear~: the former one always exists when the latter one needs to specify the range of $ u \in \Rr $ where it is defined. But for the purpose that we have set, a real function is more suited.
\end{remark}
We now prove theorem \ref{Thm:Intro}.

\begin{proof}
Remember the change of probability of the Gaussian measure~: for all $ u \in \Rr $ and for $ X_\gamma \sim \Ns(0, \gamma^2) $
\begin{align*}
\frac{  e^{ u X_\gamma } }{ \Esp{ e^{ u X_\gamma } } } \bullet \Pp_{ X_\gamma } = \Pp_{X_\gamma + u \gamma^2 }
\end{align*}

Hence, for all $ u \in \Rr $
\begin{align*}
\frac{  \Esp{ e^{ u \He_\gamma(\Phi) } } }{ \Esp{ e^{ u X_\gamma } } } = \frac{ \Esp{ e^{ u X_\gamma } \Phi\prth{ \frac{X_\gamma}{\gamma^2} }   }  }{  \Esp{ e^{ u X_\gamma } } \Esp{ \Phi\prth{ \frac{X_\gamma}{\gamma^2} } }  } = \frac{ \Esp{  \frac{ e^{ u X_\gamma }  }{\Esp{ e^{ u X_\gamma } }} \Phi\prth{ \frac{X_\gamma}{\gamma^2} }   }  }{   \Esp{ \Phi\prth{ \frac{X_\gamma}{\gamma^2} } }  } = \frac{ \Esp{ \Phi\prth{ \frac{X_\gamma}{\gamma^2} + u } } }{ \Esp{ \Phi\prth{ \frac{X_\gamma}{\gamma^2} } } }
\end{align*}

Now, since $ \Phi $ is integrable, dominated convergence allows to exchange $ \lim_{\gamma \to +\infty } $ and expectation. One has moreover $ X_\gamma/\gamma^2 \eqlaw X_1 / \gamma $. This last quantity tends to $ 0 $ in distribution, hence, using the continuity of $ \Phi $, we have, locally uniformly in $ u \in \Rr $
\begin{align*}
\frac{  \Esp{ e^{ u \He_\gamma(\Phi) } } }{ \Esp{ e^{ u X_\gamma } } } \tendvers{\gamma}{ + \infty} \frac{ \Phi(u) }{\Phi(0) } = \Phi(u)
\end{align*}
as $ \Phi(0) = 1 $.
\end{proof}


\begin{remark}
It is enough to suppose that $ \Phi \in L^1(\Pp_{G/\gamma + u}) $ for all $ \gamma, u \in \Rr $.
\end{remark}


\begin{remark} What happens if the conditions of theorem \ref{Thm:InterprModG} are not fullfilled ? Can we still construct an equivalent of the distribution $ \Hs(\gamma, \Phi) $ ? This question is important since the moments conjecture described in \cite{KeatingSnaith} uses a function $ \Phi $ that can be written as $ \Phi = \Phi_A \Phi_M $ where $ \Phi_A $ (arithmetic factor) has no singularity on $ (-1, +\infty) $, but $ \Phi_M $ (matrix factor) is singular in the neighbourhood of $ -1 $ ; note that both functions are bounded on $ (K, +\infty) $ for all $ K > -1 $ by known asymptotics of the involved functions taken in a real variable (see \cite{KeatingSnaith}).

One way to proceed in this case is to use a cutoff function and a diagonal extraction procedure. For instance, one can use $ \Phi_K := \Phi \cdot \Un_{(K, +\infty)} $ and form the distribution $ \Hs(\gamma, \Phi_K ) $. The diagonal extraction in $ (\gamma, K) $ can ultimately define a sequence that converges in the mod-Gaussian-Laplace sense to $ \Phi $ on its total interval of definition. 
\end{remark}


\subsection{Mod-* convergence with infinitely divisible distributions} In order to understand the general mechanism at stake in section \ref{Section:LaplaceSetting}, we now generalize the construction \eqref{Def:VariableHIntro} to the setting of infinitely divisible distributions, that is, values of a L\'evy process at fixed time $ \gamma > 0 $. Since we consider Laplace transform, we restrict ourselves to positive random variables, i.e. one-dimensional marginals of subordinators. These distributions are characterised by a triplet $ ( \operatorname{k}, \operatorname{d}, \Pi ) $ satisfying $ \operatorname{k}, \operatorname{d} \geq 0 $ and $ \int_{\Rr_+} (1 \wedge x) \Pi(dx) < \infty $ (see e.g. \cite{SatoLevy}). If $ X_\gamma $ is a random variable having such a distribution, the L\'evy-Kintchine formula gives
\begin{align}\label{Eq:LevyKintchineSub}
\Esp{ e^{ - \theta X_\gamma} } = \exp\prth{ -\gamma \Lambda_X(\theta) } \ \ \mbox{with } \ \ \ \Lambda_X(\theta) = \operatorname{k} + \operatorname{d}  \theta + \int_0^{+ \infty } \prth{ 1 - e^{- \theta u } } \Pi(du)
\end{align}

We define mod-L\'evy-Laplace convergence by the following

\begin{definition} Let $ (Z_n)_n $ be a sequence of positive random variables such that $ \Esp{ e^{ -\theta Z_n } } < \infty $ for all $ \theta \in \Rr $ and let $ (\gamma_n)_n $ be a sequence of strictly positive real numbers. $ (Z_n)_n $ is said to converge in the mod-L\'evy-Laplace sense at speed $ (\gamma_n)_n $ if, locally uniformly in $ x \in \Rr_+ $
\begin{align*}
\frac{ \Esp{ e^{ - \Upsilon(x) Z_n} } }{ \Esp{ e^{ - \Upsilon(x) X_{\gamma_n} } } } \tendvers{n }{ + \infty } \Phi\prth{ x }
\end{align*}
where $ \Phi : \Rr_+ \to \Rr_+ $ is a continuous function satisfying $ \Phi(\Lambda'_X(0)) = 1 $, where $ X_{\gamma_n} $ is a random variable distributed according to \eqref{Eq:LevyKintchineSub} with the additional hypothesis $ \int_{ \Rr_+ } e^{ \alpha u } u \Pi(du) < \infty $ for all $ \alpha \geq 0 $ and (setting $ \inf \emptyset := \infty $)
\begin{align*}
\Upsilon(x) := \inf\ensemble{  y \in \Rr \ / \  \Lambda'_X(x) \geq y  }
\end{align*}
\end{definition}


The analogue to theorem \ref{Thm:Intro} is then given by

\begin{theorem}\label{Thm:InfDiv} Let $ \Phi : \Rr_+ \to \Rr_+ $ be a continuous positive integrable function on $ \Rr $ such that $ \Phi(\Lambda'_X(0)) = 1 $. Let $ \gamma > 0 $ and $ X_\gamma $ defined by \eqref{Eq:LevyKintchineSub}. Define the distribution of a random variable $ \Xe_\gamma(\Phi) $ by the change of probability 
\begin{align}\label{Def:VariableXinfdiv}
\Pp_{ \Xe_\gamma(\Phi) } := \frac{ \Phi\prth{ \frac{X_\gamma}{\gamma} }  }{ \Esp{ \Phi\prth{ \frac{X_\gamma}{\gamma} } } }  \bullet \Pp_{ X_\gamma }
\end{align}

Then, locally uniformly in $ x \in \Rr_+ $
\begin{align*}
\frac{  \Esp{ e^{ - \Upsilon(x) \Xe_\gamma(\Phi) } } }{ \Esp{ e^{ - \Upsilon(x) X_\gamma } } } \tendvers{ \gamma }{ + \infty }  \Phi(x)
\end{align*}
\end{theorem}


\begin{proof} Define for all $ y \in \Rr $
\begin{align*}
\Pp_{ X^{ (y) }_\gamma  } := \frac{ e^{ -y X_\gamma } }{ \Esp{ e^{ -y X_\gamma } } }  \bullet \Pp_{  X_\gamma   } 
\end{align*}

Since we have supposed $ \int_{ \Rr_+ } e^{ \alpha u } u \Pi(du) < \infty $ for all $ \alpha \geq 0 $, this random variable is well-defined. In particular, 
\begin{align*}
\Esp{ e^{-\theta X^{ (y) }_\gamma } } = e^{ -\gamma \, \prth{ \Lambda_X(\theta + y) - \Lambda_X(\theta) } } = \exp\prth{ \operatorname{d} \theta + \int_{\Rr_+} \prth{1 - e^{-\theta u } } e^{ - y u }  \Pi(du) }
\end{align*}

Hence, $ X^{ (y) }_\gamma $ is infinitely divisible of triplet $ (0, \operatorname{d}, \Pi^{ (y) } ) $ with
\begin{align*}
\Pi^{ (y) }(du) :=  e^{ - y u }  \Pi(du)
\end{align*}

For all $ y \in \Rr $, one has the duality relation
\begin{align}\label{Eq:DualiteModLevy}
\frac{ \Esp{ e^{ - y \Xe_{\gamma}(\Phi) } }  }{ \Esp{ e^{ - y X_{\gamma } } } } = \frac{ \Esp{ \Phi\prth{  \frac{ X_\gamma^{ (y) }  }{\gamma} } } }{ \Esp{ \Phi\prth{ \frac{ X_{\gamma } }{\gamma} } } }
\end{align}

Using the L\'evy-Kintchine formula \eqref{Eq:LevyKintchineSub}, one has
\begin{align*}
\frac{X_\gamma^{ (y) }}{\gamma } \cvlaw{\gamma}{+ \infty }  \operatorname{d}  + \int_0^{+ \infty } u \, \Pi^{ (y) }(du) =  \operatorname{d}  + \int_0^{+ \infty } u \, e^{ -u y }\, \Pi(du) = \Lambda'_X(y)
\end{align*}

Since we have supposed that $ \Lambda'_X(y) < \infty $ for all $ y \in \Rr_- $, and the case $ y \in \Rr_+ $ coming from the definition of a L\'evy measure, we have $ \Lambda'_X(y) < \infty $ for all $ y \in \Rr $.

Using the integrability and continuity of $ \Phi $, dominated convergence and this last convergence in distribution, we thus have
\begin{align*}
\frac{ \Esp{ e^{ - y \Xe_{\gamma}(\Phi) } }  }{ \Esp{ e^{ - y X_{\gamma } } } } \tendvers{\gamma}{ + \infty } \frac{ \Phi\prth{ \Lambda'_X(y) } }{ \Phi\prth{ \Lambda'_X(0) } }
\end{align*}

Last, remark that $ \Lambda'_X $ is decreasing on $ \Rr $ as a Laplace transform of the measure $ u \Unens{u \geq 0} \Pi(du) $ and its inverse bijection $ \Upsilon $ is well-defined. Setting $ y = \Upsilon(x) $, and using $ \Phi\prth{ \Lambda'_X(0) } = 1 $, one thus has the result.
\end{proof}


\begin{example} The Poisson distribution corresponds to $ (\operatorname{k}, \operatorname{d}, \Pi ) = (0, 0, \delta_1 )  $. Several examples of combinatorial random variables converging in the mod-Poisson-Laplace sense were analysed in \cite{HNNZ, Hwang1, KowalskiNikeghbali1}. Since $ \Lambda'_X(y) = e^{-y} $ and $ \Upsilon(x) = -\log x $, the construction of a random variable $ \Xe_\gamma(\Phi) $ uses the characteristic function $ x \mapsto \Esp{ x^Z } $ and one has, for $ P_\gamma $ Poisson-distributed of expectation $ \gamma > 0 $, locally uniformly in $ x \in \Rr_+ $
\begin{align*}
\frac{ \Esp{ x^{  \Xe_{\gamma}(\Phi) } }  }{ \Esp{ x^{ P_{\gamma } } } }  =  \frac{ \Esp{ \Phi\prth{  \frac{ P_{x \gamma}  }{\gamma} } } }{ \Esp{ \Phi\prth{ \frac{ P_{\gamma } }{\gamma} } } } \tendvers{\gamma}{ + \infty } \Phi(x)
\end{align*}
\end{example}


\begin{example} The Dickman distribution (see e.g. \cite{ArratiaBarbourTavare}) corresponds to $ (\operatorname{k}, \operatorname{d}, \Pi ) = (0, 0, \Pi_D )  $ with
\begin{align*}
\Pi_D(du) := \Unens{0 \leq u \leq 1} \frac{du}{u}
\end{align*}

Since $ \Pi_D $ has a compact support, the conditions of theorem \ref{Thm:InfDiv} are fullfilled. One has
\begin{align*}
\Lambda'_D(x) = \int_0^{+ \infty} u e^{-xu } \Pi_D(du) = \int_0^1 e^{-xu } du = \frac{1 - e^{-x} }{x}
\end{align*}
and one needs to invert this last bijection to get the corresponding $ \Upsilon(x) $.

The Dickman distribution occurs as the limiting distribution of random variables of the type $ \frac{1}{n} \sum_{k = 1}^n k Z_k $ with independent random variables $ (Z_k)_k $ such that $ \Prob{Z_k = 0} = 1/k $ ; for instance $ Z_k $ is Poisson, Bernoulli or Geometrically distributed of parameter $ 1/k $. It also arises in the framework of Poisson-Dirichlet point processes~: the fluctuations of the maximum length of the cycles of a random uniform permutation are Dickman-distributed (see e.g. \cite{ArratiaBarbourTavare}).
\end{example}


\begin{remark} An alternative method to find an analogue of theorem \ref{Thm:Intro} consists in using $ \Phi \circ \Upsilon $ in place of $ \Phi $ if $ \Upsilon $ is continuous and if $ \Prob{ \abs{\Upsilon( X^{(y)}_\gamma / \gamma)  } < \infty} = 1 $ for all $ \gamma > 0 $ and $ y \in \Rr $. 

In such a case, $ \Upsilon( X^{(y)}_\gamma / \gamma) \to y $ in distribution when $ \gamma \to + \infty $ and 
\begin{align*}
\frac{ \Esp{ e^{ - y \Xe_{\gamma}(\Phi \circ \Upsilon) } }  }{ \Esp{ e^{ - y X_{\gamma } } } } \tendvers{\gamma}{ + \infty } \frac{ \Phi \circ \Upsilon\prth{ \Lambda'_X(y) } }{ \Phi \circ \Upsilon\prth{ \Lambda'_X(0) } } = \Phi(y)
\end{align*}

Such a construction does not apply to the Poisson distribution since $ \Upsilon(x) = -\log x $ and $ \Prob{ \abs{ \Upsilon( P_\gamma/\gamma ) } = \infty} = \Prob{ P_\gamma = 0 } \neq 0 $.
\end{remark}

\section{Stein's method}\label{Chapitre:Stein}

\subsection{Reminder of Stein's method with zero-bias}

As explained in the introduction, there are several versions of Stein's method. The version we present here, in addition to illustrate the philosophy of the method, gives the easiest bound to have a normal approximation for sums of symmetric independent random variables, which is the example of interest in section~\ref{Chapitre:ExempleIID}.


\begin{definition}[Zero-bias transform, \cite{GoldsteinReinert, Ross}] Let $ X $ be a real random variable with zero expectation such that $ \Esp{ X^2 } < \infty $. The \textit{zero-bias transform} of $ X $ is the random variable $ X\zzer $ defined by the following identity 
\begin{align}\label{Def:ZeroBias}
\Esp{ Xf(X) } = \Esp{X^2} \Esp{ f'\prth{ X\zzer } } \ \ \  \forall f \in \He
\end{align}
$ \He $ being the space of continuously differentiable functions $ f $ such that $ \Esp{ \abs{ X f(X) } } < \infty $ and $ \Esp{ \abs{ f'(X) } } < \infty $.
\end{definition}

One useful property of such a transform is the following identity in distribution (see e.g. \cite{GoldsteinReinert, Ross})~: if $ Z_n = \sum_{k = 1}^n X_k $ with $ (X_k)_k $ a sequence of i.i.d. random variables, then
\begin{align}\label{Eq:ZeroBiaisEqLaw}
Z_n\zer \eqlaw Z_n + \prth{ X\zer_I - X_I } = \sum_{ k \neq I } X_k   +   X\zer_I 
\end{align}
where $ I \sim \Us\prth{ \intcrochet{1, n} } $ is a random variable independent of $ (X_k)_k $ and $ (X\zer_k)_k $, those two last sequences being independent, and $ (X\zer_k)_k $ being a sequence of i.i.d. random variables distributed according to the zero-bias distribution of $ X $.


This last property implies, for $ f \in \He $
\begin{align*}
\abs{ \Esp{ f'(Z_n) - \frac{Z_n}{\Esp{Z_n^2}} f(Z_n) } } =  \abs{ \Esp{ f'(Z_n) - f'\prth{ Z_n\zer }  }  } \leq  \norm{f''}_\infty \Esp{ \abs{ Z_n - Z_n\zer } } 
\end{align*}
that is, 
\begin{align}\label{Ineq:ZeroBoundsOnOperator}
\abs{ \Esp{ f'(Z_n) - \frac{Z_n}{\Esp{Z_n^2}} f(Z_n) } } \leq \norm{f''}_\infty \Esp{ \abs{ X\zer_I - X_I } } 
\end{align}

Hence, if the quantity $ \abs{ \Esp{ f'(Z_n) - Z_n f(Z_n)/\Esp{Z_n^2} } } $ is of interest to understand the behaviour of $ (Z_n)_n $, we have a useful bound.

$ $

Stein's method for Gaussian approximation consists in solving the following ``Stein's equation'' 
\begin{align}\label{Eq:SteinEquation}
\Le_{\Ns(0, 1) } f = h - \Esp{ h(G) }
\end{align}
with $ G \sim \Ns(0, 1) $, $ h $ a function such that $ \Esp{ \abs{h(G)} } < \infty $ and $ \Le_{\Ns(0, 1) } $ the operator defined for a differentiable function $f$ by  
\begin{align}\label{Def:SteinOperator}
\Le_{\Ns(0, \sigma^2) } f(x) :=  f'(x) - \frac{x}{\sigma^2} f(x)
\end{align}

This operator has the particularity that it characterises the Gaussian distribution in the following way (see \cite{SteinApprox}, p. 21)
\begin{align*}
X \sim \Ns(0,1) \ \ \  \Longleftrightarrow  \ \ \  \forall f \in \He, \ \ \Esp{ \Le_{\Ns(0, 1) } f (X) } = 0 
\end{align*}

A function of the form $  h - \Esp{ h(G) } $ is precisely in the image of the operator $ \Le_{\Ns(0, 1) } $. This allows to define the \textit{pseudo-inverse} $ \Le_{\Ns(0, 1) }\inv $  of the operator on such functions, that is, the inverse of the operator on its image, in addition to impose $ \Le_{\Ns(0, 1) }\inv (h -\Esp{ h(G) } ) $ to be the solution of \eqref{Eq:SteinEquation} that vanishes at infinity . Writing
\begin{align*}
h_G & :=  h - \Esp{ h(G) } \\
f   & := \Le_{\Ns(0, 1) }\inv  h_G
\end{align*}
we get
\begin{align*}
\abs{ \Esp{ h(W) } - \Esp{ h(G) } }  = \Esp{ h(W)  - \Esp{ h(G) } } = \Esp{ h_G(W) } = \Esp{ \Le_{\Ns(0, 1) }\Le_{\Ns(0, 1) }\inv h_G(W) }
\end{align*}

Injecting this equality in \eqref{Ineq:ZeroBoundsOnOperator}, and setting $D$ for the operator of differentiation, we get 
\begin{align*}
\abs{ \Esp{ h(W) } - \Esp{ h(G) } }  \leq  \norm{ D^2 \Le_{\Ns(0, 1) }\inv h_G }_\infty \Esp{ \abs{ W\zzer - G} }
\end{align*}


In the Gaussian case, setting $ f_G(x) := e^{-x^2/2} /\sqrt{2 \pi } $, we have (see \cite{SteinApprox} p. 15) 
\begin{align*}
\Le_{\Ns(0, 1) }\inv h_G(x) = \frac{ \Esp{ h_G(G) \Unens{G \leq x } } }{ f_G(x) }
\end{align*}

This allows prove the following inequalities (see \cite{SteinApprox} p. 25)
\begin{align*}
\norm{ \Le_{\Ns(0, 1) }\inv h_G }_\infty     & \leq   \sqrt{ \frac{ \pi }{2} } \norm{ h_G  }_\infty \\
\norm{ D \Le_{\Ns(0, 1) }\inv h_G }_\infty   & \leq   2 \norm{ h_G  }_\infty \\
\norm{ D^2 \Le_{\Ns(0, 1) }\inv h_G }_\infty & \leq   2 \norm{ h ' }_\infty 
\end{align*}
that lead to the Stein's bound 
\begin{align}\label{Ineq:SteinBoundsApprox}
\abs{ \Esp{ h(W) } - \Esp{ h(G) } }  \leq 2 \norm{ h'  }_\infty  \Esp{ \abs{ W\zzer - G} } 
\end{align}

This inequality allows to bound a particular distance between $ W $ and $ G $. Set
\begin{align*}
\He & := \ensemble{ h \in \Ce^1(\Rr) /  \norm{ h }_\infty \leq 1, \norm{ h' }_\infty \leq 1, \lim \!\!\! \hphantom{x}_{\pm \infty} h = 0 }  \\
d_{\He }(X, Y) & := \sup_{h \in \He  } \abs{ \Esp{h(X)} -\Esp{ h(Y) }  }
\end{align*}
then, \eqref{Ineq:SteinBoundsApprox} implies that 
\begin{align*}
d_{\He }(W, G) \leq 2 \Esp{ \abs{ W\zzer - G} } 
\end{align*}

\subsection{A Stein's operator for the penalised Gaussian distribution}

In order to apply Stein's method to a sequence of random variables converging in the mod-Gaussian sense with parameters $ ((\gamma_n)_n, \Phi) $, we first describe a characteristic operator $ \Le_{H_\gamma} $ of $ H_\gamma \sim \Hs(\Phi, \gamma) $.

\begin{theorem}[Characteristic operator for a penalisation of the Gaussian distribution]\label{Thm:OperateurCarac} Let $ \Phi $ satisfying the hypotheses of theorem \ref{Thm:InterprModG} and $  H_\gamma \sim \Hs(\Phi, \gamma) $. Suppose that $ \Phi > 0 $ on $ \Rr $. Set
\begin{align*}
\Psi(x)						   & := \log \Phi(x) \\
\kappa_\gamma(x) 			   & := \frac{x^2}{2 \gamma^2} - \Psi\prth{ \frac{x}{\gamma^2 } } \\
\rho_\gamma(x)   			   & := \frac{x}{\gamma^2} - \frac{1}{\gamma^2} \Psi'\prth{ \frac{x}{\gamma^2 } } = \kappa_\gamma'(x) \\
\widetilde{\He}_{\Phi, \gamma}  & := \ensemble{ h \in \Ce^1_m \ / \ \Esp{ \abs{ h'(H_\gamma) } } < \infty, \lim \!\!\! \hphantom{x}_{\pm \infty} h = 0 }
\end{align*}
where $ \Ce^1_m $ is the space of continuous and piecewise continuously differentiable functions $ f : \Rr \to \Rr $.

Suppose moreover that 

\begin{enumerate}

\item $  \kappa_\gamma(x) \tendvers{x}{ \pm \infty} + \infty $, 

\item $ \forall x \in \Rr $, $ \rho_\gamma'(x) \geq 0 $.

\end{enumerate}

Then, a characteristic operator of $ H_\gamma $ on $ \widetilde{\He}_{\Phi, \gamma} $ is $ \Le_{\Phi, \gamma} $ defined for all $ h \in \widetilde{\He}_{\Phi, \gamma} $ by
\begin{align}\label{Def:SteinOpMult}
\Le_{\Phi, \gamma} h(x) :=  h'(x) -  \prth{ \frac{x}{\gamma^2} - \frac{1}{\gamma^2} \Psi'\prth{ \frac{x}{\gamma^2} }  } h(x)
\end{align}
\end{theorem}


\begin{proof} 
By inverting the first order differential operator $ \Le_{\Phi, \gamma} = D - \rho_\gamma $, it is easily proven that for functions $ h_\gamma := h - \Esp{ h(H_\gamma) } $ with $ h \in \widetilde{\He}_{\Phi, \gamma} $, we have
\begin{align*}
\Le_{\Phi, \gamma} g =  h_\gamma \  \Longleftrightarrow \  g(x) = \Le_{\Phi, \gamma}\inv h_\gamma (x) =  \frac{1}{f\hgam(x)} \int_{- \infty }^x \!\!\! f\hgam(y) h_\gamma(y) dy
\end{align*}
$ f\hgam $ being given in \eqref{Eq:HLebesgueDensite}, i.e. $ f\hgam = e^{ - \kappa_\gamma } / c_\gamma $. Using the random variable $H_\gamma$, we can write
\begin{align}\label{Eq:SteinSolutionProba}
\Le_{\Phi, \gamma}\inv h_\gamma(x) =  \frac{  \Esp{ h_\gamma(H_\gamma) \Unens{ H_\gamma \leq x }  }   }{ f\hgam(x) } = - \frac{  \Esp{ h_\gamma(H_\gamma) \Unens{ H_\gamma \geq x }  }   }{ f\hgam(x) }
\end{align}
the last equality coming from the fact that $ \Esp{ h_\gamma(H_\gamma) } = 0 $.


Now, if a random variable $Y$ is such that for all $ h \in \widetilde{\He}_{\Phi, \gamma} $, $ \Esp{ \Le_{\Phi, \gamma} h(Y) } = 0 $, it is true in particular for the function 
\begin{align*}
h_{x, \gamma} := \Le_{\Phi, \gamma}\inv ( u_x - \Esp{ u_x(H_\gamma) } ) \ \ \ \mbox{ with } \ \ \ u_x : y \mapsto \Unens{ y \leq x }
\end{align*}

Let us prove that $  h_{x, \gamma} \in \widetilde{\He}_{\Phi, \gamma} $.  By \eqref{Eq:SteinSolutionProba}, we can write
\begin{align*}
h_{x, \gamma} (y) = \frac{  \Esp{ \prth{ \Unens{ H_\gamma \leq x } - \Prob{ H_\gamma \leq x } } \Unens{ H_\gamma \leq y }  }   }{ f\hgam(y) }  = \frac{\operatorname{cov}\prth{ \Unens{ H_\gamma \leq x }, \Unens{ H_\gamma \leq y } } }{f\hgam(y) } 
\end{align*}

Let $ x \geq 0 $. As $ H_\gamma \eqlaw -H_\gamma $, we only have to consider this case. 
\begin{align*}
\Prob{ H_\gamma \geq x } & = \int_x^{+ \infty} e^{ - \kappa_\gamma(u) } \frac{du}{c_\gamma} \ \ \ \mbox{with $ c_\gamma $ given in \eqref{Eq:HLebesgueDensite} }   \\
						 & \leq \int_x^{+ \infty}  \frac{\kappa'_\gamma(u) }{\kappa'_\gamma(x) } e^{ - \kappa_\gamma(u) } \frac{du}{c_\gamma} \ \ \ \mbox{since } \ \kappa''_\gamma(x) = \rho_\gamma'(x) \geq 0 \mbox{ by hypothesis } (1) \\
						 & =  \frac{ e^{ - \kappa_\gamma(x) } }{ c_\gamma \kappa'_\gamma(x) } 
\end{align*}

We can rewrite this inequality as 
\begin{align}\label{Ineq:Hgamma}
\Prob{ H_\gamma \geq x } \leq  \frac{ f\hgam(x) }{ \rho_\gamma(x) } 
\end{align}

Hence, setting
\begin{align*}
u_{x, \gamma}(y) :=  u_x(y)  - \Esp{ u_x(H_\gamma) } = \Unens{y \leq x} - \Prob{ H_\gamma \leq x }
\end{align*}
we have $ u_{x, \gamma}(y)  \leq 2 $ and
\begin{align*}
\abs{ h_{x, \gamma} (y) } = \abs{ \Le_{\Phi, \gamma}\inv u_{x, \gamma}(y) } \leq \frac{  \Esp{ \abs{ u_{x, \gamma}(H_\gamma) } \Unens{ H_\gamma \leq y }  }   }{ f\hgam(y) }  \leq 2 \frac{ \Prob{ H_\gamma \leq y } }{ f\hgam(y) } \leq \frac{2}{  \kappa'_\gamma(y) } \tendvers{y}{+ \infty } 0 
\end{align*}
as  $  \kappa_\gamma(x) \tendvers{x}{ \pm \infty} + \infty $ by hypothesis $ (2) $.

Using $ \Le_{\Phi, \gamma} = D - \rho_\gamma  $  and Stein's equation $ \Le_{\Phi, \gamma} h_{x, \gamma} =  u_{x, \gamma} $, we have 
\begin{align*}
D h_{x, \gamma} = \rho_\gamma h_{x, \gamma} + u_{x, \gamma} 
\end{align*}
which implies
\begin{align*}
\abs{ D h_{x, \gamma} (y) } & \leq \abs{ u_{x, \gamma}(y) } + \abs{ \rho_\gamma(y) h_{x, \gamma}(y) } \\
               & \leq 2 \prth{ 1 + \abs{ \rho_\gamma(y) \frac{ \Prob{ H_\gamma \geq y } }{  f\hgam(y) } } } \\
               & \leq 4 \quad \mbox{ by \eqref{Ineq:Hgamma} }
\end{align*}

This implies that $ \Esp{ \abs{h_{x, \gamma}'(H_\gamma) } } < \infty $, i.e. $  h_{x, \gamma} \in \widetilde{\He}_{\Phi, \gamma} $.  Then, for all $ x \in \Rr $,
\begin{align*}
0 = \Esp{ \Le_{\Phi, \gamma} h_{x, \gamma}(Y) } = \Esp{ u_x(Y) - \Esp{ u_x(H_\gamma) } }  = \Prob{ Y \leq x } - \Prob{ H_\gamma \leq x } 
\end{align*}
that is : $ Y \eqlaw H_\gamma $.


$ $

Reciprocally, let us prove that for all $ h \in \widetilde{\He}_{\Phi, \gamma} $, 
\begin{align*}
\Esp{ \Le_{\Phi, \gamma} h(H_\gamma) } = 0 
\end{align*}

As $ -\rho_\gamma = D\log(f\hgam) = f'\hgam / f\hgam $ we have, by the Fubini theorem 
\begin{align*}
\Esp{ h'(H_\gamma) }  & = \int_\Rr h'(u) f\hgam(u) du = \int_{-\infty}^0 h'(u) f\hgam(u) du +  \int_0^{+\infty} h'(u) f\hgam(u) du \\
					 & = \int_{-\infty}^0 h'(u) \int_{-\infty}^u f'\hgam(v) dv du  +  \int_0^{+\infty} h'(u) \int_v^{+\infty}-f'\hgam(v)dv du \\
					 & = - \int_{-\infty}^0 \!  \int_{-\infty}^0 \Unens{ v \leq u \leq 0 } \rho_\gamma(v) f\hgam(v) h'(u) dudv + \\
					 &  \ \ \ \ \ \ \ \ \ \ \ \  \int_0^{+\infty} \! \int_0^{+\infty} \Unens{ v \geq u \geq 0 } \rho_\gamma(v) f\hgam(v)h'(u) dudv \\
					 & = - \int_{-\infty}^0 \int_v^0 h'(u)du \, \rho_\gamma(v) f\hgam(v) dv + \int_0^{+\infty} \! \int_0^v h'(u)du \, \rho_\gamma(v) f\hgam(v) dv \\
					 & = \int_\Rr \prth{ h(v) - h(0) } \rho_\gamma(v) f\hgam(v) dv = \Esp{ h(H_\gamma)\rho_\gamma(H_\gamma) } - h(0) \Esp{ \rho_\gamma(H_\gamma) }
\end{align*}

Now, by symmetry of $ H_\gamma $ and by parity of $ \rho_\gamma $, we have $ \Esp{ \rho_\gamma(H_\gamma) } = 0 $. Hence, 
\begin{align*}
\Esp{ h'(H_\gamma) } = \Esp{ \rho_\gamma(H_\gamma)h(H_\gamma) } = \Esp{ \frac{1}{\gamma^2}\crochet{ H_\gamma -  \Psi'\prth{ \frac{ H_\gamma }{\gamma^2} } } h(H_\gamma)}
\end{align*}
\end{proof}


\begin{remark} We thus have 
\begin{align*}
Y \eqlaw H_\gamma  & \Longleftrightarrow    \Esp{ \Le_{\Phi, \gamma} h(Y) } = 0 \ \ \ \forall h \in \widetilde{\He}_{\Phi, \gamma}  \\
			      & \Longleftrightarrow       \Esp{ h'(Y) } =  \frac{1 }{\gamma^2 }\Esp{  \crochet{ Y  - \Psi'\prth{ \frac{Y}{\gamma^2} } } h(Y) }  \ \ \forall h \in \widetilde{\He}_{\Phi, \gamma}  
\end{align*}

The usual characterisation of the law $ \Ns(0, \gamma^2) $ can be recovered taking $ \Phi = 1 $ in the last formula, that is $ \Psi' = 0 $. If we think of $ \gamma $ as a parameter going to $ + \infty $, we have a small correction to the Gaussian distribution that takes the form of a small perturbation of the characteristic operator $ \Le_{\Ns(0, \gamma^2)} $, i.e. $ \Le_{\Phi, \gamma} = \Le_{\Ns(0, \gamma^2)} +  \Psi'(\cdot / \gamma^2 ) /\gamma^2  $.
\end{remark}

\subsection{Perturbation of the Gaussian operator and Edgeworth expansion}

For a random variable $ X $, denote by $ \phi_ X $ its Fourier transform 
\begin{align*}
\phi_X (u) := \Esp{ e^{i u X} }
\end{align*}

Note $ \phi_X (u) := \Fe f_X(u) $ if $ X $ admits a Lebesgue-density $ f_X $.

Let $ H_\gamma \sim \Hs(\Phi, \gamma) $, and suppose that $ \int_\Rr \abs{ f_\gamma(x) }^2 dx < \infty $. Denote $ \Phi_\gamma := \Phi( \cdot / \gamma) $ and let $ (\He_k)_k $ be the unnormalised Hermite polynomials defined by their Rodrigues form
\begin{align*}
\He_k(y) :=  e^{ y^2/2} \prth{ - \frac{d}{dy} }^k e^{-y^2/2}
\end{align*}

They satisfy $ \Esp{\He_k(G) \He_\ell(G) } = h_k \Unens{k = \ell} $ for $ G \sim \Ns(0, 1) $ and $ h_k > 0 $. Thus,
\begin{align*}
\phi_ {H_\gamma} (u) := \int_\Rr e^{iu x} f_\gamma(x) dx = \int_\Rr e^{iu x} \Phi \prth{ \frac{x}{\gamma^2} } e^{ -\frac{1}{2}  \prth{ \frac{x}{\gamma} }^2 } \frac{dx}{ c_\gamma } = \frac{\gamma}{c_\gamma } \int_\Rr \Phi_\gamma(y) e^{ i u \gamma y -\frac{y^2}{2}   } dy 
\end{align*}

Since $ \int_\Rr \abs{ \Phi_\gamma(x) }^2 e^{ -x^2 } dx =  \prth{ \frac{c_\gamma}{\gamma} }^2 \int_\Rr \abs{ f_\gamma(x) }^2 dx < \infty $, we can write
\begin{align*}
\Phi_\gamma  =  \sum_{k \geq 0} h_k\inv \Esp{ \He_k(G) \Phi_\gamma(G) } \He_k 
\end{align*}
which implies
\begin{align*}
f_\gamma (\gamma y)  =  \frac{1}{c_\gamma} \sum_{k \geq 0} h_k\inv \Esp{ \He_k(G) \Phi_\gamma(G) } \He_k (y) e^{-y^2/2} 
\end{align*}

A development of the form 
\begin{align*}
g_{\gamma}^{(k)} (y)  =  e^{-(y / \gamma)^2/2} \sum_{ \ell = 0}^k  a_\ell(\gamma)  \He_k \prth{ \frac{y}{\gamma} } 
\end{align*}
is said to be an \textit{Edgeworth expansion} of a random variable. For such a development, truncated at a certain order, there is no possibility to obtain a probability density due to the sign changes of the Hermite polynomials. Without truncation, the function $ y \mapsto g_{\gamma}^{(\infty)} (y)   e^{(y / \gamma)^2/2 } $ can still be positive and we get a probabilistic penalisation. From this point of view, mod-Gaussian convergence is a non-truncated Edgeworth expansion.


\subsection{Approximation by signed measures}

Using the Rodrigues form of the Hermite polynomials, we have 
\begin{align*}
f_\gamma (\gamma y) = c_\gamma\inv \sum_{k \geq 0} h_k\inv \Esp{ \He_k(G) \Phi_\gamma(G) } \prth{ -\frac{d}{dy} }^k e^{-y^2/2}
\end{align*}

Taking the Fourier transform and using the fact that $ \Fe (Df) (u) = -iu \, \Fe f(u) $, we get
\begin{align*}
\phi_{H_\gamma} (u)  & =  c_\gamma\inv\sum_{k \geq 0} h_k\inv \Esp{ \He_k(G) \Phi_\gamma(G) } (iu/ \gamma
)^k \Fe\prth{x \mapsto e^{-(x/\gamma)^2/2} } (u) \\
					 & = \gamma c_\gamma\inv \sqrt{2 \pi } e^{ - \gamma^2 u^2 / 2 } \sum_{k \geq 0} h_k\inv \Esp{ \He_k(G) \Phi_\gamma(G) } (iu/\gamma)^k
\end{align*}
the last development being convergent in $ L^2( e^{-x^2 } dx) $ since $ \Phi_\gamma $ belongs to this space. As $ \phi_{\gamma G}(u) =  e^{ - \gamma^2 u^2 / 2 } $, denoting by 
\begin{align*}
\widetilde{\Phi}_\gamma (u) := \gamma \sqrt{2 \pi } \sum_{k \geq 0} h_k\inv \Esp{ \He_k(G) \Phi_\gamma(G) } (iu/\gamma)^k
\end{align*}
we get 
\begin{align}\label{Eq:FourierLoiH}
\phi_{H_\gamma} (u) = \widetilde{\Phi}_\gamma (u) \phi_{\gamma G}(u)
\end{align}

In particular, if we know that locally uniformly in $u \in \Rr $, $ \Phi(u) = \lim_{\gamma \to \infty} \phi_{H_\gamma}(u) / \phi_{\gamma G}(u) $ exists, then, locally uniformly in $u \in \Rr $
\begin{align}\label{Eq:FourierLimitePhi}
\Phi(u) = \lim_{\gamma \to \infty} \widetilde{\Phi}_\gamma (u)
\end{align}
i.e. $ \widetilde{\Phi}_\gamma $ is an approximation of $ \Phi $.

The construction of a random variable whose distribution satisfies \eqref{Eq:FourierLoiH} and \eqref{Eq:FourierLimitePhi} is not always possible. In the flavour of \cite{BarbourAl}, one can be interested in a signed measure approximation of sequences converging in the mod-Gaussian sense. A case of interest is to suppose that $ \widetilde{\Phi}_\gamma $ can be approximated by
\begin{align*}
\Phi^\sharp_\gamma(u) := e^{ P(i\theta) } 
\end{align*}
$ P $ being a polynomial satisfying the symmetry condition of theorem \ref{Thm:InterprModG}, i.e. $ P(i\theta) = P(\theta) $ for all $ \theta \in \Rr $, and $ P(0) = 0 $. Let  $ \mu_\gamma := \Fe\inv\prth{ \Phi^\sharp_\gamma \, \phi_{\gamma G}  } $ be the signed measure obtained by inverting equation \eqref{Eq:FourierLoiH}. Then, a Stein operator of $ \mu_\gamma $ can be defined by
\begin{align*}
\Le_{\Phi^\sharp, \gamma} := \Le_{ \Ns(0, \gamma^2) } - P'\prth{ -\frac{d}{dx} }
\end{align*}

Indeed, suppose $ \Fe \mu_\gamma \in L^1 $ with density $ f \in \Ce^1 $. Such an operator satisfies $ \int_\Rr \Le_{\Phi^\sharp, \gamma} g \cdot f = 0 $ for all functions $ g $ of class $ \Ce^\infty \cap L^1 $ that vanish at $ \pm \infty $. By integration by parts, this amounts to the following equation which is the analogue of Stein's equation for positive measures
\begin{align*}
\Le^*_{\Phi^\sharp, \gamma} f (x) := x f(x) + \gamma^2 f'(x) - P'\prth{ \frac{d}{dx} } f(x) = 0
\end{align*}

Taking the Fourier transform, setting $ \hat{f} := \Fe f $ and using $ \Fe\prth{ x \mapsto xf(x) }(\xi) = -i \frac{d}{d\xi} \hat{f} (\xi) $ and $ \Fe\prth{ f' }(\xi) = (-i\xi) \hat{f} (\xi) $, we get
\begin{align*}
-i \frac{d}{d\xi} \hat{f}(\xi) + \prth{ \gamma^2 (-i \xi) - P'(i\xi) } \hat{f}(\xi) = 0
\end{align*}

The integration of this equation gives 
\begin{align*}
\hat{f}(\xi) = \hat{f}(0) \exp\prth{ -\gamma^2 \frac{\xi^2}{2} + P(  i \xi ) - P(0) } = \hat{f}(0)  \exp\prth{ -\gamma^2 \frac{\xi^2}{2} + P( \xi )  }
\end{align*}

By Fourier inversion, 
\begin{align*}
f(x) = \hat{f}(0) \int_\Rr \exp\prth{ -i \xi x  -\gamma^2 \frac{\xi^2}{2} + P(  \xi ) } \frac{d\xi}{2\pi}
\end{align*}
which is (proportional to) the density of the measure $ \mu_\gamma $. 

$ $

In this setting, Stein's method with an operator such as $ \Le_{\Phi^\sharp, \gamma} $ allows to approximate $ \Pp_{H_\gamma} $ with $ \mu_\gamma $. Such a procedure is to relate to \cite{BarbourAl} where Kolmogorov approximations in the Poisson setting were found with respect to a signed measure, and to \cite{BarbourAsymptotic} where this type of correction to $ \Le_{ \Ns(0, \gamma^2) } $ is discussed in details (see also \cite{Rotar}). 

In \cite{BarbourAsymptotic}, a perturbation of $ \Le_{ \Ns(0, 1) } $ is done with a polynomial in the operator of differentiation, and this polynomial is the truncation of the cumulant generating series. As one can see on example \ref{Ex:SommeIID}, in the case of the sum of i.i.d. random variables, the limiting function $ \Phi $ is obtained by a suitable renormalisation of the cumulant function.

\begin{example} In the case of example \ref{Ex:SommeIID}, as an approximation of  $ \widetilde{\Phi}_\gamma $ is $ \Phi^\sharp_\gamma = \Phi ( \cdot / \gamma^2 ) : x \mapsto e^{-C x^4/(4 \gamma^8 ) } $, the condition is fullfilled and one can have an Edgeworth expansion by using a suitable truncation of the Taylor development of $ e^{-C x^4/4 } $ in addition to a signed measure approximation of density $ x \mapsto \int_\Rr \exp\prth{ -i \xi x  -\gamma^2 \frac{\xi^2}{2} -C  \xi^4/(4 \gamma^8  ) } \frac{d\xi}{2\pi} $.
\end{example}

\section{The sum of i.i.d. symmetric random variables}\label{Chapitre:ExempleIID}

We develop the important example of the sum of i.i.d. random variables. In order to agree with theorem \ref{Thm:InterprModG}, we only consider the symmetric case.


\subsection{A mod-Gaussian approximation theorem}

\begin{theorem}[Mod-Gaussian bounds for the sum of i.i.d. random variables]\label{Thm:MASSmodG} Let $ (X_k)_k $ be a sequence of i.i.d. symmetric random variables having the same law as $X$ such that  $ \Esp{X} = 0 $, $ \Esp{X^2} = 1 $ and $ \Esp{X^4} < 3 $. Set 
\begin{align*}
Z_n        & :=  \frac{1}{ n^{1/4} } \sum_{k = 1}^n X_k  \\
\gamma_n   & :=  n^{1/4} \\
C		   & :=  \frac{3 - \Esp{X^4} }{6} > 0 \\
\Phi_C(x)  & :=  e^{ - C \frac{x^4}{4 } } \\
c_1        & :=  \sqrt{ 2 \pi  } \, \Esp{ \Phi_C(G) } \  \mbox{with } \ G \sim \Ns(0,1)
\end{align*}

Define moreover
\begin{align*}
\He^1_\Phi & := \ensemble{ h \in \Ce^1(\Rr) /  \norm{ h }_\infty \leq 1, \norm{ h' }_\infty \leq 1, \lim_{\pm \infty} h = 0, \int_\Rr \abs{h'} < \infty }  \\
d_{\He^1_\Phi}(X, Y) & := \sup_{h \in \He^1_\Phi } \abs{ \Esp{h(X)} -\Esp{ h(Y) }  }
\end{align*}

Let $ H_n \sim \Hs(\Phi_C, \gamma_n) $ and $ h \in \He^1_\Phi $. Set also $ \sigma_{1, 3} := \Esp{ \abs{X} } \vee \Esp{ \frac{ \abs{X}^3  }{2} } $. Then, 
\begin{align}\label{Dtv:BorneSteinSansEspace}
\begin{aligned}
\ \  \abs{ \Esp{h(Z_n)} -\Esp{ h(H_n) } }  \leq   \frac{4 \sqrt{ 2(1-C) }}{\gamma_n} \norm{h'}_\infty  +  \frac{ 4 }{\gamma_n^2 } \norm{ h }_\infty  \prth{   C c_1 \sigma_{1, 3}    +  \frac{ 1 }{\gamma_n^2 } }
\end{aligned}
\end{align}

In particular, for $n$ large enough
\begin{align}\label{Dtv:BorneSteinEspaceH}
\begin{aligned}
d_{\He^1_\Phi}(Z_n, H_n)  & \leq  \frac{4 \sqrt{2(1 - C)} }{\gamma_n} + O\prth{ \frac{1}{ \gamma_n^2 } } \\
d_{\He^1_\Phi} \prth{ \frac{Z_n}{\gamma_n} , \frac{H_n}{\gamma_n}} & \leq  \frac{4 \sqrt{2(1 - C)} }{\gamma_n^2 } + O\prth{ \frac{1}{ \gamma_n^4 } }
\end{aligned}
\end{align}

\end{theorem}

$ $

\begin{proof}

We start with the usual Stein's argument~: for  $ h \in \He^1_\Phi $ 
\begin{align*}
h\hn(x) := h(x) - \Esp{ h(H_n) } = \Le\hn\Le_{H_n}\inv h\hn(x)
\end{align*}

Setting $ g := \Le_{H_n}\inv h\hn $, we get 
\begin{align*}
\Esp{h(Z_n)} - \Esp{ h(H_n) } = \Esp{ \Le_{H_n}g(Z_n) } = \Esp{ g'(Z_n) - \rho_{\gamma_n} (Z_n) g(Z_n) }
\end{align*}

We recall that
\begin{align*}
\rho_{\gamma_n}(x) := \frac{1}{\gamma_n^2 } \prth{ x - (\log \Phi_C)' \prth{ \frac{x}{\gamma_n^2} } } =  \frac{x}{\gamma_n^2} + C \frac{x^3}{\gamma_n^8}
\end{align*}

For the usual Stein's operator given by $ \Le_{\Ns(0, \gamma_n^2) } f(x) :=  f'(x) - (x/\gamma_n^2) f(x) $, we have 
\begin{align*}
\Esp{\Le_{\Ns(0, \gamma_n^2) } g(Z_n) } = \Esp{  g'(Z_n) -  \frac{Z_n}{\gamma_n^2} g(Z_n) } =  \Esp{ g'(Z_n) - g'\prth{ Z_n\zer } }
\end{align*}
where $ Z_n\zer $ is the \textit{zero-bias transformation} of $ Z_n $. By \eqref{Eq:ZeroBiaisEqLaw}, we have
\begin{align*}
Z_n\zer \eqlaw Z_n + \frac{ \gamma_n }{ \sqrt{n} } \prth{ X\zer_I - X_I } = Z_n + \frac{ 1 }{ \gamma_n } \prth{ X\zer_I - X_I }
\end{align*}
where $ I \sim \Us\prth{ \intcrochet{1, n} } $ is a random variable independent of $ (X_k)_k $ and $ (X\zer_k)_k $, those two last sequences being independent, and $ (X\zer_k)_k $ is a sequence of i.i.d. random variables distributed according to the zero-bias distribution of $ X $.


In particular, 
\begin{align*}
\abs{ \Esp{\Le_{\Ns(0, \gamma_n^2) } g(Z_n) } }  & \leq   \Esp{ \abs{g'(Z_n) - g'\prth{ Z_n\zer }  }  } \\
						  & \leq \norm{g''}_\infty \Esp{ \abs{ Z_n - Z_n\zer } } = \frac{1}{\gamma_n} \norm{g''}_\infty \Esp{ \abs{ X\zer_I - X_I } } 
\end{align*}

Moreover, by independence and the i.i.d. property,
\begin{align*}
\Esp{ \abs{ X\zer_I - X_I } } = \frac{1}{n} \sum_{k = 1}^n \Esp{ \abs{ X\zer_k - X_k } } = \Esp{ \abs{ X\zer - X } }
\end{align*}


For the perturbative operator, we get
\begin{align*}
 \Esp{(\Le_{H_n} - \Le_{\Ns(0, \gamma_n^2) }) g(Z_n) }  =   -\frac{C}{\gamma_n^8} \Esp{ Z_n^3 g(Z_n)  } =: -\frac{C}{\gamma_n^8} \Esp{ Z_n \tilde{g}(Z_n)  }
\end{align*}
with $ \tilde{g}(x) := x^2 g(x) $. Using the same technique, we write
\begin{align*}
\Esp{ Z_n \tilde{g}(Z_n)  } = \gamma_n^2 \Esp{ \tilde{g}'\prth{ Z_n\zer } } = \gamma_n^2 \Esp{ 2 Z_n\zer g\prth{ Z_n\zer } + \prth{ Z_n\zer }^2 g'\prth{ Z_n\zer } }
\end{align*}

Hence, 
\begin{align*}
\abs{ \Esp{(\Le_{H_n} - \Le_{\Ns(0, \gamma_n^2) }) g(Z_n) } }      & \leq     \frac{C}{\gamma_n^8} \prth{  2 \gamma_n^2 \Esp{ \abs{ Z_n\zer } }  \norm{g}_\infty  +  \gamma_n^2   \Esp{ \abs{ Z_n\zer }^2 }  \norm{g'}_\infty    } \\
          & \leq   \frac{C}{ \gamma_n^6 }  \Bigg( 2 \norm{g}_\infty \frac{ \gamma_n }{ \sqrt{n} } \crochet{ (n-1) \Esp{ \abs{X} } + \Esp{ \abs{ X\zzer } } }   \\  
          &  \quad\quad  + \norm{g'}_\infty \prth{ \frac{ \gamma_n }{ \sqrt{n} } }^2  \crochet{ (n-1) \Esp{X^2}  + \Esp{ \abs{ X\zzer }^2 } } \Bigg) \\
          & \leq   \frac{C}{ \gamma_n^6 } \prth{ 2 \norm{g}_\infty \frac{n}{\gamma_n } \Esp{ \abs{X} } \vee \Esp{ \abs{X\zzer } }  + \frac{n}{\gamma_n^2 } \norm{g'}_\infty  \Esp{ \abs{X}^2 } \vee \Esp{ \abs{ X\zzer }^2 }   }
\end{align*}

Using \eqref{Def:ZeroBias} for well-choosen functions, we can prove that  
\begin{align*}
\Esp{ \abs{  X\zzer }  }       & =   \frac{ \Esp{ \abs{X}^3 }  }{2 \Esp{ X^2 } } \\
\Esp{ \prth{ X\zzer }^2   }    & =   \frac{ \Esp{ X^4 }  }{3 \Esp{ X^2 } }
\end{align*}

As $ \gamma_n^4 = n $, $ \Esp{X^2} = 1 $, and $ \Esp{ X^4 } = 3(1 - 2 C) < 3 $, using the independence of $ X $ and $ X\zzer $, we get
\begin{align*}
\Esp{ \abs{ X - X\zzer } } \leq \sqrt{\Esp{ \abs{ X - X\zzer }^2 } } = \sqrt{\Esp{ X^2 + (X\zzer)^2 } } = \sqrt{ 2(1-C) } 
\end{align*}
and finally
\begin{align*}
\abs{ \Esp{\Le_{H_n}g(Z_n) } }  \leq     \frac{\sqrt{ 2(1-C) }}{\gamma_n} \norm{g''}_\infty  +    \frac{2 C}{\gamma_n^3 } \norm{g}_\infty  \Esp{ \abs{X} } \vee \Esp{ \frac{ \abs{X}^3  }{2} }   +  \frac{ C}{\gamma_n^4 } \norm{g'}_\infty
\end{align*}

$ $

The last part of Stein's method consists in using the estimates for $ \norm{D^k g}_\infty = \norm{D^k \Le_{H_n}\inv h_{ H_n } }_\infty $ with $ k \in \ensemble{0, 1, 2} $, which is done in lemma \ref{Lemma:BorneSteinOperateur} in the appendix. Substituting the inequalities \eqref{Ineq:EstimeesBornees} to \eqref{Ineq:EstimeesAC} in the last inequality and using $ \sigma_{1, 3} := \Esp{ \abs{X} } \vee \Esp{ \frac{ \abs{X}^3  }{2} } $, we get
\begin{align}\label{Ineq:BornePresqueFinale}
\abs{ \Esp{h(Z_n)} -\Esp{ h(H_n) } }  \leq   \frac{4 \sqrt{ 2(1-C) }}{\gamma_n} \norm{h'}_\infty  +  \frac{ 2 }{\gamma_n^2 } \norm{ h_{ H_n } }_\infty  \prth{   C c_1 \sigma_{1, 3}  +  \frac{ 1 }{\gamma_n^2 } }
\end{align}

Last, using $ \norm{ h_{ H_n } }_\infty \leq 2 \norm{ h }_\infty $, we get the desired bound.
\end{proof}


\begin{remark}\label{Rk:BonScaling} With the suitable rescaling of $ Z_n $ and $ H_n $, we get
\begin{align*}
\abs{ \Esp{ h\prth{ \frac{Z_n}{\gamma_n} }  } -\Esp{ h\prth{ \frac{H_n}{\gamma_n} } } }  \leq  \frac{1}{\gamma_n^2 } \prth{ 4 \sqrt{ 2(1-C) } \norm{h'}_\infty  + \norm{ h  }_\infty O(1)  }     = O\prth{ \frac{1}{ \sqrt{n} } }
\end{align*}
which corresponds to the classical Berry-Ess\'een bound for the CLT, with no improvement. A natural question is thus to ask whether this is optimal or not. This bound will be improved in theorem \ref{Thm:MASSmodG2}, thus showing that we indeed gain a considerable factor by approximating with this new random variable.
\end{remark}

\subsection{A Kolmogorov approximation} 

Following the steps of Stein (\cite{SteinApprox} p. 36), we have the

\begin{corollary}[Kolmogorov bounds] Let $ (X_k)_k $ satisfying the hypothesis of theorem \ref{Thm:MASSmodG}. Then
\begin{align}\label{Ineq:KolmogorovBound}
\dkol(Z_n, H_n) := \sup_{x \in \Rr} \abs{ \Prob{Z_n \leq x} - \Prob{H_n \leq x} }  \leq  \frac{ 4  (1-C)^{1/4} c_1^{-1/2} }{\gamma_n}  + O\prth{ \frac{1}{\gamma_n^2} }
\end{align}

\end{corollary}

$ $

\begin{proof}

Set 
\begin{align*}
h_{x, \delta }(y) := \Unens{ y \leq x } + \prth{ 1 - \frac{y-x}{\delta} } \Unens{x \leq y \leq x + \delta }
\end{align*}

For all $x$, we have $ h_{x - \delta, \delta }(y)  \leq \Unens{ y \leq x } \leq h_{x, \delta }(y) \leq \Unens{ y \leq x + \delta } $, which implies that 
\begin{align}\label{Ineq:Dtv}
\Esp{ h_{x - \delta, \delta }(Z_n) }  \leq \Prob{ Z_n \leq x } \leq  \Esp{ h_{x, \delta }(Z_n) } 
\end{align}

Moreover, we have
\begin{align*}
\norm{ h_{x, \delta } - \Esp{h_{x, \delta }(H_n) } }_\infty   & \leq  1   \\
\norm{ h'_{x, \delta }  }_\infty  							& \leq   \frac{1}{\delta}
\end{align*}

Using these inequalities, inequality \eqref{Ineq:BornePresqueFinale}, and $ \alpha :=  \Esp{ \abs{X} } \vee \Esp{ \frac{ \abs{X}^3  }{2} } \ c_1 C   +  1 $, we get
\begin{align*}
\abs{ \Esp{h_{x, \delta }(Z_n)} - \Esp{ h_{x, \delta }(H_n) } }  \leq   \frac{4 \sqrt{ 2(1-C) }}{\gamma_n} \frac{1}{\delta}  +  \frac{ 2 \alpha }{\gamma_n^2 } 
\end{align*}

By \eqref{Ineq:Dtv}, we have
\begin{align*}
\Prob{ Z_n \leq x }   & \leq   \Esp{ h_{x, \delta }(Z_n) } \\
					 & \leq    \Esp{ h_{x, \delta }(H_n) } + \frac{4 \sqrt{ 2(1-C) }}{\gamma_n \delta} +  \frac{ 2 \alpha }{\gamma_n^2 } \\
					 & =     \Prob{ H_n \leq x } + \Esp{ \prth{ 1 - \frac{H_n - x}{ \delta } } \Unens{ 0 \leq \frac{H_n - x}{ \delta } \leq 1 }  }  + \frac{4 \sqrt{ 2(1-C) }}{\gamma_n \delta} +  \frac{ 2 \alpha }{\gamma_n^2 } \\
					 & \leq  \Prob{ H_n \leq x } +  \Prob{ 0 \leq H_n - x \leq \delta }  + \frac{4 \sqrt{ 2(1-C) }}{\gamma_n \delta} +  \frac{ 2 \alpha }{\gamma_n^2 } \\
					 & \leq  \Prob{  H_n \leq x } +  \frac{\delta}{c_{\gamma_n}}  + \frac{4 \sqrt{ 2(1-C) }}{\gamma_n \delta} +  \frac{ 2 \alpha }{\gamma_n^2 } 
\end{align*}

This last inequality comes from the fact that
\begin{align*}
\Prob{ 0 \leq H_\gamma - x \leq \delta }  =  \int_0^\delta e^{ - P_{q_\gamma}(y-x) } \frac{dy}{c_\gamma} \leq \frac{ \delta }{c_\gamma} \sup_{y \in \Rr_+ } \ensemble{ e^{ - P_{q_\gamma}(y-x) } } = \frac{ \delta }{c_\gamma}
\end{align*}

Optimising in $ \delta $ the LHS of the former inequality gives 
\begin{align*}
\delta = \sqrt{  c_{\gamma_n} \frac{4 \sqrt{ 2(1-C) }}{\gamma_n } }
\end{align*}
and the optimal value 
\begin{align*}
\Prob{ Z_n \leq x } - \Prob{  H_n \leq x }  \leq  2 \sqrt{   \frac{4 \sqrt{ 2(1-C) }}{ c_{\gamma_n} \gamma_n } }   +  \frac{ 2 \alpha }{\gamma_n^2 } 
\end{align*}

As by \eqref{Eq:HLebesgueDensite}, $ c_\gamma = \gamma \sqrt{2 \pi } \Esp{ e^{ -C \frac{G^4}{4 \gamma^4} } } \geq \gamma \sqrt{ 2 \pi } \Esp{ e^{ -C \frac{G^4}{4} } } = \gamma c_1  $ for $ \gamma \geq 1 $, we finally have
\begin{align*}
\Prob{ Z_n \leq x } - \Prob{  H_n \leq x }  \leq   \frac{ 4 ( (1-C)/c_1^2 )^{1/4} }{\gamma_n}   +  \frac{ 2 \alpha }{\gamma_n^2 } 
\end{align*}

The corresponding lower bound follows from the same manipulations using the lower bound in \eqref{Ineq:Dtv}.
\end{proof}


\subsection{Beyond the classical Berry-Ess\'een speed of convergence}

As noticed in remark \ref{Rk:BonScaling}, the bound \eqref{Ineq:BornePresqueFinale} is not optimal since a suitable rescaling gives the same speed of convergence as the usual CLT. Using the knowledge of mod-Gaussian convergence of the sequence, we can improve on this bound in the following way~:

\begin{theorem}[Mod-Gaussian bounds for the sum of i.i.d. random variables]\label{Thm:MASSmodG2} Let $ (X_k)_k $ be a sequence satisfying the hypotheses of theorem \ref{Thm:MASSmodG}. Define
\begin{align*}
 & \He^2_\Phi:= \ensemble{ h \in \Ce^2(\Rr) /  \norm{ h }_\infty \leq 1, \norm{ h' }_\infty \leq 1, \norm{ h'' }_\infty \leq 1, \lim_{\pm \infty} h = \lim_{\pm \infty} h' = 0, \int_\Rr \abs{h'} < \infty , \int_\Rr \abs{h''} < \infty}  \\
 & d_{\He^2_\Phi}(X, Y) := \sup_{h \in \He^2_\Phi } \abs{ \Esp{h(X)} -\Esp{ h(Y) }  }
\end{align*}

Then, for all $ h \in \He^2_\Phi $ and for all $ n \geq 1 $ (with $ \gamma_n := n^{\frac{1}{4}} $), we have the following bound
\begin{align}\label{Dtv:BorneSteinImprovedSansEspace}
\begin{aligned}
\abs{ \Esp{h(Z_n)} -\Esp{ h(H_n) } }  \leq   \prth{3 + 2C + \frac{35 C}{ \gamma_n^4 }} \frac{2 - 3C }{\gamma_n^2} \norm{h''}_\infty +  \frac{66 \, C}{\gamma_n^3 } \norm{h'}_\infty 
\end{aligned}
\end{align}

In particular,
\begin{align}\label{Dtv:BorneSteinImprovedEspaceH}
\begin{aligned}
d_{\He^2_\Phi}(Z_n, H_n)  & \leq   \frac{(3 + 2C)(2 - 3C ) }{\gamma_n^2} + O\prth{ \frac{1}{ \gamma_n^3} } \\
d_{\He^2_\Phi}\prth{ \frac{Z_n}{\gamma_n}, \frac{H_n}{\gamma_n} } & \leq  \frac{(3 + 2C)(2 - 3C ) }{\gamma_n^4} + O\prth{ \frac{1}{ \gamma_n^8} } 
\end{aligned} 
\end{align}
\end{theorem}

$ $

Before proving the theorem, observe the difference between \eqref{Dtv:BorneSteinSansEspace} and \eqref{Dtv:BorneSteinImprovedSansEspace}~: we need one more degree of regularity in the functions that we use. This could lead to a problem, as one can check that $ D^3 \Le_{H_n}\inv $ is not bounded on the unit sphere of $ L^\infty $. The trick is to consider a surrogate function, namely to invert $ \Le_{H_n} $ on the space orthogonal to $ (1, X) $ where $ X(x) := x $. By restricting $ D^3 \Le_{H_n}\inv $ to this space, we get a continuous operator for the $ L^\infty $ topology. But the final Stein bound will use, in the same lines as the proof of \eqref{Dtv:BorneSteinSansEspace}, the norm of the single derivative of such a projected function, and using the classical estimate that relates it with the norm of the initial function will result in an additional linear term, a priori not bounded. The trick will be to remove the projection on $ X $ by utilizing one more derivative, loosing in the procedure a factor $ \gamma_n $, see the second estimate of \eqref{Ineq:EstimeesAC}. By rescaling the random variables, though, we have improved the classical Stein bound (see remark \ref{Rk:ImprovedSteinBound}).


\begin{proof}

For $ h \in \He^2_\Phi $, define the following $ L^2(\Pp_{H_n}) $-projection on $ (1, X) $
\begin{align*}
\widehat{h}\hn(x) :=  h(x) - \Esp{ h(H_n) } - x \Esp{ h'(H_n) }
\end{align*}

As $ \Esp{H_n} = \Esp{Z_n} = 0 $, one can write
\begin{align*}
\Esp{h(Z_n)} - \Esp{ h(H_n) } =  \Esp{ \widehat{h}\hn(Z_n) }
\end{align*}


As $ \Esp{ \widehat{h}\hn (H_n) }= 0 $, one can invert $ \Le_{H_n} $ on $ \widehat{h}\hn $ and define 
\begin{align*}
g := \Le_{H_n}\inv \widehat{h}\hn
\end{align*}

We thus have
\begin{align*}
\Esp{h(Z_n)} - \Esp{ h(H_n) } = \Esp{ \Le_{H_n}g(Z_n) } = \Esp{ g'(Z_n) - \rho_{\gamma_n} (Z_n) g(Z_n) }
\end{align*}

We treat the case of the Stein's operator $ \Le_{\Ns(0, \gamma_n^2) } $ like before, with
\begin{align*}
\Esp{\Le_{\Ns(0, \gamma_n^2) } g(Z_n) } = \Esp{  g'(Z_n) -  \frac{Z_n}{\gamma_n^2} g(Z_n) } =  \Esp{ g'(Z_n) - g'\prth{ Z_n\zer } }
\end{align*}

Setting 
\begin{align*}
\Delta Z_n := Z_n\zer - Z_n = \frac{X_I\zer - X_I}{\gamma_n}
\end{align*}
we have
\begin{align*}
g'(Z_n\zer) -  g'(Z_n)  - \Delta Z_n g''(Z_n)  & =  \int_0^1 g''\prth{ Z_n + w \Delta Z_n } \Delta Z_n \, dw - \Delta Z_n g''(Z_n) \\
					& = \Delta Z_n \int_0^1 \int_0^1 g'''\prth{ Z_n + w u \Delta Z_n  }  w \Delta Z_n du  \, dw \\
					& =  \frac{(\Delta Z_n)^2}{2} \int_0^1 \int_0^1 g'''\prth{ Z_n + u \sqrt{v} \Delta Z_n  } du \, dv 
\end{align*}


By independence of the terms and since $ \Esp{X_I\zer} = 0 $ as seen in \eqref{Def:ZeroBias} taking $ f(x) = x^2 $, we have
\begin{align*}
\Esp{ \Delta Z_n g''(Z_n) } & = \frac{1}{\gamma_n} \Esp{ (X_I\zer - X_I) \, g'' \prth{ \frac{1}{\gamma_n} \sum_{k = 1}^n X_k  } } \\
						   & = \frac{1}{\gamma_n} \Esp{ X_I\zer } \Esp{ g'' \prth{ \sum_{k = 1}^n X_k  }  }  - \frac{1}{\gamma_n} \Esp{ X_I  g'' \prth{ \frac{1}{\gamma_n} \sum_{k = 1}^n X_k   } } \\ 
						  & = -\frac{1}{\gamma_n}  \Esp{ \frac{1}{n} \sum_{k = 1}^n X_k  g'' \prth{ \frac{1}{\gamma_n} \sum_{k = 1}^n X_k } }
\end{align*}
by integrating on $ I $ which is independent of $ (X_k)_k $. As $ n = \gamma_n^4 $, we hence have
\begin{align*}
\Esp{ \Delta Z_n g''(Z_n) }  & = -\frac{1}{\gamma_n^4}  \Esp{ S_n g'' \prth{ S_n   } } = -\frac{1}{\gamma_n^4} \Esp{S_n^2} \Esp{ g'''(Z_n\zer) } = -\frac{1}{\gamma_n^2} \Esp{ g'''(Z_n\zer) }
\end{align*}

$ $

Let $ U, V \sim \Us\prth{[0,1]} $ be two independent random variables. Then,
\begin{align*}
\abs{ \Esp{\Le_{\Ns(0, \gamma_n^2) } g(Z_n) } } & = \frac{1}{\gamma_n^2 } \abs{ \Esp{ \frac{(X_I - X\zer)^2}{2} g'''\prth{ Z_n + 2 V \sqrt{U} \Delta Z_n } -  g'''\prth{ Z_n\zer } } } \\
					& \leq  \frac{1}{\gamma_n^2 } \norm{ g''' }_\infty \prth{ \Esp{ \frac{(X_I - X_I\zer)^2}{2} } + 1 }
\end{align*}

By independence of $ X $ and $ X\zzer $ and by the i.i.d. property
\begin{align*}
\Esp{ ( X_I - X_I\zer )^2 } = \frac{1}{n} \sum_{k = 1}^n \Esp{ ( X_k - X\zer_k )^2 } = \Esp{ ( X - X\zzer )^2 }  = \Esp{ X^2 + (X\zzer)^2 }  = 2(1-C)  
\end{align*}
which finally gives
\begin{align*}
\abs{ \Esp{\Le_{\Ns(0, \gamma_n^2) } g(Z_n) } } \leq \frac{2 - C}{\gamma_n^2 } \norm{ g''' }_\infty
\end{align*}

$ $

For the perturbative operator, we get
\begin{align*}
 \Esp{(\Le_{H_n} - \Le_{\Ns(0, \gamma_n^2) }) g(Z_n) }  =   -\frac{C}{\gamma_n^8} \Esp{ Z_n^3 g(Z_n)  } =: -\frac{C}{\gamma_n^8} \Esp{ Z_n \tilde{g}(Z_n)  } \ \ \mbox{with} \ \ \  \tilde{g}(x) := x^2 g(x)
\end{align*}

Using iteratively the $0$-bias transform, we have
\begin{align*}
\Esp{ Z_n \tilde{g}(Z_n)  } & = \gamma_n^2 \Esp{ \tilde{g}'\prth{ Z_n\zer } }  \\
						   & =  \gamma_n^2 \Esp{ 2 Z_n\zer g\prth{ Z_n\zer } + \prth{ Z_n\zer }^2 g'\prth{ Z_n\zer } } \\
						   & = \gamma_n^2 \Esp{ 2 \Esp{ \prth{ Z_n\zer }^2 } g'\prth{ Z_n\zerR } + \prth{ Z_n\zer }^2 g'\prth{ Z_n\zer } } 
\end{align*}
where $ Z_n\zerR $ is the $0$-bias transform of $ Z_n\zer $. This last random variable is well-defined since $0$-biasing preserves the property of being of expectation $0$ for symmetric random variables, as seen in \eqref{Def:ZeroBias} taking $ f(x) = x^2 $.

As $ Z_n\zer \eqlaw \frac{1}{\gamma_n} \prth{ \sum_{k \neq I} X_k + X_I\zer } $ and $ \Esp{ \prth{ X\zzer }^2  } = 1 - 2C $, we have 
\begin{align*}
\Esp{ \prth{ Z_n\zer }^2 } = \frac{1}{\gamma_n^2} \prth{ (n-1) \Esp{ X^2 } + \Esp{ \prth{ X\zzer }^2  } } = \frac{1}{\gamma_n^2} \prth{ \gamma_n^4 - 2C }
\end{align*}
which gives 
\begin{align*}
\abs{ \Esp{(\Le_{H_n} - \Le_{\Ns(0, \gamma_n^2) }) g(Z_n) } }     \leq   \frac{ C}{\gamma_n^8} \, 3 \gamma_n^2 \Esp{ \prth{ Z_n\zer }^2 } \norm{g'}_\infty = 3C \frac{\gamma_n^4 - 2C }{\gamma_n^8 } \norm{g'}_\infty  \leq \frac{3C}{\gamma_n^4 } \norm{g'}_\infty 
\end{align*}

$ $

Finally
\begin{align*}
\abs{ \Esp{\Le_{H_n}g(Z_n) } }  \leq     \frac{2 - C }{\gamma_n^2} \norm{g'''}_\infty  +    \frac{3C}{\gamma_n^4 } \norm{g'}_\infty 
\end{align*}

$ $

To conclude with Stein's methodology, we need the estimates $ \norm{D^k g}_\infty = \big\vert\! \big\vert D^k \Le_{H_n}\inv \widehat{h}_{  \vphantom{\Le_H\inv} H_n } \big\vert\! \big\vert_\infty $ with $ k \in \ensemble{1, 3} $. Nevertheless, as we have used $ g = \Le_{H_n}\inv \widehat{h}_{  \vphantom{\Le_H\inv} H_n } $, we need to use the bound \eqref{Ineq:EstimeesAC} namely
\begin{align*}
\norm{ D \Le_{H_n}\inv \widehat{h}_{  \vphantom{\Le_H\inv} H_n } }_\infty \leq 11 \gamma \norm{  D \widehat{h}_{  \vphantom{\Le_H\inv} H_n } }_\infty = 11 \gamma \norm{  h' - \Esp{ h'(H_n) } }_\infty \leq 22 \gamma \norm{  h'  }_\infty
\end{align*}

Using in addition \eqref{Ineq:EstimeesD3}, namely
\begin{align*}
\norm{ D^3 \Le_{H_n}\inv \widehat{h}_{  \vphantom{\Le_H\inv} H_n } }_\infty \leq  \prth{ 3 + 2C + \frac{12 C }{\gamma_n^4 } } \norm{ h'' }_\infty =: A \norm{ h'' }_\infty
\end{align*}
we obtain
\begin{align}\label{Ineq:BornePresqueFinale2}
\abs{ \Esp{h(Z_n)} -\Esp{ h(H_n) } }  \leq   A \frac{2 - 3C }{\gamma_n^2} \norm{h''}_\infty +  \frac{66 C}{\gamma_n^3 } \norm{ h' }_\infty 
\end{align}
which is the desired bound.
\end{proof}


\subsection{Beyond the classical Kolmogorov approximation} 

As a corollary of theorem \ref{Thm:MASSmodG2}, we have the following

\begin{corollary}[Kolmogorov bounds]\label{Corollary:KolmogorovMASmodG2} Let $ (X_k)_k $ satisfying the hypothesis of theorem \ref{Thm:MASSmodG}. Then
\begin{align}\label{Ineq:KolmogorovBound2}
\dkol(Z_n, H_n) := \sup_{x \in \Rr} \abs{ \Prob{Z_n \leq x} - \Prob{H_n \leq x} }  \leq  \frac{ 2 (3 + 2C)^{1/3}  (2 - 3C)^{1/3}   }{ c_1^{2/3} \ \ \gamma_n^{4/3} } + O\prth{ \frac{ 1}{\gamma_n^{8/3} } }
\end{align}
%
%
%

\end{corollary}


\begin{proof}
Set 
\begin{align*}
Q(t) & := 2 t^3 - 3 t^2 + 1 \\
D(t) & :=  \Unens{ t \leq 0 } + Q(t) \Unens{ 0 \leq t \leq 1 } \\
h_{x, \delta }(y) & := D\prth{ \frac{y-x}{\delta} }
\end{align*}

For all $x$, we have $ h_{x - \delta, \delta }(y)  \leq \Unens{ y \leq x } \leq h_{x, \delta }(y) $, which implies that 
\begin{align}\label{Ineq:Dtv2}
\Esp{ h_{x - \delta, \delta }(Z_n) }  \leq \Prob{ Z_n \leq x } \leq  \Esp{ h_{x, \delta }(Z_n) } 
\end{align}

Moreover, we have
\begin{align*}
\norm{ h_{x, \delta } - \Esp{h_{x, \delta }(H_n) } }_\infty  & \leq  1   \\
\norm{ h'_{x, \delta }  }_\infty                             & \leq   \frac{1}{\delta } \\
\norm{ h''_{x, \delta }  }_\infty                            & \leq   \frac{1}{\delta^2}
\end{align*}

Using these inequalities and \eqref{Ineq:BornePresqueFinale2}, and setting $ \alpha_n :=  \prth{ 3 + 2C + \frac{12 C}{ \gamma_n^4 } } \frac{ 2 - 3C }{ \gamma_n^2 } $ and $ \beta_n := \frac{ 66 C }{\gamma_n^3 }  $, we get
\begin{align*}
\abs{ \Esp{h_{x, \delta }(Z_n)} - \Esp{ h_{x, \delta }(H_n) } }  \leq     \frac{ \alpha_n }{\delta^2}  +   \frac{ \beta_n }{\delta} 
\end{align*}

By \eqref{Ineq:Dtv2}, we have
\begin{align*}
\Prob{ Z_n \leq x }  & \leq   \Esp{ h_{x, \delta }(Z_n) } \\
					& \leq   \Esp{ h_{x, \delta }(H_n) } + \frac{ \alpha_n }{ \delta^2 } +  \frac{ \beta_n }{ \delta } \\
					& =   \Prob{ H_n \leq x } + \Esp{ Q\prth{  \frac{H_n - x}{ \delta } }  \Unens{ 0 \leq \frac{H_n - x}{ \delta } \leq 1 }  }  + \frac{ \alpha_n }{ \delta^2 } +  \frac{ \beta_n }{ \delta }  \\
					& \leq  \Prob{ H_n \leq x } +  \Prob{ 0 \leq H_n - x \leq \delta }  + \frac{ \alpha_n }{ \delta^2 } +  \frac{ \beta_n }{ \delta }  \\
					& \leq \Prob{  H_n \leq x } +  \frac{\delta}{ \gamma_n c_1 }  + \frac{ \alpha_n }{ \delta^2 } +  \frac{ \beta_n }{ \delta } 
\end{align*}

Optimising in $ \delta $ the LHS of the former inequality gives\footnote{In fact, we optimise $ \delta \mapsto \delta / (\gamma_n c_1 ) +   \alpha_n/ \delta^2   $ as an analysis of the right power of $ \gamma_n $ if one looks for $ \delta = \gamma_n^\kappa $ for a certain $ \kappa $ reveals that the term $ \beta_n / \delta $ is already much smaller.} 
\begin{align*}
\delta = \prth{  \frac{ 2\alpha_n c_1 }{ \gamma_n } }^{1/3}
\end{align*}
and the optimal value 
\begin{align*}
\Prob{ Z_n \leq x } - \Prob{  H_n \leq x }  \leq  \frac{2^{1/3} + 2^{-2/3}}{c_1^{2/3} } \prth{  \frac{ \alpha_n   }{  \gamma_n^4 }  }^{1/3}   +  \frac{ 66 C }{ (2\alpha_n c_1)^{1/3}  \gamma_n^{3 - 1/3} } 
\end{align*}

Finally, using $ 2^{1/3} + 2^{-2/3} \approx 1,\!889 \leq 2 $ and $ \alpha_n =  \frac{ (3 + 2C) (2 - 3C) }{ \gamma_n^2 } + O(\gamma_n^{-6} ) $ we have
\begin{align*}
\Prob{ Z_n \leq x } - \Prob{  H_n \leq x }  \leq  \frac{ 2 }{ c_1^{2/3} } \frac{ (3 + 2C)^{1/3}  (2 - 3C)^{1/3}   }{ \gamma_n^{4/3} }   +  O \prth{ \frac{ 1}{\gamma_n^{8/3} } }
\end{align*}

The corresponding lower bound follows from the same manipulations.
\end{proof}


\subsection{Last remarks}

\begin{remark} 
The zero-bias transform is not characteristic of the distribution $  \Hs(\Phi_C, \gamma) $. Indeed, the Stein's equation \eqref{Eq:SteinEquation} characteristic of the Gaussian distribution is equivalent to the fixed point equation in distribution 
\begin{align*}
X \sim \Ns(0, 1) \ \  \Longleftrightarrow \ \  X \eqlaw X\zzer
\end{align*}
but we do not characterise the distribution $ \Hs(\Phi_C, \gamma) $ with such a transformation. 

A natural transformation would be the following $ C $-bias transform, defined for a random variable $ W $ such that $ \Esp{W} = \Esp{W^3} = 0 $, $ \Esp{W^2} =  \gamma^2 $ and $ \Esp{W^4} < \infty $  and for all absolutely continuous functions $f$ satisfying $ \Esp{ \abs{ W^3 f(W) } } < \infty $ by
\begin{align}\label{Def:CTransfo}
\Esp{  f'\prth{ W\zzerc }  } = \prth{ \gamma^2 + \frac{4C}{\gamma^6 } \Esp{W^4} }\inv \Esp{  \rho_C(W) f(W)  }
\end{align}
with
\begin{align*}
\rho_C(x) := x + \frac{4C}{\gamma^6 } x^3
\end{align*}

The distribution of $ W\zzerc $ is absolutely continuous with respect to Lebesgue measure and has for density
\begin{align*}
f_{ W\zzerc }(x) = \Esp{ \rho_C(W) \Unens{W \geq x } }
\end{align*}

The proof of such a result is the same as in the case of the zero-bias transform, and we refer to \cite{GoldsteinReinert} or \cite{Ross} for the details.

We remark that we recover the zero-bias transform letting $ C \to 0 $. Moreover, the translation of \eqref{Def:CTransfo} in terms of a fixed point equation in law is
\begin{align*}
X \sim \Hs(\Phi_C, \gamma) \ \  \Longleftrightarrow \ \  X \eqlaw X\zzerc
\end{align*}

Unfortunately, due to the non-linearity of $ \rho_C $, the application to sums of i.i.d. random variables fails~: the $ C $-bias transform of $ S_n $ is not immediate to find, and the replacement at random of one term of the sum by an independent $ C $-biased term does not give the result.
\end{remark}


\begin{remark} \label{Rk:ImprovedSteinBound}
With the suitable rescaling of $ h $, hence of $ Z_n $ and $ H_n $, we get
\begin{align*}
\abs{ \Esp{ h\prth{ \frac{Z_n}{\gamma_n} }  } -\Esp{ h\prth{ \frac{H_n}{\gamma_n} } } }  \leq  \frac{1}{\gamma_n^4 } \prth{A (2 - 3C)  \norm{h''}_\infty +  66 C \norm{h' }_\infty } = O\prth{ \frac{1}{n} }
\end{align*}
which corresponds to an improvement of the classical Berry-Ess\'en bound. This can be understood as an additive correction to the usual norm by writing
\begin{align*}
\abs{ \Esp{ h\prth{ \frac{Z_n}{\gamma_n} }  } -\Esp{ h\prth{ \frac{H_n}{\gamma_n} } } } = \abs{ \Esp{ h\prth{ \frac{Z_n}{\gamma_n} }  } - \Esp{ h\prth{ G } }  + \operatorname{Corr}(n, h) } 
\end{align*}
with 
\begin{align*}
\operatorname{Corr}(n, h) := \Esp{ h\prth{ G } - h\prth{ \frac{H_n}{\gamma_n} } } = \Esp{ \int_{H_n/\gamma_n}^{G} h'(x)dx } = \Esp{  h'\prth{ H_n/\gamma_n + U \Delta_n } \Delta_n  }
\end{align*}
where $ U \sim \Ue([0, 1]) $, $ G $ and $ H_n $ are independent and $ \Delta_n := G - H_n/\gamma_n $.

The search for an additional correction that would give a faster approximation was also developed in \cite{ElKarouiJiao} with financial applications such as an approximation of the price of CDOs. A comparison between the two corrective terms would be interesting for the applications. 
\end{remark}


\section{Appendix : Stein's estimates}\label{Chapitre:SteinEstimates}

\subsection{Overview and main definitions}

We develop here the equivalent of the Stein's estimates that are relevant in our case by carefully adapting the steps of Stein \cite{SteinApprox}. To this goal, we first prove some in subsection \ref{Subsection:BasicEstimates} some basic estimates on the tail of the $ \Phi^4(a, b) $ distribution defined by the density $ \exp\prth{- a x^2/2 - b x^4/4 } / \Ze_{a, b} $ with $ \Ze_{a, b} := \int_\Rr \exp\prth{- a x^2/2 - b x^4/4 } dx $. We then find an integral form for the relevant operators that enter into the composition of the Stein operator and its derivatives in subsection \ref{Subsection:IntegralRepresentations}, and we finally prove in subsection \ref{Subsection:OperatorNormsEstimates}  the Stein's estimates that take the form of an operator norm estimate of the type $ \vert\vert D^k \Le\inv D^{-k'} \vert\vert_{L^\infty \to L^\infty} < \infty $ for certain integers $ k, k' $. Here, $D$ is the operator of differentiation and $ \Le $ the operator of interest. The details of these last proofs use the basic estimates.

We will adopt a set of general conventions and definitions throughout this whole chapter. We set
\begin{align}\label{Def:a_and_b}
a := \frac{1}{\gamma^2}, \qquad b := \frac{ C }{\gamma^8 }
\end{align}

Define 
\begin{align}\label{Def:z_gamma}
z_\gamma := \int_\Rr e^{ -a \frac{x^2}{2} - b \frac{x^4}{4} } dx 
\end{align}

Note that the constant $ z_\gamma $ writes with $ G_\gamma \sim \Ns(0, \gamma^2) $ as
\begin{align*} 
z_\gamma = \gamma \sqrt{2\pi }\, \Esp{ e^{ - C \frac{ G_\gamma^4 }{4 \gamma^8} }  } = \gamma \sqrt{2\pi }\, \Esp{ e^{ - C \frac{ G_1^4 }{4 \gamma^4} }  } = \gamma \sqrt{2\pi } \prth{ 1 - \frac{3C }{4 \gamma^4} + O\prth{ \frac{1}{\gamma^8 } } }
\end{align*}

In particular, 
\begin{align*} 
1 - \frac{3C }{4 \gamma^4} \leq \frac{z_\gamma}{ \gamma \sqrt{2\pi } } \leq 1
\end{align*}

We will constantly consider the random variable $ H_\gamma $ of distribution
\begin{align}\label{Def:H_gamma}
H_\gamma \sim \frac{1}{z_\gamma } e^{ -a \frac{x^2}{2} - b \frac{x^4}{4} } dx
\end{align}

We define
\begin{align}\label{Def:GeneralDefinitions}
\begin{aligned}
f_\gamma(x) & := \Esp{ \delta_0(H_\gamma - x) } = \frac{1}{z_\gamma } e^{ -a \frac{x^2}{2} - b \frac{x^4}{4} }  \\
F_\gamma(x) & := \Prob{H_\gamma \leq x } , \qquad  \overline{F}_\gamma(x) := \Prob{ H_\gamma \geq x } \\
\psi_\gamma(x) & := \Esp{H_\gamma \Unens{ H_\gamma \geq x } }  = - \Esp{H_\gamma \Unens{ H_\gamma \leq x } } \\
\varphi_\gamma(x) & := \Esp{ (x - H_\gamma)_+ } = \int_{-\infty}^x F_\gamma \\
\overline{\varphi}_\gamma (x) & := \Esp{ ( H_\gamma - x)_+ } = \int_x^{+\infty}  \overline{F}_\gamma \\ 
\rho_\gamma(x) & := a x + b x^3 , \qquad \widetilde{\rho}_\gamma(x)  := a  + b x^2 =  \rho_\gamma(x)/x \\
B_\gamma(x) & := \Bgammax{x} \\
D_\gamma(x) & := 2 \rho_\gamma'(x) + \rho_\gamma(x)^2 \\
V_\gamma(x) & := \frac{1}{2} (x^2 + \sigma_\gamma^2), \quad \sigma_\gamma^2 := \Esp{H_\gamma^2 } \\
G_\gamma(x) & := 1 + \rho_\gamma(x) \frac{ F_\gamma(x) }{f_\gamma(x) } , \qquad \overline{G}_\gamma(x)   := 1 - \rho_\gamma(x) \frac{ \overline{F}_\gamma(x) }{f_\gamma(x) } 
\end{aligned}
\end{align}

We also define
\begin{align}\label{Def:OperatorsAndFunctions}
\begin{aligned}
D & : f \in \He^1_\Phi  \mapsto f' \\
\Le_\gamma & : f \in \He^1_\Phi \mapsto f' - \rho_\gamma f \\
X & : x \in \Rr \hspace{+0.3cm} \mapsto x \\
h_\gamma & : x \in \Rr \hspace{+0.3cm} \mapsto h(x) - \Esp{ h(H_\gamma) } \\
\widehat{h}_\gamma & : x \in \Rr \hspace{+0.3cm} \mapsto h_\gamma(x) - x \Esp{ h'(H_\gamma) } = h(x) - \Esp{ h(H_\gamma) } - x \Esp{ h'(H_\gamma) }
\end{aligned}
\end{align}

Note that $ \Le_\gamma $ is invertible on $ \ensemble{ f \in  \He^2_\Phi : \, \Esp{f(H_\gamma) } = 0 } $, and in particular, one can define $ \Le_\gamma\inv h_\gamma $ and $ \Le_\gamma\inv \widehat{h}_\gamma $ as the solutions $ (g_\gamma, \widehat{g}_\gamma) $ of the Stein equations $ \Le_\gamma g_\gamma = h_\gamma $ and $ \Le_\gamma \widehat{g}_\gamma = \widehat{h}_\gamma $ that vanish in $ \pm \infty $, namely
\begin{align}\label{Def:SolutionsSteinEquationsPhiFour}
\begin{aligned}
\Le_\gamma\inv h_\gamma (x) & = \frac{\Esp{ h_\gamma(H_\gamma) \Unens{ H_\gamma \geq x } } }{ f_\gamma(x) } \\
\Le_\gamma\inv \widehat{h}_\gamma (x) & = \frac{\Esp{ \widehat{h}_\gamma(H_\gamma) \Unens{ H_\gamma \geq x } } }{ f_\gamma(x) } 
\end{aligned}
\end{align}

Last, we define
\begin{align*}
\chi_\gamma(x) := \int_{-\infty}^x \varphi_\gamma(t) dt & := \Esp{ \int_{-\infty}^x (t - H_\gamma)_+ dt} = \Esp{ \int_{-\infty}^x (t - H_\gamma) \Unens{ t - H_\gamma \geq 0 } dt} \\
                 & = \Esp{ \int_{H_\gamma}^x (t - H_\gamma) dt \Unens{ x \geq  H_\gamma  }  } = \Esp{ \frac{ (x- H_\gamma)^2 }{2 } \Unens{ x \geq H_\gamma} } \\
                 & = \frac{x^2}{2} \Prob{ H_\gamma \leq x} - x \Esp{ H_\gamma \Unens{ H_\gamma \leq x } } + \frac{1}{2} \Esp{ H_\gamma^2 \Unens{ H_\gamma \leq x } }
\end{align*}
and
\begin{align*}
\overline{\chi}_\gamma(x) := \int_x^{+\infty} \overline{ \varphi }_\gamma(t) dt & := \Esp{ \int_x^{+\infty} ( H_\gamma - t )_+ dt} = \Esp{ \int_x^{+\infty} ( H_\gamma - t) \Unens{ H_\gamma - t \geq 0 } dt} \\
                 & = \Esp{ \int_x^{H_\gamma} ( H_\gamma - t) dt \Unens{   H_\gamma \geq x }  } = \Esp{ \frac{ ( H_\gamma - x)^2 }{2 } \Unens{  H_\gamma \geq x } } \\
                 & = \frac{x^2}{2} \Prob{ H_\gamma \geq x} - x \Esp{ H_\gamma \Unens{ H_\gamma \geq x } } + \frac{1}{2} \Esp{ H_\gamma^2 \Unens{ H_\gamma \geq x } }
\end{align*}

\subsection{Basic estimates} 
\label{Subsection:BasicEstimates}


\begin{lemma}\label{Lemma:BasicEstimates1} For $ \gamma > 0 $, we have  
\begin{align}\label{Ineq:QGaussienneDilateePlus}
\forall x > 0, \ \ \ \ \  \overline{F}_\gamma(x) \leq \frac{ f_\gamma(x) }{ \rho_{\gamma}( x ) }
\end{align}
\begin{align}\label{Ineq:QGaussienneDilateeMoins}
\forall x < 0, \ \ \ \ \  F_\gamma(x)  \leq \frac{ f_\gamma(x) }{ \rho_{\gamma}( \abs{x} ) }
\end{align}

\end{lemma}


\begin{proof}
As $ \rho_\gamma'(x) = a + 3 b x^2 > 0 $, the function $ \rho_\gamma $ is strictly increasing on $ \Rr $ and we can write for $ x > 0 $
\begin{align*}
\overline{F}_\gamma(x) = \Prob{ H_\gamma \geq x } = \int_x^{+ \infty } f_\gamma(y) dy \leq \int_x^{+ \infty } \frac{ \rho_\gamma(y) }{ \rho_\gamma(x) } f_\gamma(y) dy = \frac{ f_\gamma(x)  }{ \rho_\gamma(x) }
\end{align*}
using $ f_\gamma' = - \rho_\gamma f_\gamma $ (easily seen with the definition of $ f_\gamma $) and $ \lim_{\pm \infty } \rho_\gamma f_\gamma = 0 $.

The second inequality is nothing but the first one where $x$ has been replaced by $ -x $, using the fact that $ F_\gamma(-x) = \overline{F}_\gamma(x) $ and $ \rho_\gamma(-x) = \rho_\gamma(x) $, $ f_\gamma(-x) = f_\gamma(x) $.
\end{proof}

$ $


\begin{lemma}\label{Lemma:BasicEstimates2} We have for all $ x \in \Rr $ and all $ \gamma > 0 $
\begin{align}\label{Ineq:QGaussienneXRho}
\psi_\gamma(x)  \leq \frac{ x f_\gamma(x) }{ \rho_\gamma( x ) }
\end{align}
\end{lemma}


\begin{proof}
The fact that this last inequality is symmetric compared for example to \eqref{Ineq:QGaussienneDilateeMoins} comes from the fact that $ R(x) := \rho_{\gamma}( x ) / x = a + b x^2  = R(\abs{x}) $ and the fact that the function on the left hand side is odd.

As $ (f_\gamma / R)' = - f_\gamma \prth{ R'/R^2 + \rho_{\gamma}/R } $, and $ \lim_{+ \infty}f_\gamma / R = 0 $, we have
\begin{align*}
\frac{ x f_\gamma(x) }{ \rho_{\gamma}( x ) } =  \int_x^{+ \infty } \prth{ R'/R^2 + \rho_{\gamma}/R } f_\gamma = \Esp{  \prth{ R'/R^2 + \rho_{\gamma}/R } (H_\gamma) \Unens{ H_\gamma \geq x } }
\end{align*}

But $ R(x) := \rho_{\gamma}( x ) / x $, so $ \rho_{\gamma}(x)/R (x) = x $, hence, setting $ r(x) := R'(x) /R(x)^2  $, we get
\begin{align*}
\frac{ x f_\gamma(x) }{ \rho_{\gamma}( x ) } = \Esp{   r(H_\gamma) \Unens{ H_\gamma \geq x } }  + \Esp{ H_\gamma \Unens{ H_\gamma \geq x } } =: \Esp{  r(H_\gamma) \Unens{ H_\gamma \geq x } } + \psi_\gamma(x)
\end{align*}
and it remains to show that
\begin{align*}
\Esp{ r(H_\gamma) \Unens{ H_\gamma \geq x } } \geq~0 
\end{align*}
which is clearly the case for $ x \geq 0 $ since $ R'(x) = 2bx $.

For $ x \leq 0 $, set $ r := R'/ R^2 $. As $ r(-x) = - r(x) $ and $ H_\gamma \eqlaw -H_\gamma $, one has $ r(H_\gamma) \eqlaw - r(H_\gamma) $ and in particular, $ \Esp{ r(H_\gamma) } = 0 $. Then
\begin{align*}
\Esp{ r(H_\gamma) \Unens{ H_\gamma \geq x } } & = \Esp{ r(H_\gamma) \Unens{ H_\gamma \geq -\abs{ x } } } = - \Esp{ r(-H_\gamma) \Unens{ -H_\gamma \leq \abs{ x } } } \\
			& =  - \Esp{ r(H_\gamma) \Unens{ H_\gamma \leq \abs{ x } } }  \ \ \mbox{as} \ \ r(H_\gamma) \eqlaw -r(H_\gamma) \\
			& =   \Esp{ r(H_\gamma) \Unens{ H_\gamma \geq \abs{ x } } }  \ \ \mbox{as} \ \ \Esp{r(H_\gamma)} = 0 \\
			& \geq  r(\abs{ x }) \Prob{ H_\gamma \geq \abs{ x } } \geq 0
\end{align*}
which concludes the proof.
\end{proof}

\begin{remark}
Replacing the function $ r $ by the function $ X : x \mapsto x $ in the last equalities gives $ \Esp{ H_\gamma \Unens{ H_\gamma \geq x } } \geq 0 $ for all $ x \in \Rr $.
\end{remark}

$ $


\begin{lemma}\label{Lemma:BasicEstimates3} We have for all $ x \in \Rr $  and all $ \gamma > 0 $
\begin{align}\label{Ineq:QGaussienneGnedenko}
\Prob{  H_\gamma \leq x }  & \geq -\frac{ \rho_{\gamma}(x)  }{ \rho'_{\gamma}(x) + \rho_{\gamma}(x)^2 } f_\gamma(x) \\
\Prob{  H_\gamma \geq x }  & \geq  \ \ \frac{ \rho_{\gamma}(x)  }{ \rho'_{\gamma}(x) + \rho_{\gamma}(x)^2 } f_\gamma(x) 
\end{align}
\end{lemma}


\begin{proof}

We first remark that
\begin{align*}
\Prob{  H_\gamma \leq x }  + \frac{ \rho_{\gamma}(x)  }{ \rho'_{\gamma}(x) + \rho_{\gamma}(x)^2 } f_\gamma(x) & =  \int_{- \infty }^x  \prth{ f_\gamma(u) + \frac{d}{du} \prth{ \frac{ \rho_{\gamma} f_\gamma  }{ \rho'_{\gamma} + \rho_{\gamma}^2 } }(u)  } du \\
\Prob{  H_\gamma \geq x }  - \frac{ \rho_{\gamma}(x)  }{ \rho'_{\gamma}(x) + \rho_{\gamma}(x)^2 } f_\gamma(x) & =  \int_x^{+ \infty }  \prth{ f_\gamma(u) + \frac{d}{du} \prth{ \frac{ \rho_{\gamma} f_\gamma  }{ \rho'_{\gamma} + \rho_{\gamma}^2 } }(u)  } du
\end{align*}

Hence, it is sufficient to prove that 
\begin{align*}
1 + \frac{1}{f_\gamma(x) } \frac{d}{dx} \prth{ \frac{ \rho_{\gamma} f_\gamma  }{ \rho'_{\gamma} + \rho_{\gamma}^2 } }(x) \geq 0
\end{align*}

But
\begin{align*}
1 + \frac{1}{f_\gamma } \prth{ \frac{ \rho_{\gamma} f_\gamma  }{ \rho'_{\gamma} + \rho_{\gamma}^2 } }' & = 1 + \prth{ \frac{ \rho_{\gamma}    }{ \rho'_{\gamma} + \rho_{\gamma}^2 } }' - \rho_\gamma \prth{ \frac{ \rho_{\gamma}  }{ \rho'_{\gamma} + \rho_{\gamma}^2 } } = \prth{ \frac{ \rho_{\gamma}    }{ \rho'_{\gamma} + \rho_{\gamma}^2 } }'  +  \frac{ \rho_{\gamma}'    }{ \rho'_{\gamma} + \rho_{\gamma}^2 } \\ 
			& =  \frac{1}{ \rho_{\gamma}  } \prth{ \frac{ \rho_{\gamma}^2    }{ \rho'_{\gamma} + \rho_{\gamma}^2 } }' = \frac{1}{ \rho_{\gamma}  } \prth{ \frac{ 1    }{ 1 -  \prth{ 1 / \rho_{\gamma}  }' } }'   \\
			& =  \frac{1}{ \rho_{\gamma}  }  \prth{\frac{1}{ \rho_{\gamma}  } }'' \frac{ 1    }{ \prth{ 1 -  \prth{ 1 / \rho_{\gamma}  }' }^2 }  = \frac{ 2(\rho_{\gamma}')^2  - \rho_{\gamma} \rho_{\gamma}''   }{ \rho_{\gamma}^4 \prth{ 1 -  \prth{ 1 / \rho_{\gamma}  }' }^2 }  
\end{align*}

It is now sufficient to prove that $ 2(\rho_{\gamma}')^2  - \rho_{\gamma} \rho_{\gamma}'' \geq 0 $. Setting $ Y :=  x^2 $, we have
\begin{align*}
2(\rho_{\gamma}')^2(x) - \rho_{\gamma}(x) \rho_{\gamma}''(x)  & = 2 (a + 3 b x^2 )^2 - 6 b x (a x + b x^3 ) = 2 (a + 3 b Y )^2 - 6 (a b Y + b^2 Y^2 ) \\
			& = 12 (bY)^2 + 6 abY + 2a^2 = 2a^2 \prth{ 6 \prth{ \frac{bY}{a} }^2 + 3 \prth{ \frac{bY}{a} } + 1 }
\end{align*}

As the function $ t \mapsto 6 t^2 + 3t + 1 $ is positive on $ \Rr $, we finally have the result.
\end{proof}

$ $


\begin{lemma}[Inequality on $ \varphi_\gamma + \overline{\varphi}_\gamma $]\label{Lemma:InequalityPhiPlusPhiBar} Recall that $ \widetilde{\rho}_\gamma(x) := \rho_\gamma(x)/x = a + b x^2 $. Then, for all $ x \in \Rr $, we have 
\begin{align}\label{Ineq:InequalityPhiPlusPhiBar}
\varphi_\gamma(x) + \overline{\varphi}_\gamma(x) \leq 2 \frac{f_\gamma(x) }{ \widetilde{\rho}_\gamma(x) }
\end{align}
\end{lemma}


\begin{proof}
We have
\begin{align*}
\varphi_\gamma(x) + \overline{\varphi}_\gamma(x) & = 2 \psi_\gamma(x) + x \prth{ F_\gamma(x) - \overline{F}_\gamma(x) } \\
                & \leq 2 \frac{f_\gamma(x) }{ \widetilde{\rho}_\gamma(x) } + x \prth{ F_\gamma(x) - \overline{F}_\gamma(x) } \quad\mbox{ by \eqref{Ineq:QGaussienneXRho} }
\end{align*}

If $ x < 0 $, we have $ x \prth{ F_\gamma(x) - \overline{F}_\gamma(x) } =  - \abs{x} \prth{ F_\gamma(- \abs{x}) - \overline{F}_\gamma(- \abs{x}) } = \abs{x} \prth{ F_\gamma( \abs{x}) - \overline{F}_\gamma( \abs{x}) } $ using $ F_\gamma(-x) = \overline{F}_\gamma(x) $ which is equivalent to $ H_\gamma \eqlaw - H_\gamma $.

Using $ F_\gamma  + \overline{F}_\gamma = 1 $, we have $ F_\gamma  - \overline{F}_\gamma = 1 - 2 \overline{F}_\gamma $. Moreover, $ 1 - 2 \overline{F}_\gamma \leq \varepsilon \overline{F}_\gamma $ iff $ \overline{F}_\gamma \geq \frac{1}{2 + \varepsilon } $, namely, iff $ x \geq \overline{F}_\gamma\inv\prth{ \frac{1}{2 + \varepsilon } } \geq 0 $. For all $x \in \Big[\overline{F}_\gamma\inv\prth{ \frac{1}{2 + \varepsilon } }, +\infty \Big[ $, we thus have
\begin{align*}
\varphi_\gamma(x) + \overline{\varphi}_\gamma(x) & \leq 2 \frac{f_\gamma(x) }{ \widetilde{\rho}_\gamma(x) } + \varepsilon x  \overline{F}_\gamma(x)   \leq (2 + \varepsilon) \frac{f_\gamma(x) }{ \widetilde{\rho}_\gamma(x) }  \quad\mbox{ by \eqref{Ineq:QGaussienneDilateePlus} }
\end{align*}

The result is valid for all $ \varepsilon > 0 $, hence, by continuity of $ \overline{F}_\gamma\inv $, one can pass to the limit $ \varepsilon \to 0 $ and since $ \overline{F}_\gamma\inv(1/2) = 0 $, one gets the desired result.
\end{proof}


\begin{lemma} Set $ M = \gamma / 2 $. Then, for all $ x \in \Rr  $
\begin{align}\label{Ineq:XGbar}
\overline{F}_\gamma(x) \geq  \frac{x - M}{ (x + \gamma) \rho_\gamma(x) } f_\gamma(x)
\end{align}
\end{lemma}


\begin{proof}
The inequality is obvious on $ \Rr_- $ as $ \overline{F}_\gamma(x) \geq 0 $, $ f_\gamma(x) \geq 0 $, $ (x + \gamma) \rho_\gamma(x) \geq x \rho_\gamma(x) \geq 0 $ and $ x - M \leq 0 $. Now, for $ x > 0 $
\begin{align*}
\overline{F}_\gamma(x) - \frac{x - M}{ (x + \gamma) \rho_\gamma(x) } f_\gamma(x) & = \int_x^{ +\infty } \prth{ f_\gamma(t) + \frac{d}{dt } \prth{ \frac{t - M}{ (t + \gamma) \rho_\gamma(t) } f_\gamma(t) }  } dt \\
              & = \int_x^{ +\infty } f_\gamma(t) \prth{ 1 +  \frac{d}{dt } \prth{ \frac{t - M}{ (t + \gamma) \rho_\gamma(t) } }  - \rho_\gamma(t) \frac{t - M}{ (t + \gamma) \rho_\gamma(t) }   } dt \\
              & = \int_x^{ +\infty } f_\gamma(t) \frac{r_\gamma(t)}{ (t + \gamma)^2 \rho_\gamma(t)^2 } dt
\end{align*}
with, using \eqref{Def:a_and_b}
\begin{align*}
r_\gamma(x) & := b^2 \prth{  M  + \frac{ 1 }{ \sqrt{a} }  } x^7 + \frac{ b^2 }{ \sqrt{a} } \prth{ M  + \frac{ 1}{ \sqrt{a} }  } x^6 + 2 \sqrt{a} b (  M \sqrt{a} + 1 ) x^5 +  b( 2 M \sqrt{a}   - 1 ) x^4 \\
          & \quad + \prth{ M ( a^2 + 4 b) + a^{3/2} - 2 \frac{ b }{ \sqrt{a} } } x^3 + M \prth{  a^{3/2} + 3 \frac{  b }{ \sqrt{a} } } x^2 + 2 M a x + M\sqrt{a}
\end{align*}

In order to have this last polynomial positive on $ \Rr_+ $, we need to choose $M$ so that
\begin{align*}
2 M \sqrt{a}   - 1  \geq 0 \ & \Longleftrightarrow M \geq \frac{1 }{ 2 \sqrt{a} } = \frac{\gamma}{2} \\
M ( a^2 + 4 b) + \frac{ a^2 - 2b }{ \sqrt{a} } \geq 0 \ & \Longleftrightarrow M \geq  \sqrt{a}  \frac{ 2b - a^2 }{ a^2 + 4 b } = \sqrt{a} \prth{ -1 + \frac{ 6b }{a^2 + 4b } } \geq -\sqrt{a} 
\end{align*}

Thus, for $ M = \gamma / 2 $, we have the result.
\end{proof}

\begin{remark}
Using $ \overline{G}_\gamma $ defined in \eqref{Def:GeneralDefinitions}, the inequality \eqref{Ineq:XGbar} is equivalent for all $ x \in \Rr $ to 
\begin{align}\label{Ineq:XGbarBis}
(x + \gamma) \overline{G}_\gamma(x) \leq \frac{ 3 \gamma }{2 }
\end{align}
\end{remark}

$ $

The forthcoming estimates required computer-helped computations. They were performed using the software Maple 9 and guessed using Octave 4.0.0.

$ $


\begin{lemma}[Inequality on $ \psi_\gamma $]\label{Lemma:InequalityPsiGamma} Suppose that $ \gamma^4 \geq 2C \prth{  \frac{ \sqrt{15}}{3} - 1 } $ and set $ X(x) := x $. Then, we have 
\begin{align}\label{Ineq:InequalityPsiGamma}
\sign\prth{ \psi_\gamma  - f_\gamma \frac{ \rho_\gamma + X ( \rho_\gamma^2 + 2 \rho_\gamma') }{ \rho_\gamma'' + 3 \rho_\gamma \rho_\gamma' + \rho_\gamma^3 } }  = \sign(X) 
\end{align}
\end{lemma}


\begin{proof}
Recall that $ \rho_\gamma(x) =  a x + b x^3 $ with $ a := \frac{1}{\gamma^2 }$ and $ b := \frac{C}{\gamma^8} $. Define for all $ x \in \Rr^* $
\begin{align*}
\widehat{Q}_\gamma(x) & :=  \frac{ \rho_\gamma(x) + x( \rho_\gamma(x)^2 + 2 \rho_\gamma'(x) ) }{ \rho_\gamma''(x) + 3 \rho_\gamma(x) \rho_\gamma'(x) + \rho_\gamma(x)^3 } \\
            & =  \frac{ b^2 x^6 + 2 ab x^4 + (a^2  + 7 b)  x^2 + 3 a }{   b^3 x^8 + 3 a b^2 x^6 + 3b( a^2  + 3 b ) x^4 + a(a^2 + 12 b) x^2 + 3 (a^2 + 2 b) } 
\end{align*}

As $ \psi_\gamma(-x) = -\psi_\gamma(x) $, $ f_\gamma(-x) = f_\gamma(x) $ and $ \widehat{Q}_\gamma(-x) =  \widehat{Q}_\gamma(x) $, we see that $ \psi_\gamma  - f_\gamma \widehat{Q}_\gamma \leq 0 $ on $ \Rr_- $ as a sum of negative terms. It is thus enough to prove that 
for all $ x \geq 0 $, 
\begin{align*}
\psi_\gamma(x)  + f_\gamma(x) \widehat{Q}_\gamma(x) \geq 0
\end{align*}

We have
\begin{align*}
\psi_\gamma(x)  + f_\gamma(x) \widehat{Q}_\gamma(x) & = \int_x^{+\infty}  \prth{  u f_\gamma(u) - \frac{d}{du } (f_\gamma \widehat{Q}_\gamma)(u) } du, \qquad x \geq 0 
\end{align*}
hence, it is sufficient to prove that for all $ u \geq 0 $
\begin{align*}
x - \frac{1}{f_\gamma(u) } \frac{d}{du } (f_\gamma Q_\gamma)(u) \geq 0
\end{align*}

Now, using $ f_\gamma' = - \rho_\gamma f_\gamma $, we get
\begin{align*}
x - \frac{1}{f_\gamma(x) } \frac{d}{dx } (f_\gamma Q_\gamma)(x) & = x - Q_\gamma'(x) + \rho_\gamma(x) Q_\gamma(x) \\
             & =:  \frac{ 2 x \widehat{P}_\gamma(x) }{  ( b^3 x^8 + 3 ab^2 x^6 +  3b( a^2  + 3 b ) x^4 +  a ( a^2 + 12 b ) x^2 + 3 (a^2 + 2 b) )^2 } 
\end{align*}
where 
\begin{align*}
\widehat{P}_\gamma(x) =   b^6 x^{16} & + 6 a b^5 x^{14} + \left( 15 a^2 b^4 + 18 b^5 \right) x^{12} + \left( 20 a^3 b^3 + 78 a b^4 \right) x^{10} \\
              & + \left( 15 a^4 b^2 + 132 a^2 b^3 + 93 b^4 \right) x^8 + \left( 6 a^5 b + 108 a^3 b^2 + 252 a b^3 \right) x^6 \\
              & + \left( a^6 + 42 a^4 b + 234 a^2 b^2 + 120 b^3 \right) x^4 + \left( 6 a^5 + 84 a^3 b + 132 a b^2 \right) x^2 \\
              & + 3(3 a^4 + 12 a^2 b - 8 b^2)
\end{align*}

This last polynomial has non negative coefficients if and only if $ 3 a^4 + 12 a^2 b - 8 b^2 \geq 0 $, i.e. if $ 3 \prth{\frac{a^2 }{b} }^2 + 12 \frac{a^2 }{b} - 8  \geq 0 $. The polynomial $ 3 X^2 + 12 X - 8 $ is non negative on $ [\alpha, +\infty) $ with $ \alpha := -2 + \frac{2}{3} \sqrt{15} \approx 0,\! 5819 $. Thus, we have the result for $ a^2/b \geq \alpha  $, i.e. $ \gamma^4 \geq \alpha C $.
\end{proof}

$ $

\begin{lemma}[Bound on $ V_\gamma I_2 $] For all $ x \geq 0 $ and $ \gamma \geq 1 $, we have
\begin{align}\label{Ineq:BoundOnVgammaI2}
\prth{ D_\gamma(x) - B_\gamma(x) \frac{\overline{F}_\gamma(x) }{f_\gamma(x) } } V_\gamma(x) \leq 1 + \frac{18 C }{ 10 }
\end{align}
\end{lemma}


\begin{proof}
Recall that $ \rho_\gamma(x) := a x + b x^3 $ and that $ \sigma_\gamma^2 := \Esp{H_\gamma^2 } $. Let $ d > 0 $ be a constant to be chosen such that $ (D_\gamma  - B_\gamma  \overline{F}_\gamma(x) / f_\gamma ) V_\gamma \leq d $ on $ \Rr_+ $. As $ B_\gamma, D_\gamma, V_\gamma \geq 0 $ on $ \Rr_+ $, this is equivalent to $ \overline{F}_\gamma - f_\gamma (D_\gamma V_\gamma - d ) B_\gamma / V_\gamma \geq 0 $. Define 
\begin{align*}
\widetilde{Q}_\gamma := (D_\gamma V_\gamma - d ) \frac{ B_\gamma }{ V_\gamma }
\end{align*}

Then, for all $ x \geq 0 $, using $ f_\gamma' = - \rho_\gamma f_\gamma $, we get
\begin{align*}
\overline{F}_\gamma(x) - f_\gamma(x) \widetilde{Q}_\gamma(x) = \int_x^{ +\infty } f_\gamma \prth{ 1 + \widetilde{Q}'_\gamma - \rho_\gamma \widetilde{Q}_\gamma }
\end{align*}

Moreover, setting $ \sigma \equiv \sigma_\gamma $,
\begin{align*}
1 + \widetilde{Q}'_\gamma(x) - \rho_\gamma(x) \widetilde{Q}_\gamma(x)  & = \frac{ 2  \widetilde{P}_\gamma(x^2) }{ {x}^{2}  \left( {\sigma}^{2}+{x}^{2}  \right) ^{2}}  \\
               & \quad \times \frac{1}{ \prth{ b^3 x^8 + 3 a b^2 x^6 + \left( 3 a^2 b + 9 b^2 \right) x^4 + \left( a^3 + 12 ab \right) x^2+ 3( a^2 + 2 b ) }^2  }
\end{align*}
with 
\begin{align*}
\widetilde{P}_\gamma(Y) & :=  d b^{4} Y^7 + \left( b^4 d \sigma^2 + 4 a b^3 d \right) Y^6 + 2 b^2\left( \crochet{ 2 a b  \sigma^2 + 3 a^2 + 10 b } d - 15 b \right) Y^5 \\ 
                       & \quad + b \prth{ \crochet{ 6 a^2 b  \sigma^2 + 18 b^2  \sigma^2 + 4 a^3 + 48 a b } d - 15 b \crochet{ a + 4 b \sigma^2 } } Y^4 \\
                       & \quad + \left( \crochet{ 4 a^3 b\sigma^2 + 42 a b^2 \sigma^2 + a^4 + 36 a^2 b + 69 b^2 } d - 6b\crochet{ 2 a^2 - 5 b^2 \sigma^4 - 5 a b \sigma^2 + 6 b } \right) Y^3 \\
                       & \quad + \left( \crochet{ a^4\sigma^2 + 30 a^2 b\sigma^2 + 51 b^2\sigma^2 + 8 a^3 + 66 a b } d - 3\crochet{ 5 a b^2 \sigma^4 + 8 a^2 b \sigma^2 - 24 b^2 \sigma^2 + a^3 + 2 a b } \right) Y^2 \\
                       & \quad + \left( \crochet{ 6 a^3\sigma^2 + 42 ab\sigma^2 + 9 a^2 + 18 b } d - 6 \sigma^2 \crochet{ 2 a^2 b \sigma^2 + a^3  + 2 a b - 6 b^2 \sigma^2 } \right) Y \\
                       & \quad + 3 \sigma^2 ( 2 b + a^2 ) (d - \sigma^2 a) 
\end{align*}

If one can choose $ d $ in such a way that all the coefficients of this polynomial are non negative, the desired result will follow, namely $\overline{F}_\gamma(x) - f_\gamma(x) \widetilde{Q}_\gamma(x) \geq 0 $ for all $ x \geq 0 $. For that, we obtain the set of inequations 
\begin{align*}
d & \geq \frac{ 15 b  }{ 2 a b  \sigma^2 + 3 a^2 + 10 b  } \\ 
d & \geq \frac{ 15 b \crochet{ a + 4 b \sigma^2 } }{ 6 a^2 b  \sigma^2 + 18 b^2  \sigma^2 + 4 a^3 + 48 a b } \\
d & \geq \frac{ 6b \crochet{ 2 a^2 - 5 b^2 \sigma^4 - 5 a b \sigma^2 + 6 b } }{ 4 a^3 b\sigma^2 + 42 a b^2 \sigma^2 + a^4 + 36 a^2 b + 69 b^2 }                        \\
d & \geq \frac{ 3\crochet{ 5 a b^2 \sigma^4 + 8 a^2 b \sigma^2 - 24 b^2 \sigma^2 + a^3 + 2 a b }  }{ a^4\sigma^2 + 30 a^2 b\sigma^2 + 51 b^2\sigma^2 + 8 a^3 + 66 a b  } \\
d & \geq \frac{6 \sigma^2 \crochet{ 2 a^2 b \sigma^2 + a^3  + 2 a b - 6 b^2 \sigma^2 }  }{ 6 a^3\sigma^2 + 42 ab\sigma^2 + 9 a^2 + 18 b } \\
d & \geq \sigma^2 a 
\end{align*}

We have $ a = \gamma^{-2} $ and $ b = C \gamma^{-8} $. Moreover, with $ G \sim \Ns(0, 1) $, we have 
\begin{align*}
\sigma_\gamma^2 & := \Esp{ H_\gamma^2 } = \frac{ \int_\Rr  x^2 e^{ -\frac{x^2}{2 \gamma^2 } - \frac{C}{\gamma^4} \frac{x^4}{4\gamma^4} } dx }{   \int_\Rr e^{ -\frac{x^2}{2 \gamma^2 } - \frac{C}{\gamma^4} \frac{x^4}{4\gamma^4} } dx } = \frac{ \gamma^3 \int_\Rr  t^2 e^{ -\frac{t^2}{2  } - \frac{C}{\gamma^4} \frac{t^4}{4 } } dt }{  \gamma \int_\Rr e^{ -\frac{t^2}{2  } - \frac{C}{\gamma^4} \frac{t^4}{4 } } dt } = \gamma^2 \frac{ \Esp{ G^2 e^{ - \frac{C}{\gamma^4} \frac{G^4}{4 } } } }{ \Esp{ e^{ - \frac{C}{\gamma^4} \frac{G^4}{4 } } } }
\end{align*}

In particular, 
\begin{align*}
\Esp{ G^2 \prth{ 1 - \frac{C}{\gamma^4} \frac{G^4}{4 } } } \leq \frac{\sigma_\gamma^2}{ \gamma^2 } \leq \frac{ \Esp{G^2 } }{ \Esp{  1 - \frac{C}{\gamma^4} \frac{G^4}{4 }  } } \leq \Esp{G^2 } \Esp{  1 + \frac{C}{\gamma^4} \frac{G^4}{4 }  }
\end{align*}
namely
\begin{align}\label{Ineq:VariancePhi4}
1 - \frac{15}{4 }\frac{C}{\gamma^4} \leq \frac{\sigma_\gamma^2}{ \gamma^2 } \leq    1 + \frac{3}{4 }  \frac{ C}{\gamma^4}  
\end{align}

Setting $ \varsigma_\gamma^2 := \sigma_\gamma^2/\gamma^2 $ and using \eqref{Ineq:VariancePhi4} and simple inequalities in addition to $ \gamma \geq 1 $ and $ 3 \geq C $, we get
\begin{align*}
\bullet  \quad & \frac{ 15 b  }{ 2 a b  \sigma^2 + 3 a^2 + 10 b  } = \frac{ 15 C }{ 3 \gamma^4 + 2 C \varsigma_\gamma^2 + 10 C } \leq \frac{15  }{13}  \\ 
\bullet  \quad & \frac{ 15 b \crochet{ a + 4 b \sigma^2 } }{ 6 a^2 b \sigma^2 + 18 b^2  \sigma^2 + 4 a^3 + 48 a b }    =     \frac{ 15 C }{ 2 } \frac{  \gamma^4 + 4 C \varsigma_\gamma^2  }{ 2 \gamma^8 + 3 C \varsigma_\gamma^2 \gamma^4 + 24 C \gamma^4 + 9 C^2 \varsigma_\gamma^2 } \leq \frac{15 C }{4 \gamma^4 }   \\
\bullet  \quad & \frac{ 6b \crochet{ 2 a^2 - 5 b^2 \sigma^4 - 5 a b \sigma^2 + 6 b } }{ 4 a^3 b\sigma^2 + 42 a b^2 \sigma^2 + a^4 + 36 a^2 b + 69 b^2 }    =   6 C \frac{  2 \gamma^8 +  C \gamma^4 ( 5\varsigma_\gamma^2 - 6 ) + 5 C^2 \varsigma_\gamma^4  }{ \gamma^4 \left( \gamma^8 + 4  C \gamma^4 (\varsigma_\gamma^2 + 9) + C^2 (42 \varsigma_\gamma^2 + 69 )\right) } \\
& \hspace{+6.5cm} \leq  \frac{12 C}{ \gamma^4 } \\
\bullet  \quad & \frac{ 3\crochet{ 5 a b^2 \sigma^4 + 8 a^2 b \sigma^2 - 24 b^2 \sigma^2 + a^3 + 2 a b }  }{ a^4\sigma^2 + 30 a^2 b\sigma^2 + 51 b^2\sigma^2 + 8 a^3 + 66 a b  } = 3 \frac{ \gamma^8 + 8 C \varsigma_\gamma^2 \gamma^4 + 5 C^2 \varsigma_\gamma^4 + 2 C \gamma^4 - 24 C^2 \varsigma_\gamma^2 }{ \varsigma_\gamma^2 \gamma^8 + 8 \gamma^8 + 30 C \varsigma_\gamma^2 \gamma^4 + 66 C \gamma^4 + 51 C^2 \varsigma_\gamma^2 }  \\
& \hspace{+6.9cm}\leq \frac{ 1 }{3 } + \frac{87 C^3 }{ 51 C^2 } \leq 1 + \frac{18 C }{ 10 } \\
\bullet  \quad & \frac{6 \sigma^2 \crochet{ 2 a^2 b \sigma^2 + a^3  + 2 a b - 6 b^2 \sigma^2 }  }{ 6 a^3\sigma^2 + 42 ab\sigma^2 + 9 a^2 + 18 b } = \frac{ 2 \left( 3 \gamma^8 + C \varsigma_\gamma^2 \gamma^4 + 6 C \gamma^4 - 18 C^2 \varsigma_\gamma^2 \right) \varsigma_\gamma^2 }{3{\gamma}^{4} \left( 2 \varsigma_\gamma^2 \gamma^4 + 3 \gamma^4 + 14 C \varsigma_\gamma^2 + 6 C \right) } \leq 1   \\
\bullet  \quad & \sigma^2 a  = \varsigma_\gamma^2  \leq 1 + \frac{3C}{4 \gamma^4 }
\end{align*}

Finally, we choose for instance $ d = 1 + \frac{18 C }{ 10 } $ to get the result, as for $ \gamma $ big enough ($ \gamma \geq 12 C $ more precisely), the other bounds are negligible. 
\end{proof}

$ $

\begin{lemma}[Positivity of $ G'_\gamma $ and $ \overline{G}'_\gamma $ ]\label{Lemma:PositivityDerivatives} Suppose $ \gamma^4 \geq 3C $ and set $ X(x) := x $. Then, we have 
\begin{align}\label{Ineq:PositivityDerivativesOfG}
\begin{aligned}
\sign\prth{ F_\gamma  + f_\gamma \frac{ \rho_\gamma^2 + 2 \rho_\gamma'}{ \rho_\gamma'' + 3 \rho_\gamma \rho_\gamma' + \rho_\gamma^3 } } & = \sign(X) \\
\sign\prth{ \overline{F }_\gamma   - f_\gamma \frac{ \rho_\gamma^2 + 2 \rho_\gamma'}{ \rho_\gamma'' + 3 \rho_\gamma \rho_\gamma' + \rho_\gamma^3 }   } & = - \sign(X)
\end{aligned}
\end{align}
\end{lemma}


\begin{proof}
Recall from definition \eqref{Def:GeneralDefinitions} that $ \rho_\gamma (x)  =   a x + b x^3 $ with $ a := \gamma^{-2} $ and $ b := C \gamma^{-8} $. Define for all $ x \in \Rr^* $
\begin{align*}
Q_\gamma(x) & :=  \frac{ \rho_\gamma(x)^2 + 2 \rho_\gamma'(x)}{ \rho_\gamma''(x) + 3 \rho_\gamma(x) \rho_\gamma'(x) + \rho_\gamma(x)^3 } \\
            & =  \frac{ b^2 x^6 + 2 ab x^4 + (a^2 + 6 b ) x^2 + 2 a }{ x \prth{ b^3 x^8 + 3 a b^2 x^6 + 3b( a^2  + 3 b ) x^4 + a( a^2 + 12 b) x^2 + 3 (a^2 + 2 b) } } 
\end{align*}

As $ \overline{F}_\gamma(-x) = F_\gamma(x) $, $ f_\gamma(-x) = f_\gamma(x) $ and $ Q_\gamma(-x) = -Q_\gamma(x) $, we see that the two equalities are the same up to changing $ x $ into $ -x $. It is thus enough to prove the first one.

It is enough to prove the first equality on $ \Rr^*_- $ as, on $ \Rr^*_+ $, $ F_\gamma  + f_\gamma \frac{ \rho_\gamma^2 + 2 \rho_\gamma'}{ \rho_\gamma'' + 3 \rho_\gamma \rho_\gamma' + \rho_\gamma^3 } $ is clearly positive as a sum of positive terms. This amounts to prove that for all $ x < 0 $, 
\begin{align*}
F_\gamma(x)  + f_\gamma(x) Q_\gamma(x) \leq 0
\end{align*}

We have
\begin{align*}
F_\gamma(x)  + f_\gamma(x) Q_\gamma(x) & = \int_{-\infty}^x \prth{ f_\gamma(u) + \frac{d}{du } (f_\gamma Q_\gamma)(u) } du, \qquad x < 0 
\end{align*}
hence, it is sufficient to prove that for all $ u \leq 0 $
\begin{align*}
1 + \frac{1}{f_\gamma(u) } \frac{d}{du } (f_\gamma Q_\gamma)(u) \leq 0
\end{align*}

Now, using $ f_\gamma' = - \rho_\gamma f_\gamma $, we get
\begin{align*}
1 + \frac{1}{f_\gamma(x) } \frac{d}{dx } (f_\gamma Q_\gamma)(x) & = 1 + Q_\gamma'(x) - \rho_\gamma(x) Q_\gamma(x) \\
             & = -\frac{6}{x^2} \frac{ 10 b^3 x^6 + 5 ab^2 x^4 +  4b( a^2 - 3 b ) x^2 + a(a^2 + 2 b) }{  ( b^3 x^8 + 3 ab^2 x^6 +  3b( a^2  + 3 b ) x^4 +  a ( a^2 + 12 b ) x^2 + 3 (a^2 + 2 b) )^2 } 
\end{align*}

The numerator of this rational function is a polynomial with non negative coefficients if and only if $ a^2 - 3b \geq 0 $, i.e. if $ \gamma^4 \geq 3C $, hence the result.
\end{proof}


\begin{lemma}[Inequality on $ \overline{\chi}_\gamma $]\label{Lemma:BoundOnChiBarBoverf} Suppose that $ \gamma \geq 1 $. Then, for all $ x \in \Rr_+ $
\begin{align}\label{Ineq:BoundOnChiBarBoverf}
\overline{ \chi}_\gamma(x) \leq \frac{ f_\gamma(x) }{ \Bgammax{x} } 
\end{align}
\end{lemma}


\begin{proof}
Define (with $ a := \gamma^{-2} $ and $ b = C \gamma^{-8} $)
\begin{align*}
\begin{cases}
Q_0(x) \hspace{-0.2cm} & = \frac{b x^2 }{ \Bgammax{x} }  \\
Q_k(x) \hspace{-0.2cm} & =   \frac{ Q_{k - 1}(x) }{ 1 + Q_{k - 1}'(x) }, \qquad \qquad k \in \ensemble{1, 2, 3}
\end{cases}  
\end{align*}

It is easily seen that for all $ k \in \ensemble{0, 1, 2, 3} $ and for all $ x \geq 0 $, we have $ Q_k(x) \geq 0 $ and $ 1 + Q_k'(x) \geq 0 $.

Using \eqref{Ineq:InequalityPsiGamma} and $ X(x) := x $ we have on $ \Rr_+ $
\begin{align*}
\frac{\overline{\varphi}_\gamma }{ f_\gamma } = \frac{ \psi_\gamma - X \overline{F}_\gamma  }{ f_\gamma } \leq \frac{\psi_\gamma }{f_\gamma } \leq \frac{1}{ \widetilde{\rho}_\gamma } \leq \frac{1}{ b X^2 }
\end{align*}
as $ \widetilde{\rho}_\gamma(x) := a + b x^2 \geq b x^2 $.

Recall from definition \eqref{Def:GeneralDefinitions} that  $ B_\gamma(x) := \Bgammax{x} $. As $ \overline{\varphi}_\gamma > 0 $ on $ \Rr_+ $, we have
\begin{align*}
\overline{ \chi}_\gamma \frac{ B_\gamma }{ f_\gamma } = \overline{ \chi}_\gamma \frac{ B_\gamma }{ \overline{\varphi}_\gamma } \frac{ \overline{\varphi}_\gamma }{ f_\gamma } \leq \overline{ \chi}_\gamma \frac{ B_\gamma }{ \overline{\varphi}_\gamma } \frac{1}{ b X^2 } =:  \frac{  \overline{ \chi}_\gamma }{ \overline{\varphi}_\gamma } \frac{1}{ Q_0 } 
\end{align*}

As $ B_\gamma, f_\gamma \geq 0 $, the inequality \eqref{Ineq:BoundOnChiBarBoverf} is satisfied in particular if 
\begin{align*}
\frac{ \overline{ \chi}_\gamma }{ \overline{\varphi}_\gamma } \frac{1}{ Q_0 }  \leq 1     \qquad \Longleftrightarrow \qquad     \overline{ \chi}_\gamma \leq Q_0  \overline{\varphi}_\gamma 
\end{align*}

Using $ \overline{ \chi}_\gamma(x) := \int_x^{+\infty} \overline{ \varphi }_\gamma $ and $ Q_0(x)  \overline{\varphi}_\gamma(x) = -\int_x^{+\infty} ( Q_0  \overline{\varphi}_\gamma )' $ as $ \lim_{+\infty } \overline{\varphi}_\gamma = 0 $ faster than any polynomial or rational fraction such as $ Q_0 $, this last inequality is equivalent to 
\begin{align*}
\int_x^{+\infty } \prth{ (1 + Q_0') \overline{\varphi}_\gamma   + Q_0 \overline{\varphi}_\gamma' } \leq 0
\end{align*}

Using the fact that $ \overline{\varphi}_\gamma(x) = \int_x^{ +\infty } \overline{F}_\gamma $, hence that $ \overline{\varphi}_\gamma' = -\overline{F}_\gamma $, we see that this last inequality is satisfied in particular if, on $ \Rr_+ $, 
\begin{align*}
(1 + Q_0') \overline{\varphi}_\gamma   - Q_0 \overline{F}_\gamma \leq 0     \qquad \Longleftrightarrow \qquad     \overline{ \varphi}_\gamma \leq Q_1  \overline{F}_\gamma 
\end{align*}

Using $ \lim_{+\infty} Q_0 \overline{F}_\gamma = 0 $, this last inequality is equivalent for all $ x \geq 0 $ to
\begin{align*}
\int_x^{+\infty } \prth{ (1 + Q_1') \overline{F}_\gamma   + Q_1 \overline{F}_\gamma' } \leq 0
\end{align*}
and, using $ \overline{F}_\gamma' = - f_\gamma $, it is in particular satisfied if, on $ \Rr_+ $,
\begin{align*}
(1 + Q_1') \overline{F}_\gamma   - Q_1 f_\gamma \leq 0     \qquad \Longleftrightarrow \qquad     \overline{ F }_\gamma \leq Q_2  f_\gamma 
\end{align*}

Using the fact that $ \overline{F}_\gamma(x) = \int_x^{ +\infty } f_\gamma $, this last inequality is equivalent, for all $ x \geq 0 $ to  
\begin{align*}
\int_x^{+\infty } \prth{ (1 + Q_2') f_\gamma   + Q_2 f_\gamma' } \leq 0
\end{align*}
and, using $ f_\gamma' = - \rho_\gamma f_\gamma $, it is in particular satisfied if, on $ \Rr_+ $,
\begin{align*}
(1 + Q_2') f_\gamma   - Q_1 \rho_\gamma f_\gamma \leq 0     \qquad \Longleftrightarrow \qquad     \frac{1}{\rho_\gamma}  \leq Q_3 
\end{align*}
as $ f_\gamma > 0 $.

Now, a tedious computation\footnote{The author thanks the software Maple that did all the algebraic manipulations in these proofs.} shows that
\begin{align*}
Q_3(x) - \frac{1}{\rho_\gamma(x) } =  \frac{ \gamma^8 P_1(x^2) }{ x (\gamma^6 + C x^2) P_2(x^2) } 
\end{align*}
with
\begin{align*}
P_1(Y) = &  \, C^{12} Y^{16} + 12  C^{11} \gamma^{6} Y^{15} 
           +  \left( 66  C^{10} \gamma^{12} + 36  C^{11} \gamma^{8} \right)  Y^{14} 
           +  \left( 220 C^{9} \gamma^{18} + 372 C^{10} \gamma^{14} \right)  Y^{13} \\
         & + \left( 495  C^{8} \gamma^{24} + 1740  C^{9} \gamma^{20} + 468  C^{10} \gamma^{16} \right)  Y^{12} 
           +  \left( 792  C^{7} \gamma^{30} + 4860  C^{8} \gamma^{26} + 4086 C^{9} \gamma^{22} \right)  Y^{11} \\
         & +  \left( 924  C^{6} \gamma^{36} + 9000  C^{7} \gamma^{32} + 15894 C^{8} \gamma^{28} + 2646  C^{9} \gamma^{24} \right)  Y^{10} \\
         & + \left( 792  C^{5} \gamma^{42} + 11592  C^{6} \gamma^{38} + 36216 C^{7} \gamma^{34} + 19224 C^{8} \gamma^{30} \right)  Y^9 \\
         & +  \left( 495  C^{4} \gamma^{48} + 10584  C^{5} \gamma^{44} + 53424 C^{6} \gamma^{40} + 60948  C^{7} \gamma^{36} + 6138  C^{8} \gamma^{32} \right)  Y^8 \\
         &  +  \left( 220  C^{3} \gamma^{54} + 6840  C^{4} \gamma^{50} + 53172  C^{5} \gamma^{46} + 110268  C^{6} \gamma^{42} + 38040  C^{7} \gamma^{38} \right)  Y^7 \\
         & +  \left( 66  C^{2} \gamma^{60} + 3060  C^{3} \gamma^{56} + 36036  C^{4} \gamma^{52} + 124848  C^{5} \gamma^{48} + 99063  C^{6} \gamma^{44} + 5292  C^{7}  \gamma^{40} \right)  Y^6 \\
         & +  \left( 12 C \gamma^{66} + 900  C^{2} \gamma^{62} + 16344  C^{3} \gamma^{58} + 91116  C^{4} \gamma^{54} + 140286  C^{5} \gamma^{50} + 34110  C^{6} \gamma^{46} \right)  Y^5 \\
         & +  \left(  \gamma^{72} + 156 C \gamma^{68} + 4716  C^{2} \gamma^{64} + 42444  C^{3} \gamma^{60} + 116565  C^{4} \gamma^{56} + 76506 C^{5} \gamma^{52} + 7425  C^{6} \gamma^{48} \right)  Y^4 \\
         & +  \left( 12  \gamma^{74} + 774 C \gamma^{70} + 11988  C^{2} \gamma^{66} + 56988  C^{3} \gamma^{62} + 80172  C^{4} \gamma^{58} + 29376  C^{5} \gamma^{54} \right)  Y^3 \\
         & +  \left( 54  \gamma^{76} + 1818 C \gamma^{72} + 15465  C^{2} \gamma^{68} + 41616  C^{3} \gamma^{64} + 38628 
 C^{4} \gamma^{60} + 6534  C^{5} \gamma^{56} \right)  Y^2 \\
          & + \left( 108  \gamma^{78} + 1998 C \gamma^{74} + 9942  C^{2} \gamma^{70} + 19308  C^{3} \gamma^{66} + 13104  C^{4} \gamma^{62} \right) Y \\
          & + 81  \gamma^{80} + 810 C \gamma^{76} + 3015  C^{2} \gamma^{72} + 4950  C^{3} \gamma^{68} + 3024  C^{4} \gamma^{64}
\end{align*}
and 
\begin{align*}
P_2(Y) = & \,  C^{12}Y^{16} + 12  C^{11} \gamma^{6} Y^{15} +  \left( 66  C^{10} \gamma^{12} + 35  C^{11} \gamma^{8} \right) Y^{14} +  \left( 220  C^{9} \gamma^{18} + 362  C^{10} \gamma^{14} \right) Y^{13}  \\
& + \left( 495  C^{8} \gamma^{24} + 1695  C^{9} \gamma^{20} + 441  C^{10} \gamma^{16} \right) Y^{12} +  \left( 792  C^{7} \gamma^{30} + 4740  C^{8} \gamma^{26} + 3861  C^{9} \gamma^{22} \right) Y^{11} \\
& +  \left( 924  C^{6} \gamma^{36} + 8790  C^{7} \gamma^{32} + 15066  C^{8} \gamma^{28} + 2406  C^{9} \gamma^{24} \right) Y^{10} \\
& +  \left( 792  C^{5} \gamma^{42} + 11340  C^{6} \gamma^{38} + 34452  C^{7} \gamma^{34} + 17607  C^{8} \gamma^{30} \right) Y^{9} \\
& + \left( 495  C^{4} \gamma^{48} + 10374  C^{5} \gamma^{44} + 51030  C^{6} \gamma^{40} + 56259  C^{7} \gamma^{36} + 5355  C^{8} \gamma^{32} \right) Y^{8} \\
& +  \left( 220  C^{3} \gamma^{54} + 6720  C^{4} \gamma^{50} + 51030  C^{5} \gamma^{46} + 102651  C^{6} \gamma^{42} + 33873  C^{7} \gamma^{38} \right) Y^{7}  \\
& +  \left( 66  C^{2} \gamma^{60} + 3015  C^{3} \gamma^{56} + 34776  C^{4} \gamma^{52} + 117303  C^{5} \gamma^{48} + 89874  C^{6} \gamma^{44} + 4521  C^{7} \gamma^{40} \right) Y^{6}  \\
& +  \left( 12 C \gamma^{66} + 890  C^{2} \gamma^{62} + 15876  C^{3} \gamma^{58} + 86481  C^{4} \gamma^{54} + 129513  C^{5} \gamma^{50} + 30579  C^{6} \gamma^{46} \right) Y^{5} \\
& + \left(  \gamma^{72} + 155 C \gamma^{68} + 4617  C^{2} \gamma^{64} + 40737  C^{3} \gamma^{60} + 109428  C^{4} \gamma^{56} + 70542  C^{5} \gamma^{52} + 6723  C^{6} \gamma^{48} \right) Y^{4}  \\
& +  \left( 12  \gamma^{74} + 765 C \gamma^{70} + 11649  C^{2} \gamma^{66} + 54387  C^{3} \gamma^{62} + 75564  C^{4} \gamma^{58} + 27378  C^{5} \gamma^{54} \right) Y^{3}  \\
& +  \left( 54  \gamma^{76} + 1791 C \gamma^{72} + 15006  C^{2} \gamma^{68} + 40023  C^{3} \gamma^{64} + 36894  C^{4}  \gamma^{60} + 6198  C^{5} \gamma^{56} \right) Y^{2}  \\
& +  \left( 108  \gamma^{78} + 1971 C \gamma^{74} + 9753  C^{2} \gamma^{70} + 18870  C^{3} \gamma^{66} + 12768  C^{4} \gamma^{62} \right) Y  \\
& + 81  \gamma^{80} + 810 C \gamma^{76} + 3015  C^{2} \gamma^{72} + 4950  C^{3} \gamma^{68} + 3024  C^{4} \gamma^{64}
\end{align*}

As $ P_1(x), P_2(x) \geq 0 $, this concludes the proof.
\end{proof}

$ $

\subsection{Integral representations} \label{Subsection:IntegralRepresentations}

The representation \eqref{Eq:SteinSolutionProba} of the solution of the Stein's equation is at the core of the operator norm estimates that allow to conclude the proof of theorems \eqref{Thm:MASSmodG} and \eqref{Thm:MASSmodG2}. We give two other related integral representations for the functions $ h_\gamma $ and $ \widehat{h}_\gamma $ in addition to the particularisation of \eqref{Eq:SteinSolutionProba} in the case of interest.


\begin{lemma}[Integral representation of $ h_\gamma $] Let $ h $ be an absolutely continuous function and $ h_\gamma $ defined in \eqref{Def:GeneralDefinitions}. Then, 
\begin{align}\label{Eq:IntegralFormH}
h_\gamma (x)  =  \int_{-\infty}^x h'(u) F_\gamma(u) du -  \int_x^{+\infty} h'(u) \overline{F}_\gamma(u) du  
\end{align}
namely
\begin{align}\label{Def:PremierNoyau}
h_\gamma (x)  = \int_\Rr h'(u) K_H(x, u) du, \qquad  K_H(x, u) := \Esp{ \Unens{ H_\gamma < u < x } - \Unens{ H_\gamma > u > x } }  
\end{align}

\end{lemma}


\begin{proof}
Write
\begin{align*}
h_\gamma (x)  & =  \Esp{ h(x) -  h(H_\gamma) }  \\
					& = \Esp{ \int_{H_\gamma}^x h'(u) du \, \Unens{ H_\gamma < x } } + \Esp{ \int_{H_\gamma}^x h'(u) du \, \Unens{ H_\gamma > x } } \\
					& = \Esp{ \int_\Rr h'(u) \prth{ \Unens{ H_\gamma < u < x } - \Unens{ H_\gamma > u > x } } du } \\
					& = \int_\Rr h'(u) \Esp{ \Unens{ H_\gamma < u < x } - \Unens{ H_\gamma > u > x } } du =: \int_\Rr h'(u) K_H(x, u) du
\end{align*}
where the last equality comes from the Fubini theorem (since $ \int_\Rr \abs{h'} < \infty $), and where 
\begin{align*} 
K_H(x,u) := \Esp{ \Unens{ H_\gamma < u < x } - \Unens{ H_\gamma > u > x } } := F_\gamma(u) \Unens{ u < x } -  \overline{F}_\gamma(u) \Unens{ u > x } 
\end{align*}
which is another form of \eqref{Def:PremierNoyau}.
\end{proof}


\begin{lemma}[Integral expression of $ \widehat{h}_\gamma $] Let $ h \in \Ce^2_0 $ be such that $ \norm{h''}_\infty < \infty $.  Recall from \eqref{Def:GeneralDefinitions} that for all $ x \in \Rr $, $ \widehat{h}_\gamma(x) := h(x) - \Esp{ h(H_\gamma) } - x \Esp{ h'(H_\gamma) } $
%
%
Then, 
\begin{align}\label{Eq:IntegralFormHatH}
\begin{aligned}
\widehat{h}_\gamma(x) & = x \prth{ \int_{-\infty}^x h''(u) F_\gamma(u) du - \int_x^{+\infty} h''(u) \overline{F}_\gamma(u) du } \\
             & \qquad - \prth{ \int_{-\infty}^x h''(u) \varphi_\gamma(u) du - \int_x^{+\infty} h''(u) \overline{\varphi}_\gamma(u) du }
\end{aligned}
\end{align}

One can rewrite this last equality into  
\begin{align*}
\widehat{h}_\gamma(x)  = \int_\Rr h''(u) \widehat{K}_\gamma (x, u ) du, \qquad \widehat{K}_\gamma(x, u) := \Esp{ (H_\gamma + x - u) \crochet{ \Unens{ H_\gamma \geq u \geq x } - \Unens{ H_\gamma \leq u \leq x } } }  
\end{align*}
\end{lemma}


\begin{proof}
We have
\begin{align*}
\widehat{h}_\gamma(x) & := h(x) - \Esp{h(H_\gamma) } - x \Esp{ h'(H_\gamma) } \\
                  & = \int_{-\infty}^x h' F_\gamma - \int_x^{+\infty} h' \overline{F}_\gamma - x\Esp{ h'H_\gamma) } \\
                  & = \int_{-\infty}^x \prth{h'(t) - \Esp{h'(H_\gamma) } } F_\gamma(t) dt - \int_x^{+\infty} \prth{ h'(t) - \Esp{h'(H_\gamma) } } \overline{F}_\gamma(t) dt \quad \mbox{using \eqref{Eq:SumIntFH} } \\
                  & = \int_{-\infty}^x \prth{ \int_{-\infty}^t h'' F_\gamma - \int_t^{+\infty} h'' \overline{F}_\gamma } F_\gamma(t) dt   -    \int_x^{+\infty} \prth{  \int_{-\infty}^t h'' F_\gamma - \int_t^{+\infty} h'' \overline{F}_\gamma } \overline{F}_\gamma(t) dt
\end{align*}

Then, using the Fubini theorem (valid since $ \int_\Rr \abs{h''} < \infty $ and $ F_\gamma, \overline{F}_\gamma \in \crochet{0, 1} $)
\begin{align*}
\bullet \int_{-\infty}^x \prth{ \int_{-\infty}^t h''(u) F_\gamma(u) du } F_\gamma(t) dt  & = \int_{-\infty}^x h''(u) F_\gamma(u) \prth{ \int_\Rr \Unens{ u < t < x } F_\gamma(t) dt } du \hspace{+2cm}  \\
               & = \int_{-\infty}^x h''(u) F_\gamma(u) \prth{ \varphi_\gamma(x) - \varphi_\gamma(u) } du \\
\bullet \int_{-\infty}^x \prth{ \int_t^{+\infty}  h''(u) \overline{F}_\gamma(u) du } F_\gamma(t) dt  & = \int_\Rr h''(u) \overline{F}_\gamma(u) \prth{ \int_\Rr \Unens{ t < x, u } F_\gamma(t) dt } du   \\
               & = \int_\Rr h''(u) \overline{F}_\gamma(u)  \varphi_\gamma(x \wedge u) du \\
               & = \int_\Rr h''(u) \overline{F}_\gamma(u)  \prth{ \varphi_\gamma(x) \Unens{x \leq u } + \varphi_\gamma(u) \Unens{u \leq x} } du \\
               & = \varphi_\gamma(x) \int_x^{+\infty } h''(u) \overline{F}_\gamma(u)  du + \int_{-\infty }^x h''(u) \overline{F}_\gamma(u)\varphi_\gamma(u)     du
\end{align*}
\begin{align*}
\bullet \int_x^{+\infty} \prth{  \int_{-\infty}^t h''(u) F_\gamma(u) du  } \overline{F}_\gamma(t) dt & = \int_\Rr h''(u) F_\gamma(u) \prth{ \int_\Rr \Unens{ t \geq u, x } \overline{F}_\gamma(t) dt } du \\
               & = \int_\Rr h''(u) F_\gamma(u) \overline{\varphi}_\gamma(x \vee u) du \\
               & = \int_\Rr h''(u) F_\gamma(u) \prth{ \overline{\varphi}_\gamma(x) \Unens{ x > u } + \overline{\varphi}_\gamma(u) \Unens{ u > x } } du \\
               & = \overline{\varphi}_\gamma(x) \int_{-\infty}^x  h''(u) F_\gamma(u) du  +  \int_x^{+\infty} h''(u) F_\gamma(u) \overline{\varphi}_\gamma(u)  du \\
\bullet \int_x^{+\infty} \prth{  \int_t^{+\infty }  h''(u) \overline{F}_\gamma(u) du  } \overline{F}_\gamma(t) dt & = \int_x^{+\infty} h''(u) \overline{F}_\gamma(u) \prth{ \int_x^{+\infty} \Unens{ x < t < u  } \overline{F}_\gamma(t) dt } du \\	
               & = \int_x^{+\infty}  h''(u) \overline{F}_\gamma(u) \prth{ \overline{\varphi}_\gamma(x) - \overline{\varphi}_\gamma(u) } du 
\end{align*}

Thus, 
\begin{align*}
\widehat{h}_\gamma(x) & = \int_{-\infty}^x \prth{ \int_{-\infty}^t h'' F_\gamma - \int_t^{+\infty} h'' \overline{F}_\gamma } F_\gamma(t) dt   -    \int_x^{+\infty} \prth{  \int_{-\infty}^t h'' F_\gamma - \int_t^{+\infty} h'' \overline{F}_\gamma } \overline{F}_\gamma(t) dt \\
                  	& =: \int_{-\infty}^x h''(u) A_1(x, u) du + \int_x^{+\infty} h''(u) A_2(x, u) du
\end{align*}
with, using $ F_\gamma + \overline{F}_\gamma = 1 $ and $ \varphi_\gamma(x) - \overline{\varphi}_\gamma(x) = x $,
\begin{align*}
A_1(x, u) & = F_\gamma(u) \prth{ \varphi_\gamma(x) - \varphi_\gamma(u) } - \overline{F}_\gamma(u) \varphi(u) - F_\gamma(u) \overline{\varphi}_\gamma(x)  = x F_\gamma(u) - \varphi(u) \\
A_2(x, u) & = - \overline{F}_\gamma(u) \varphi_\gamma(x) - F_\gamma(u) \varphi_\gamma(u) + \overline{F}_\gamma(u) \prth{ \overline{\varphi}_\gamma(x) - \overline{\varphi}_\gamma(u) } = - \prth{ x \overline{F}_\gamma(u) + \overline{\varphi}_\gamma(u)  }
\end{align*}
which gives the result. Now, one can write
\begin{align*}
\widehat{h}_\gamma(x) & = \int_\Rr h''(u) \prth{ x F_\gamma(u) \Unens{ u \leq x } - x \overline{F}_\gamma(u) \Unens{ u \geq x }  - \varphi_\gamma(u)\Unens{ u \leq x } + \overline{\varphi}_\gamma(u) \Unens{ u \geq x } } du \\
                  	& = \int_\Rr h''(u) \Esp{ x   \Unens{ H_\gamma \leq u \leq x } - x  \Unens{ H_\gamma \geq u \geq x }  - (u - H_\gamma)_+ \Unens{ u \leq x } + (H_\gamma - u)_+ \Unens{ u \geq x } } du \\
                  	& = \int_\Rr h''(u) \Esp{ x   \Unens{ H_\gamma \leq u \leq x } - x  \Unens{ H_\gamma \geq u \geq x }  + (  H_\gamma - u ) \Unens{ H_\gamma \leq u \leq x } + (H_\gamma - u)  \Unens{ H_\gamma \geq u \geq x } } du \\
                  	& = \int_\Rr h''(u) \Esp{ (x + H_\gamma - u)  \crochet{ \Unens{ H_\gamma \leq u \leq x } -    \Unens{ H_\gamma \geq u \geq x }   } } du 
\end{align*}
hence the result.
\end{proof}


\begin{lemma}[Integral representation of $ \Le_\gamma\inv h_\gamma $] Let $ h $ be absolutely continuous and $ h_\gamma $ defined in \eqref{Def:GeneralDefinitions}. Then, 
\begin{align}\label{Eq:IntegralFormLinvH}
\Le_\gamma\inv h_\gamma (x)  =  - \frac{1}{f_\gamma(x) } \prth{  \overline{F}_\gamma(x) \int_{-\infty}^x h'(u) F_\gamma(u) du  +  F_\gamma(x)\int_x^{+\infty} h'(u) \overline{F}_\gamma(u) du  }
\end{align}
namely
\begin{align}\label{Def:DeuxiemeNoyau}
\begin{aligned}
& \Le_\gamma\inv h_\gamma (x)  = \int_\Rr h'(u) \tilde{K}_H(x,u) du, \\
& \hspace{+3cm}  \tilde{K}_H(x,u) = -\frac{1}{  f_\gamma (x)  } \Esp{ \Unens{ H^{(2)}_\gamma < u < x < H^{(1)}_\gamma } + \Unens{ H^{(2)}_\gamma > u > x > H^{(1)}_\gamma } }  
\end{aligned}
\end{align}
where $ H^{(1)}_\gamma $ and $ H^{(1)}_\gamma $ are two independent copies of $ H_\gamma $.

\end{lemma}


\begin{proof}
Using \eqref{Def:SolutionsSteinEquationsPhiFour}, we have
\begin{align*}
\Le_\gamma\inv h_\gamma (x)  & =  \frac{1}{ f_\gamma (x) }\Esp{ h_\gamma (H_\gamma) \Unens{  H_\gamma  \leq x } }  = - \frac{1}{ f_\gamma (x) }\Esp{ h_\gamma (H_\gamma) \Unens{   H_\gamma  \geq x } } \\
					& =  \frac{1}{ 2 f_\gamma (x) } \Esp{ h_\gamma (H_\gamma) \crochet{  \Unens{  H_\gamma  \leq x } -  \Unens{   H_\gamma  \geq x } }  } \\
					& =: \frac{1}{  f_\gamma (x) } \Esp{ I\prth{ x, H_\gamma } h_\gamma (H_\gamma)   } 
\end{align*}
where
\begin{align*}
I(x,y) :=  \frac{1}{2} \prth{ \Unens{y < x} - \Unens{y > x} }
\end{align*}


It follows that 
\begin{align*}
\Le_\gamma\inv h_\gamma (x)  & =  \frac{1}{ f_\gamma (x) }\Esp{ I\prth{ x, H_\gamma }  h_\gamma (H_\gamma)  }   \\
					& =  \frac{1}{ f_\gamma(x) } \Esp{  I\prth{ x, H_\gamma }  \int_\Rr  K_H(H_\gamma, u) h'(u) du } \\
					& = \int_\Rr \frac{ \Esp{ I\prth{ x, H_\gamma } K_H(H_\gamma , u) } }{ f_\gamma (x) } h'(u) \, du  =: \int_\Rr \widetilde{K}_H(x,u) h'(u) \, du
\end{align*}
where the last equality comes from the Fubini theorem. 

Let $ H^{(1)}_\gamma, H^{(2)}_\gamma $ be two independent random variables equal in law to $ H_\gamma $. Then, 
\begin{align*}
f_\gamma (x)  &  \widetilde{K}_H(x,u) :=  \Esp{ I\prth{ x, H_\gamma } K_H(H_\gamma, u) } = \Esp{ I\prth{ x, H_\gamma } K_H(H_\gamma, u) \crochet{ \Unens{ x < u } + \Unens{ u < x  } } } \\
						& = \frac{1}{2} \Esp{  \crochet{  \Unens{  H^{(1)}_\gamma  \leq x } -  \Unens{  H^{(1)}_\gamma  \geq x } } \crochet{ \Unens{ H^{(2)}_\gamma < u < H^{(1)}_\gamma } - \Unens{ H^{(2)}_\gamma > u > H^{(1)}_\gamma } } \crochet{ \Unens{ x < u } + \Unens{ u < x  } } }  \\
						& = \frac{1}{2} \Ee \bigg( 0 + \Unens{  H^{(2)}_\gamma  > u >  H^{(1)}_\gamma  > x } - \Unens{  H^{(2)}_\gamma  < u <  H^{(1)}_\gamma, \, H^{(1)}_\gamma  > x, \, x < u }  -  \Unens{  H^{(2)}_\gamma  > u > x >  H^{(1)}_\gamma } \\
						&    \qquad\quad  + \Unens{  H^{(2)}_\gamma  < u <  H^{(1)}_\gamma  < x } + 0 - \Unens{  H^{(2)}_\gamma  < u < x < H^{(1)}_\gamma }  -  \Unens{  H^{(2)}_\gamma  > u >  H^{(1)}_\gamma , \, H^{(1)}_\gamma  < x, \, x > u } \bigg) \\
						& = \frac{1}{2} \Ee \bigg( \Unens{  H^{(2)}_\gamma  > u >  H^{(1)}_\gamma  > x } - \Unens{  H^{(2)}_\gamma  < x < u <  H^{(1)}_\gamma  } - \Unens{ x < H^{(2)}_\gamma  < u <  H^{(1)}_\gamma  }  - \Unens{  H^{(2)}_\gamma  > u > x >  H^{(1)}_\gamma } \\
						&     \qquad\quad   + \Unens{  H^{(2)}_\gamma  < u <  H^{(1)}_\gamma  < x } - \Unens{  H^{(2)}_\gamma  < u < x < H^{(1)}_\gamma }  -  \Unens{  H^{(2)}_\gamma > x > u >  H^{(1)}_\gamma }  -  \Unens{  x > H^{(2)}_\gamma  > u >  H^{(1)}_\gamma } \bigg) \\
						& = - \Esp{ \Unens{ H^{(2)}_\gamma < u < x < H^{(1)}_\gamma } + \Unens{ H^{(2)}_\gamma > u > x > H^{(1)}_\gamma } }
\end{align*}
this last equality coming from the exchangeability of $ (H^{(1)}_\gamma, H^{(2)}_\gamma) $. We thus get
\begin{align*}
\tilde{K}_H(x,u) = -\frac{1}{  f_{ H_\gamma } (x)  } \Esp{ \Unens{ H^{(2)}_\gamma < u < x < H^{(1)}_\gamma } + \Unens{ H^{(2)}_\gamma > u > x > H^{(1)}_\gamma } }
\end{align*}

We can write this last operator in the following form
\begin{align*}
\widetilde{K}_H(x,u) = - \frac{1}{ f_{ \vphantom{\Le_H\inv} H_\gamma }(x) } \prth{ F_{ H_\gamma }(u) \overline{F}_{ H_\gamma }(x) \Unens{ u < x } + F_{ H_\gamma }(x) \overline{F}_{ H_\gamma }(u) \Unens{ x < u }  }
\end{align*} 
which is equivalent to \eqref{Eq:IntegralFormLinvH}.
\end{proof}


\begin{lemma}[Integral expression of $ \Le_\gamma\inv\widehat{h}_\gamma $] Let $ h \in \Ce^2_0 $ be such that $ \norm{h''}_\infty < \infty $. For all $ x \in \Rr $, we have
\begin{align}\label{Eq:IntegralFormLeInvHatH}
\begin{aligned}
\Le_\gamma\inv\widehat{h}_\gamma(x) &  = - \frac{\psi_\gamma(x)}{ f_\gamma(x) } \prth{ \int_{-\infty}^x h''(u)  F_\gamma(u) du - \int_x^{+\infty}  h''(u)  \overline{F}_\gamma(u) du } \\
             & \qquad + \frac{1}{ f_\gamma(x) } \prth{ \overline{F}_\gamma(x) \int_{-\infty}^x h''(u)  \varphi_\gamma(u) du  -    F_\gamma(x) \int_x^{+\infty} h''(u)    \overline{\varphi}_\gamma(u) du } 
\end{aligned}
\end{align} 

One can rewrite this last equality into  
\begin{align*}
& \widehat{h}_\gamma(x)  = \int_\Rr h''(u) \widehat{K}'_\gamma (x, u ) du \\
\hspace{-1cm} \mbox{with } &   \ \widehat{K}'_\gamma(x, u) := \frac{-1}{ f_\gamma(x) } \Esp{ \prth{  H^{(1)}_\gamma +  H^{(2)}_\gamma - u } \!\! \crochet{  \Unens{ H^{(1)}_\gamma \leq u \leq x \leq H^{(2)}_\gamma } -    \Unens{ H^{(1)}_\gamma \geq u \geq x \geq H^{(2)}_\gamma } } }  
\end{align*}
and where $ H^{(1)}_\gamma, H^{(2)}_\gamma $ are two independent random variables equal in law to $ H_\gamma $.
\end{lemma}


\begin{proof}
We have 
\begin{align*}
\Le_\gamma\inv\widehat{h}_\gamma(x) & = \frac{1}{ f_\gamma(x) } \Esp{ \widehat{h}_\gamma(H_\gamma) I(x, H_\gamma) }, \qquad I(x, y) := \frac{1}{2} \prth{ \Unens{y < x} - \Unens{ x < y } } \\
\widehat{h}_\gamma(x) & = \int_\Rr h''(u) \widehat{K}_\gamma(x, u) du, \qquad \widehat{K}_\gamma(x, u)  := \Esp{ (H_\gamma + x - u) \crochet{ \Unens{ H_\gamma \geq u \geq x } - \Unens{ H_\gamma \leq u \leq x } } }  
\end{align*}
hence
\begin{align*}
\Le_\gamma\inv\widehat{h}_\gamma(x) & = \frac{1}{ f_\gamma(x) } \int_\Rr h''(u) \Esp{ I(x, H_\gamma) \widehat{K}_\gamma(H_\gamma, u) }  du  \\
               & =  \int_\Rr h''(u) \widehat{K}'_\gamma(x, u)    du, \qquad \widehat{K}'_\gamma(x, u) := \frac{\Esp{ I(x, H_\gamma) \widehat{K}_\gamma(H_\gamma, u) }  }{ f_\gamma(x) }
\end{align*}

Let $ H^{(1)}_\gamma, H^{(2)}_\gamma $ be as specified. Then, 
\begin{align*}
f_\gamma (x)  &  \widehat{K}'_\gamma(x,u) :=  \Esp{ I\prth{ x, H_\gamma } \widehat{K}_\gamma(H_\gamma, u) } = \Esp{ I\prth{ x, H_\gamma } \widehat{K}_\gamma(H_\gamma, u) \crochet{ \Unens{ x < u } + \Unens{ u < x  } } } \\
						& = \frac{1}{2} \Ee \bigg( \prth{ H^{(1)}_\gamma + H^{(2)}_\gamma - u } \crochet{  \Unens{  H^{(1)}_\gamma  \leq x } -  \Unens{  H^{(1)}_\gamma  \geq x } } \\
						& \qquad \qquad \times \crochet{ \Unens{ H^{(2)}_\gamma < u < H^{(1)}_\gamma } - \Unens{ H^{(2)}_\gamma > u > H^{(1)}_\gamma } } \crochet{ \Unens{ x < u } + \Unens{ u < x  } } \bigg)  \\
						& = \frac{1}{2} \Ee \bigg( \prth{ H^{(1)}_\gamma + H^{(2)}_\gamma - u } \bigg[ 0 + \Unens{  H^{(2)}_\gamma  > u >  H^{(1)}_\gamma  > x } - \Unens{  H^{(2)}_\gamma  < u <  H^{(1)}_\gamma, \, H^{(1)}_\gamma  > x, \, x < u }  \\
						& \hspace{+5cm} -  \Unens{  H^{(2)}_\gamma  > u > x >  H^{(1)}_\gamma } + \Unens{  H^{(2)}_\gamma  < u <  H^{(1)}_\gamma  < x } + 0 \\
						&    \hspace{+5cm}   - \Unens{  H^{(2)}_\gamma  < u < x < H^{(1)}_\gamma }  -  \Unens{  H^{(2)}_\gamma  > u >  H^{(1)}_\gamma , \, H^{(1)}_\gamma  < x, \, x > u } \bigg] \bigg) \\
\end{align*}
\begin{align*}
		\hspace{+1cm}	
		     			& = \frac{1}{2} \Ee \bigg( \!\! \prth{ H^{(1)}_\gamma + H^{(2)}_\gamma - u } \!\! \bigg[ \Unens{  H^{(2)}_\gamma  > u >  H^{(1)}_\gamma  > x } - \Unens{  H^{(2)}_\gamma  < x < u <  H^{(1)}_\gamma  } - \Unens{ x < H^{(2)}_\gamma  < u <  H^{(1)}_\gamma  }  \\
						&     \hspace{+4.5cm}   - \Unens{  H^{(2)}_\gamma  > u > x >  H^{(1)}_\gamma } +  \Unens{  H^{(2)}_\gamma  < u <  H^{(1)}_\gamma  < x } - \Unens{  H^{(2)}_\gamma  < u < x < H^{(1)}_\gamma }  \\
						&     \hspace{+4.5cm} -  \Unens{  H^{(2)}_\gamma > x > u >  H^{(1)}_\gamma }  -  \Unens{  x > H^{(2)}_\gamma  > u >  H^{(1)}_\gamma } \bigg] \bigg) \\
						& = - \Esp{ \!  \prth{ H^{(1)}_\gamma + H^{(2)}_\gamma - u } \!\! \crochet{ \Unens{ H^{(2)}_\gamma < u < x < H^{(1)}_\gamma } + \Unens{ H^{(2)}_\gamma > u > x > H^{(1)}_\gamma } } }
\end{align*}
the last equality coming from the exchangeability of $ (H^{(1)}_\gamma, H^{(2)}_\gamma) $ and the fact that the function $ (A, B) \mapsto A + B - u $ is symmetric.

We can write this last operator as
\begin{align*}
- f_\gamma (x) \widehat{K}'_\gamma(x,u) & := \Esp{ \!  \prth{ H^{(1)}_\gamma + H^{(2)}_\gamma - u } \!\! \crochet{ \Unens{ H^{(2)}_\gamma < u < x < H^{(1)}_\gamma } + \Unens{ H^{(2)}_\gamma > u > x > H^{(1)}_\gamma } } } \\
              & =  \Esp{ \!  H^{(1)}_\gamma   \crochet{ \Unens{ H^{(2)}_\gamma < u < x < H^{(1)}_\gamma } + \Unens{ H^{(2)}_\gamma > u > x > H^{(1)}_\gamma } } } \\
              & \qquad + \Esp{ \!  H^{(2)}_\gamma   \crochet{ \Unens{ H^{(2)}_\gamma < u < x < H^{(1)}_\gamma } + \Unens{ H^{(2)}_\gamma > u > x > H^{(1)}_\gamma } } } \\
              & \qquad - u \, \Esp{  \Unens{ H^{(2)}_\gamma < u < x < H^{(1)}_\gamma } + \Unens{ H^{(2)}_\gamma > u > x > H^{(1)}_\gamma }   } \\
              & = \Unens{u < x} \prth{ F_\gamma(u) \psi_\gamma(x) - \psi_\gamma(u) \overline{F}_\gamma(x) - u F_\gamma(u) \overline{F}_\gamma(x) } \\
              & \qquad + \Unens{ u > x } \prth{ - \overline{F}_\gamma(u) \psi_\gamma(x) + \psi_\gamma(u) F_\gamma(x) - u \overline{F}_\gamma(u) F_\gamma(x) } \\
              & = \Unens{u < x} \prth{ F_\gamma(u) \psi_\gamma(x) - \varphi_\gamma(u) \overline{F}_\gamma(x) }  + \Unens{ u > x } \prth{ - \overline{F}_\gamma(u) \psi_\gamma(x) +  \overline{\varphi}_\gamma(u)   F_\gamma(x) }
\end{align*}
which is equivalent to \eqref{Eq:IntegralFormLeInvHatH}.
\end{proof}


\subsection{Operator norms estimates} \label{Subsection:OperatorNormsEstimates}

This is the main technical tool in Stein's methodology to come back to the initial function introduced in the Stein's equation. This amounts to prove that a particular operator is bounded on the unit sphere of the relevant topology, namely $ \vert\vert D^k \Le\inv D^{-k'} \vert\vert_{L^\infty \to L^\infty} < \infty $ for certain $ k, k' $ that we detail. As the operator involves a random variable $ H_\gamma $ whose parameter tends to $ +\infty $, we need a precise estimate of these norms as a function of $ \gamma $. This differs dramatically from the usual Stein's method where an abstract boundedness would be enough if one is not concerned about the optimal constant. In what follows, the problem of finding the optimal constant will not be tackled.

\begin{lemma}[Operator norms estimates]\label{Lemma:BorneSteinOperateur}  
For $ h \in \He^2_\Phi $, recall the definitions of $ \Le_\gamma $, $h_\gamma $ and $ \widehat{h}_\gamma $ given in \eqref{Def:OperatorsAndFunctions} and the definitions of $ \Le_\gamma\inv h_\gamma $ and $ \Le_\gamma\inv \widehat{h}_\gamma $ given in \eqref{Def:SolutionsSteinEquationsPhiFour}. Then, we have, for all $ \gamma \geq 3C $
\begin{enumerate}
\item If $h$ is bounded, i.e. $ \norm{h}_\infty < \infty $, 
\begin{align}\label{Ineq:EstimeesBornees}
\begin{aligned}
\norm{ \Le_\gamma\inv h_\gamma }_\infty   & \leq   \gamma \, \sqrt{\frac{ \pi }{2} } \norm{ h_\gamma }_\infty \\
\norm{ D\Le_\gamma\inv h_\gamma }_\infty  & \leq   2 \norm{ h_\gamma }_\infty  
\end{aligned}  
\end{align}
\item If $h$ is absolutely continuous, i.e. $ \norm{Dh}_\infty < \infty $ and $ \int_\Rr \abs{h'} < \infty $,
\begin{align}\label{Ineq:EstimeesAC}
\begin{aligned}
\norm{D \Le_\gamma\inv h_\gamma }_\infty     & \leq  11 \gamma   \norm{ D h  }_\infty  \\
\norm{D^2 \Le_\gamma\inv h_\gamma }_\infty   & \leq    4 \norm{ D h  }_\infty  
\end{aligned}
\end{align}
\item If in addition $ h' $ is absolutely continuous, i.e. $ \norm{D^2 h}_\infty < \infty $ and $ \int_\Rr \abs{h''} < \infty $,  
\begin{align}\label{Ineq:EstimeesD3}
\norm{D^3 \Le_\gamma\inv \widehat{h}_\gamma }_\infty   \leq   \prth{  3 + 2C + \frac{ 12 C}{ \gamma^4 } } \norm{ D^2 h  }_\infty   
\end{align}
\end{enumerate}

\end{lemma}

$ $

\noindent We start by proving the first part of \eqref{Ineq:EstimeesBornees}, namely 
\begin{align*}
\norm{ \Le_\gamma\inv h_\gamma }_\infty   \leq   \gamma \, \sqrt{\frac{ \pi }{2} } \norm{ h_\gamma }_\infty 
\end{align*}

\begin{proof}

Using the representation \eqref{Eq:SteinSolutionProba}, we have for all $ x \in \Rr $
\begin{align*}
\abs{\Le_\gamma\inv h_\gamma(x) }  = \frac{ \Esp{  h_\gamma(H_\gamma)    \Unens{ H_\gamma \geq x } } }{ f_\gamma(x) } & \leq  \frac{ \Esp{  \abs{ h_\gamma(H_\gamma) }  \Unens{ H_\gamma \geq x } } }{ f_\gamma(x) }    \leq  \norm{ h_\gamma }_\infty \frac{ \Esp{  \Unens{ H_\gamma \geq x } } }{ f_\gamma(x) }  
\end{align*}
namely
\begin{align*}
\abs{\Le_\gamma\inv h_\gamma(x) }  \leq \norm{ h_\gamma }_\infty   \frac{ \overline{F}_\gamma(x) }{  f_\gamma(x)  } 
\end{align*}

As $ H_\gamma \eqlaw - H_\gamma $, we have $ \overline{F}_\gamma(-x) = F_\gamma(x) $ and in particular, $ \abs{\Le_\gamma\inv h_\gamma(x) } \leq \norm{ h_\gamma }_\infty  F_\gamma(\abs{x}) / f_\gamma(x) $ if $ x \leq 0 $. The fact that $ \overline{F}_\gamma / f_\gamma $ and $ F_\gamma / f_\gamma $ reach their maximum in $ 0 $ respectively on $ \Rr_+ $ and $ \Rr_- $ follows from \eqref{Ineq:QGaussienneDilateePlus} and \eqref{Ineq:QGaussienneDilateeMoins}, and the fact that (using $ f_\gamma' = - \rho_\gamma f_\gamma $)
\begin{align*}
\prth{ \frac{ \overline{F}_\gamma }{  f_\gamma }  }' & = -1 + \rho_\gamma \frac{ \overline{F}_\gamma }{  f_\gamma } \\
\prth{ \frac{ F_\gamma }{  f_\gamma }  }' & =  1 + \rho_\gamma \frac{  F_\gamma }{  f_\gamma } 
\end{align*}

Hence, by \eqref{Ineq:QGaussienneDilateePlus} and \eqref{Ineq:QGaussienneDilateeMoins}, $ \overline{F}_\gamma / f_\gamma $ is decreasing on $ \Rr_+ $ and $ F_\gamma /  f_\gamma $ is increasing on $ \Rr_- $ ; they thus reach their maxima in $0$ on these respective sets. Using the fact that $ F_\gamma(0) = \overline{F}_\gamma(0) = \frac{1}{2} $ that is a simple corollary of $ H_\gamma \eqlaw - H_\gamma $, we get
\begin{align*}
\frac{\abs{\Le_\gamma\inv h_\gamma(x) } }{ \norm{ h_\gamma }_\infty } \leq   \frac{ 1/2 }{ f_\gamma(0) } = \frac{ z_\gamma }{ 2  } = \frac{ \gamma }{2 } \sqrt{2\pi } \, \Esp{ e^{-C G^4/(4 \gamma^4) } } \leq \gamma \sqrt{ \frac{2}{\pi } }
\end{align*}
as $ \Esp{  \exp\prth{-C G^4/(4 \gamma^4) } } \leq 1 $ (note that the optimal constant is in fact $ z_1/2 $ as $ \gamma \in \Rr_+ \mapsto \Esp{ \exp\prth{- \frac{C}{4} \frac{G^4 }{\gamma^4} } } $ is decreasing). 
\end{proof}

$ $

\noindent We now prove the second part of \eqref{Ineq:EstimeesBornees}, namely 
\begin{align*}
\norm{ D\Le_\gamma\inv h_\gamma }_\infty  \leq   2 \norm{ h_\gamma }_\infty  
\end{align*}

\begin{proof}

As $ \Le_\gamma\inv h_\gamma $ is the solution of the Stein's equation $ \Le_\gamma \Le_\gamma\inv h_\gamma = h_\gamma$, we have
\begin{align*}
D\Le_\gamma\inv h_\gamma = h_\gamma + \rho_\gamma \Le_\gamma\inv h_\gamma
\end{align*}

Thus, for $ x \geq 0 $
\begin{align*}
\abs{ D\Le_\gamma\inv h_\gamma (x) }    &  \leq   \abs{ h_\gamma(x) }  +  \abs{ \rho_\gamma(x) } \abs{ \Le_\gamma\inv h_\gamma(x) }     \\
						 & \leq  \norm{ h_\gamma }_\infty \prth{ 1 + \sup_{x > 0 } \ensemble{ \rho_\gamma(x)  \frac{ \Prob{ H_\gamma \geq x } }{ f_\gamma(x) }  } } \\
						 & \leq   2 \norm{ h_\gamma }_\infty  \quad \mbox{by \eqref{Ineq:QGaussienneDilateePlus} }
\end{align*}

We proceed in the same way for $ x < 0 $.
\end{proof}

$ $

\noindent We prove the first part of \eqref{Ineq:EstimeesAC}, namely 
\begin{align*}
\norm{D \Le_\gamma\inv h_\gamma }_\infty  \leq  11 \gamma   \norm{ D h  }_\infty  
\end{align*}

\begin{proof}
We can write 
\begin{align*}
D \Le_\gamma\inv h_\gamma(x) & = h_\gamma (x) + \rho_\gamma(x) \Le_\gamma\inv h_\gamma(x) \\
              & = \int_{-\infty}^x \!\! h' F_\gamma - \int_x^{+\infty}  \!\! h' F_\gamma - \rho_\gamma(x) \frac{ \overline{F}_\gamma(x) }{f_\gamma(x) } \int_{-\infty}^x \! \! h' F_\gamma -  \rho_\gamma(x) \frac{ F_\gamma(x) }{f_\gamma(x) } \int_x^{+\infty} \!\! h' \overline{F}_\gamma \\
              & = \prth{1 - \rho_\gamma(x) \frac{ \overline{F}_\gamma(x) }{f_\gamma(x) } } \int_{-\infty}^x \! \! h' F_\gamma - \prth{1 + \rho_\gamma(x) \frac{ F_\gamma(x) }{f_\gamma(x) } } \int_x^{+\infty} \!\! h' \overline{F}_\gamma 
\end{align*}

We know that $ 1 - \rho_\gamma  \overline{F}_\gamma  / f_\gamma \geq 0 $ on $ \Rr_+ $ by \eqref{Ineq:QGaussienneDilateePlus}. On $ \Rr_- $, we use the fact that $ x = -\abs{x } $ and  $ \rho_\gamma(-x) = - \rho_\gamma(x) $ to get $ 1 - \rho_\gamma(-\abs{x }) \overline{F}_\gamma(x)  / f_\gamma(x) = 1 + \rho_\gamma( \abs{x }) \overline{F}_\gamma(x)  / f_\gamma(x) \geq 0 $ as $  \overline{F}_\gamma(x)  ,  f_\gamma(x), \rho_\gamma(\abs{x}) \geq 0$. In the same way, using \eqref{Ineq:QGaussienneDilateeMoins}, we have $ 1 + \rho_\gamma   F_\gamma  / f_\gamma   \geq 0 $ on $ \Rr $. We can thus write
\begin{align*}
\abs{ D \Le_\gamma\inv h_\gamma(x) } \leq \norm{h'}_\infty \crochet{ \prth{1 - \rho_\gamma(x) \frac{ \overline{F}_\gamma(x) }{f_\gamma(x) } } \int_{-\infty}^x \! \!   F_\gamma + \prth{1 + \rho_\gamma(x) \frac{ F_\gamma(x) }{f_\gamma(x) } } \int_x^{+\infty} \!\! \overline{F}_\gamma }
\end{align*}

Recall from definition \eqref{Def:GeneralDefinitions} that $ \varphi_\gamma := \int_{-\infty}^\cdot F_\gamma $ and $ \overline{\varphi}_\gamma := \int_\cdot^{+\infty} \overline{F}_\gamma $ and set
\begin{align*}
\tau_\gamma(x) & := \prth{1 - \rho_\gamma(x) \frac{ \overline{F}_\gamma(x) }{f_\gamma(x) } } \varphi_\gamma(x) + \prth{1 + \rho_\gamma(x) \frac{ F_\gamma(x) }{f_\gamma(x) } }  \overline{\varphi}_\gamma(x) 
\end{align*}

Using $ \psi_\gamma(x) := \Esp{H_\gamma \Unens{ H_\gamma \geq x } } $ defined in \eqref{Def:GeneralDefinitions}, the Fubini theorem and $ \Esp{H_\gamma} = 0 $, we have
\begin{align}\label{Eq:PhiAndBarPhiWithPsi}
\begin{aligned}
\varphi_\gamma(x) = \int_{-\infty}^x F_\gamma   & =  \int_\Rr \Esp{  \Unens{H_\gamma \leq u \leq x }  } du  = \Esp{ (x-H_\gamma)_+ } =  \Esp{ H_\gamma \Unens{H_\gamma \geq x } } + x F_{ H_\gamma }(x) \\
\overline{\varphi}_\gamma(x) = \int_x^{+\infty} \overline{F}_\gamma  & =  \int_\Rr \Esp{  \Unens{ x \leq u \leq H_\gamma }  } du  = \Esp{ ( H_\gamma - x)_+ } =  \Esp{ H_\gamma \Unens{H_\gamma \geq x } } - x \overline{F}_{ H_\gamma }(x)
\end{aligned}
\end{align}
which implies 
\begin{align*}
\tau_\gamma(x) =  \prth{1 - \rho_\gamma(x) \frac{ \overline{F}_\gamma(x) }{f_\gamma(x) } } \prth{ \psi_\gamma(x) + x F_\gamma(x)  } + \prth{1 + \rho_\gamma(x) \frac{ F_\gamma(x) }{f_\gamma(x) } } \prth{ \psi_\gamma(x) - x \overline{F}_\gamma(x) } 
\end{align*}

Using the fact that $ f_\gamma(-x) = f_\gamma(x) $, $ \rho_\gamma(-x) = -\rho_\gamma(x) $ and $ \overline{F}_\gamma(-x) = F_\gamma(x) $, it is easily seen that
\begin{align*}
\tau_\gamma(-x) = \tau_\gamma(x)
\end{align*}

It is thus enough to prove that $ \tau_\gamma $ is bounded on $ \Rr_+ $ and, by symmetry, we will have the result. Using \eqref{Ineq:QGaussienneXRho} and the positivity of $ 1 - \rho_\gamma  \overline{F}_\gamma  / f_\gamma $ and $ 1 + \rho_\gamma  F_\gamma  / f_\gamma $ on $ \Rr $, we get
\begin{align*}
\tau_\gamma(x) & \leq   \prth{1 - \rho_\gamma(x) \frac{ \overline{F}_\gamma(x) }{f_\gamma(x) } } \prth{ x \frac{f_\gamma(x) }{\rho_\gamma(x) } + x F_\gamma(x)  } + \prth{1 + \rho_\gamma(x) \frac{ F_\gamma(x) }{f_\gamma(x) } } \prth{ x \frac{f_\gamma(x) }{\rho_\gamma(x) } - x \overline{F}_\gamma(x) } \\
               & = 2x \frac{f_\gamma(x) }{\rho_\gamma(x) } \prth{ 1 + \rho_\gamma(x) \frac{ F_\gamma(x) }{f_\gamma(x) } } \prth{ 1 - \rho_\gamma(x) \frac{ \overline{F}_\gamma(x) }{ f_\gamma(x) } } \\
               & =: \frac{2x }{ \rho_\gamma(x) } f_\gamma(x) G_\gamma(x) \overline{G}_\gamma(x) \qquad \mbox{using \eqref{Def:GeneralDefinitions} }
\end{align*}


Let $ \varepsilon > 0 $. We now study two cases~:

$ $

\noindent\textbf{\underline{Case $ x \in \left[ \varepsilon, +\infty \right) $ }  } \linebreak
There exists $ \beta_\varepsilon \equiv \beta_{\varepsilon, \gamma } $ such that, for all $ x \geq \varepsilon $, 
\begin{align*}
1 + \rho_\gamma(x) \frac{ F_\gamma(x) }{f_\gamma(x) } \leq (1 + \beta_\varepsilon ) \rho_\gamma(x) \frac{ F_\gamma(x) }{f_\gamma(x) } 
\end{align*}

Indeed, as $ \rho_\gamma F_\gamma / f_\gamma $ is increasing, since, for all $ x \in \Rr_+ $
\begin{align*}
G_\gamma'(x) & = \prth{ \rho'_\gamma(x) + \rho_\gamma(x)^2 } \frac{ F_\gamma(x) }{f_\gamma(x) } + \rho_\gamma(x) \geq 0 \quad \mbox{by \eqref{Ineq:QGaussienneGnedenko} } \\
\end{align*}
this amounts to take $ \beta(\varepsilon) $ equal to 
\begin{align*}
\beta(\varepsilon) = \max_{x \geq \varepsilon } \ensemble{  \frac{f_\gamma(x) }{ F_\gamma(x) \rho_\gamma(x) } } =  \frac{f_\gamma(\varepsilon ) }{ F_\gamma(\varepsilon ) \rho_\gamma(\varepsilon ) } \leq  \frac{f_\gamma( 0 ) }{ F_\gamma( 0 ) \rho_\gamma(\varepsilon ) }  \leq  \frac{2}{\sqrt{ \pi} }   \frac{ \gamma }{\varepsilon }
\end{align*}

We have used the fact that $ F_\gamma(0) = 1/2 $ and $ f_\gamma(0) = z_\gamma\inv $ with $ z_\gamma = \gamma \sqrt{2 \pi } \, \Ee( \exp\prth{ - C G^4/ (4 \gamma^4) } ) \geq  \gamma \sqrt{2 \pi } ( 1 - \frac{3C}{4 \gamma^4 } ) \geq \gamma \sqrt{ \pi } $ for $ \gamma \geq 1 $ (as $ C \leq 3 $).

We thus have 
\begin{align*}
\tau_\gamma(x) & \leq \frac{2x }{ \rho_\gamma(x) } f_\gamma(x) G_\gamma(x) \overline{G}_\gamma(x) \\
               & \leq \frac{2x }{ \rho_\gamma(x) } f_\gamma(x)(1 + \beta_\varepsilon ) \rho_\gamma(x) \frac{ F_\gamma(x) }{f_\gamma(x) }  \overline{G}_\gamma(x) \\
               & \leq 2 \prth{ 1 +  \frac{2}{\sqrt{\pi} } \frac{ \gamma }{\varepsilon } } x \overline{G}_\gamma(x) \quad \mbox{ as } F_\gamma \leq 1
\end{align*}

Moreover, we have $ x \overline{G}_\gamma(x) \leq (x + \gamma) \overline{G}_\gamma(x) \leq 3 \gamma / 2 $ using \eqref{Ineq:XGbarBis}. We can thus conclude that for all $ x \geq \varepsilon $, 
\begin{align*}
\tau_\gamma(x) \leq  3 \gamma \prth{ 1  +   \frac{2}{\sqrt{\pi} } \frac{\gamma}{\varepsilon } } 
\end{align*}

\noindent\textbf{\underline{Case $ x \in \crochet{0, \varepsilon} $ }  } \linebreak
Using the monotonicity of $ f_\gamma $, $ G_\gamma $ and $ \overline{G}_\gamma = G_\gamma(- \cdot ) $, we obtain, for all $ x \in \crochet{0, \varepsilon} $
\begin{align*}
\tau_\gamma(x) & \leq \frac{ 2 x }{ \rho_\gamma(x) } f_\gamma(x) G_\gamma(x) \overline{G}_\gamma(x) =  \frac{ 2 \gamma^2  }{ 1 + C x^2 / \gamma^6 } f_\gamma(x) G_\gamma(x) \overline{G}_\gamma(x) \\
               & \leq  2 \gamma^2  f_\gamma(0) G_\gamma(\varepsilon ) \overline{G}_\gamma(0) \\
               & = \frac{2 \gamma }{ z_\gamma / \gamma } \prth{1 + \rho_\gamma(\varepsilon) \frac{F_\gamma(\varepsilon) }{ f_\gamma(\varepsilon) }} \\
               & \leq \frac{2 \gamma }{ \sqrt{\pi } } \prth{1 +  \frac{ \rho_\gamma(\varepsilon) }{ f_\gamma(\varepsilon) } }
\end{align*}
and
\begin{align*}
\frac{ \rho_\gamma(\varepsilon) }{ f_\gamma(\varepsilon) } =  z_\gamma \frac{\varepsilon}{\gamma^2} \prth{ 1 + C \frac{\varepsilon^2}{\gamma^6 } } e^{ \frac{\varepsilon^2}{2 \gamma^2 } + C \frac{ \varepsilon^4 }{4 \gamma^8 } }
\end{align*}

If $ \varepsilon = o(\gamma) $, this last quantity is a $ o(1) $ when $ \gamma \to +\infty $. If $ \varepsilon = O(\gamma) $, this quantity is a $ O(1) $, more precisely, if $ \varepsilon = \gamma $, we obtain
\begin{align*}
\frac{ \rho_\gamma(\gamma) }{ f_\gamma(\gamma) }  \leq 2 \sqrt{2\pi e }
\end{align*}

\noindent\textbf{\underline{General case : $ x \in \Rr_+ $ }  }  \linebreak
Setting $ \varepsilon = \gamma $, we thus get
\begin{align*}
\tau_\gamma(x) & \leq \max\ensemble{ \max_{ [0, \gamma] } \tau_\gamma ,  \max_{ [ \gamma, + \infty) } \tau_\gamma } \\
               & \leq \gamma \max\ensemble{ 3\prth{ 1 + \frac{2}{\sqrt{\pi} } } , \frac{2 }{ \sqrt{\pi } } \prth{1 +  2 \sqrt{2\pi e }  } } \\
               & \leq 11 \gamma
\end{align*}
as $ \frac{2 }{ \sqrt{\pi } } \prth{1 +  2 \sqrt{2\pi e }  } \approx 10,\!455 $ and $ 3\prth{ 1 + \frac{2}{\sqrt{\pi} } } \approx 6,\!385 $.
\end{proof}

\begin{remark}
Numerical simulations support the fact that $ \tau_\gamma $ is decreasing on $ \Rr_+ $, hence that the optimal value for the constant is $ \tau_\gamma(0) = \Esp{\abs{H_\gamma } } = \gamma/2 + O(\gamma\inv) $. Moreover, the function $ f_\gamma G_\gamma \overline{G}_\gamma $ is numerically seen to be decreasing, with value in $ 0 $ equal to $ z_\gamma\inv $. Using $ 2x/\rho_\gamma(x) \leq 2 \gamma^2 $, one would thus get the constant $ \sqrt{2/\pi} \approx 0,\!79788 $.
\end{remark}

$ $

\noindent We now prove the second part of \eqref{Ineq:EstimeesAC}, namely 
\begin{align*}
\norm{D^2 \Le_\gamma\inv h_\gamma }_\infty \leq 4 \norm{ D h  }_\infty 
\end{align*}

\begin{proof}
We start by expressing $ D^2 \Le_{H_n}\inv h\hn $ in terms of $ h' $. First, we differentiate $ \Le_\gamma \prth{ \Le_\gamma\inv h_\gamma  } = h_\gamma $ to get
\begin{align}\label{Eq:DeriveeSeconde}
D^2 \Le_{H_\gamma}\inv h_{ \vphantom{\Le_H\inv} H_\gamma}   & =  D\prth{ \rho_\gamma \Le_{H_\gamma}\inv h_{ \vphantom{\Le_H\inv} H_\gamma} + h_{ \vphantom{\Le_H\inv} H_\gamma}  }  \notag \\
              &  =  \rho'_\gamma \Le_{H_\gamma}\inv h_{ \vphantom{\Le_H\inv} H_\gamma} + \rho_\gamma D \Le_{H_\gamma}\inv h_{ \vphantom{\Le_H\inv} H_\gamma} + h' \notag \\
						& =  \prth{ \rho'_\gamma + \rho_\gamma^2 } \Le_{H_\gamma}\inv h_{ \vphantom{\Le_H\inv} H_\gamma} + \rho_\gamma  h_{ \vphantom{\Le_H\inv} H_\gamma} + h' 
\end{align}


We already have by \eqref{Def:PremierNoyau} 
\begin{align*}
h_\gamma (x) = \int_\Rr h'(u) \Esp{ \Unens{ H_\gamma < u < x } - \Unens{ H_\gamma > u > x } } du =: \int_\Rr h'(u) K_H(x, u) du
\end{align*}
and by \eqref{Def:DeuxiemeNoyau}
\begin{align*}
\Le_\gamma\inv h_\gamma (x)  = \int_\Rr   \frac{ \Esp{ I\prth{ x, H_\gamma } K_H(H_\gamma , u) } }{ f_{ \vphantom{\Le_H\inv} H_\gamma } (x) } h'(u) \, du  =: \int_\Rr \widetilde{K}_H(x,u) h'(u) \, du
\end{align*}
with $ \widetilde{K}_H $ defined in \eqref{Def:DeuxiemeNoyau}.

From \eqref{Eq:DeriveeSeconde}, \eqref{Def:PremierNoyau} and \eqref{Def:DeuxiemeNoyau}, and setting 
\begin{align*}
K \star h (x) := \int_\Rr K(x,y) h(y) dy
\end{align*}
we get
\begin{align*}
D^2 \Le_\gamma\inv h_\gamma   & =  h' + \prth{ \rho'_\gamma + \rho_\gamma^2 } \Le_\gamma\inv h_\gamma  + \rho_\gamma  h_\gamma \\
           & = h' + \prth{ \rho'_\gamma + \rho_\gamma^2 } \tilde{K}_H \star h' + \rho_\gamma K_H \star h'\\
		   & =: h' + K \star h' \\
		   & \hspace{+2cm} \mbox{with } \ K(x,y) := \prth{ \rho'_\gamma + \rho_\gamma^2 }(x) \tilde{K}_H(x,y) + \rho_\gamma(x) K_H(x,y)
\end{align*}

This last operator also writes with \eqref{Eq:IntegralFormH} and \eqref{Eq:IntegralFormLinvH} as
\begin{align*}
K \star h'(x) & :=  \int_\Rr  \prth{  \prth{ \rho'_\gamma(x) + \rho_\gamma^2(x) } \tilde{K}_H (x,y)  + \rho_\gamma(x)  K_H(x,y)  } h'(y) dy \\
			  & =    \prth{ \rho_\gamma - \frac{\rho'_\gamma + \rho_\gamma^2 }{ f_\gamma } \overline{F}_\gamma }\!\! (x) \!\! \int_{-\infty}^x \!\!\! F_\gamma (y) h'(y)  dy \\ 
			  & \qquad - \prth{ \rho_\gamma + \frac{\rho'_\gamma + \rho_\gamma^2 }{ f_\gamma } F_\gamma }\!\!(x) \!\! \int_x^{+\infty} \!\!\! \overline{F}_\gamma (y) h'(y) dy  \\
			  & =: -K^{(-)} \star h'(x) - K^{(+)} \star h'(x)
\end{align*}

Using \eqref{Ineq:QGaussienneGnedenko} and the obvious fact that $ F_\gamma $ and $ \overline{F}_\gamma $ are positive, we have 
\begin{align*}
\abs{ K \star h'(x) } \leq  K^{(-)} \star \abs{ h' } (x) + K^{(+)} \star \abs{ h' } (x) \leq \norm{ h'}_\infty \prth{ K^{(-)} + K^{(+)} } \star \Un(x)
\end{align*}
with
%
%
%
%
%
%
\begin{align*}
\prth{ K^{(-)} + K^{(+)} } \star \Un =  -\prth{ \rho_\gamma - \frac{\rho'_\gamma + \rho_\gamma^2 }{ f_\gamma } \overline{F}_\gamma } \varphi_\gamma + \prth{ \rho_\gamma + \frac{\rho'_\gamma + \rho_\gamma^2 }{ f_\gamma } F_\gamma } \overline{\varphi}_\gamma 
\end{align*}
where $ \varphi_\gamma := \int_{-\infty}^\cdot F_\gamma $ and $ \overline{\varphi}_\gamma := \int_\cdot^{+\infty} \overline{F}_\gamma $ were defined in \eqref{Def:GeneralDefinitions}.

Using \eqref{Eq:PhiAndBarPhiWithPsi}, we have
\begin{align}\label{Eq:SumIntFH}
\varphi_\gamma(x) - \overline{\varphi}_\gamma(x) =  x \prth{F_\gamma (x) + \overline{F}_\gamma (x) }   = x  
\end{align}

We thus deduce, using $ \psi_\gamma(x) := \Esp{ H_\gamma \Unens{ H_\gamma \geq x } } $
\begin{align*}
\prth{ K^{(-)} + K^{(+)} } \star \Un(x)  & =  \frac{\rho'_\gamma (x) + \rho_\gamma^2(x) }{ f_\gamma(x) } \psi_\gamma(x)  - x \rho_\gamma (x)  \\
					  & \leq  \frac{\rho'_\gamma (x) + \rho_\gamma^2(x) }{ \rho_\gamma (x)/x }  - x \rho_\gamma (x)        \ \ \mbox{using \eqref{Ineq:QGaussienneXRho} } \\
					  & =  \frac{ x \rho'_\gamma (x)  }{ \rho_\gamma (x) } = 3 - \frac{1}{1 +  C x^2/\gamma^6 } \\
					  & \leq  3
\end{align*}

Finally, we get
\begin{align*}
\abs{ D^2 \Le_\gamma\inv h_\gamma (x) }  & =  \abs{ h'(x) + K \star h'(x) } \leq \abs{ h'(x) } + \abs{  K \star h'(x) } \\
					& \leq   \norm{ h' }_\infty \prth{ 1 + \prth{ K^{(-)} + K^{(+)} } \star \Un(x) } \\
					& \leq   4 \norm{ h' }_\infty
\end{align*}
which concludes the proof.
\end{proof}

$ $

We finally prove \eqref{Ineq:EstimeesD3}, namely
\begin{align*}
\norm{D^3 \Le_\gamma\inv \widehat{h}_\gamma}_\infty   \leq   \prth{  3 + 2C + \frac{ 12 C}{ \gamma^4 } } \norm{ D^2 h  }_\infty   
\end{align*}

\begin{proof}
Using $ \widehat{h}_\gamma(x) := h(x) - \Esp{ h(H_\gamma) } - x \Esp{ h'(H_\gamma) } $, we can write
\begin{align*}
D^3 \Le_\gamma\inv \widehat{h}_\gamma & =   D D^2 \Le_\gamma\inv \widehat{h}_\gamma = D\prth{ \widehat{h}_\gamma' + (\rho_\gamma' + \rho_\gamma^2) \Le_\gamma\inv\widehat{h}_\gamma + \rho_\gamma \widehat{h}_\gamma } \\
             &  = h'' +  (\rho_\gamma'' + 3 \rho_\gamma \rho_\gamma' + \rho_\gamma^3 ) \Le_\gamma\inv\widehat{h}_\gamma + (2 \rho_\gamma' + \rho_\gamma^2) \widehat{h}_\gamma + \rho_\gamma \widehat{h}_\gamma'
\end{align*}

Now, using \eqref{Eq:IntegralFormH} applied to $ \widehat{h}_\gamma' $, \eqref{Eq:IntegralFormHatH} and \eqref{Eq:IntegralFormLeInvHatH}, we have
\begin{align*}
\widehat{h}_\gamma'(x) & = h'(x) - \Esp{h'(H_\gamma) } = \int_{-\infty}^x h'' F_\gamma - \int_x^{+\infty} h'' \overline{F}_\gamma \\
\widehat{h}_\gamma(x) & = x \prth{ \int_{-\infty}^x h''(u) F_\gamma(u) du - \int_x^{+\infty} h''(u) \overline{F}_\gamma(u) du } \\
             & \qquad - \prth{ \int_{-\infty}^x h''(u) \varphi_\gamma(u) du - \int_x^{+\infty} h''(u) \overline{\varphi}_\gamma(u) du } \\
\Le_\gamma\inv\widehat{h}_\gamma(x) &  = - \frac{\psi_\gamma(x)}{ f_\gamma(x) } \prth{ \int_{-\infty}^x h''(u)  F_\gamma(u) du - \int_x^{+\infty}  h''(u)  \overline{F}_\gamma(u) du } \\
             & \qquad + \frac{1}{ f_\gamma(x) } \prth{ \overline{F}_\gamma(x) \int_{-\infty}^x h''(u)  \varphi_\gamma(u) du  -    F_\gamma(x) \int_x^{+\infty} h''(u)    \overline{\varphi}_\gamma(u) du }
\end{align*}
hence, using $ X(x) := x $ as in definition \eqref{Def:OperatorsAndFunctions}, and $ B_\gamma := \rho_\gamma'' + 3 \rho_\gamma \rho_\gamma' + \rho_\gamma^3 $, $ D_\gamma := 2 \rho_\gamma' + \rho_\gamma^2 $ as in definition \eqref{Def:GeneralDefinitions}, we obtain
\begin{align*}
D^3 \Le_\gamma\inv \widehat{h}_\gamma  - h'' & =  (\rho_\gamma'' + 3 \rho_\gamma \rho_\gamma' + \rho_\gamma^3 ) \Le_\gamma\inv\widehat{h}_\gamma + (2 \rho_\gamma' + \rho_\gamma^2) \widehat{h}_\gamma + \rho_\gamma \widehat{h}_\gamma' =  B_\gamma \Le_\gamma\inv\widehat{h}_\gamma + D_\gamma \widehat{h}_\gamma + \rho_\gamma \widehat{h}_\gamma' \\
              & = \frac{ B_\gamma}{ f_\gamma } \prth{ -  \psi_\gamma \crochet{ \int_{-\infty}^\cdot  \!\!\!\! h'' F_\gamma -  \int_\cdot^{+\infty}  \!\!\!\! h'' \overline{F}_\gamma }  + \overline{F}_\gamma \int_{-\infty}^\cdot  \!\!\!\! h'' \varphi_\gamma - F_\gamma  \int_\cdot^{+\infty}  \!\!\!\! h'' \overline{\varphi}_\gamma  }      \\
              & \qquad + D_\gamma \prth{  X \crochet{ \int_{-\infty}^\cdot  \!\!\!\! h'' F_\gamma  -    \int_\cdot^{+\infty}  \!\!\!\! h'' \overline{F}_\gamma }  -   \int_{-\infty}^\cdot  \!\!\!\! h'' \varphi_\gamma +  \int_\cdot^{+\infty}  \!\!\!\! h'' \overline{\varphi}_\gamma   } \\
              & \qquad + \rho_\gamma \crochet{ \int_{-\infty}^\cdot  \!\!\!\! h'' F_\gamma  -  \int_\cdot^{+\infty}  \!\!\!\! h'' \overline{F}_\gamma  }  \\
              & = \crochet{ \int_{-\infty}^\cdot h'' F_\gamma  - \int_\cdot^{+\infty} h'' \overline{F}_\gamma  }\prth{ \rho_\gamma + X D_\gamma - \frac{\psi_\gamma}{f_\gamma } B_\gamma } \\
              & \quad + \prth{ \int_{-\infty}^\cdot h'' \varphi_\gamma } \prth{  - D_\gamma + \frac{\overline{F}_\gamma}{f_\gamma } B_\gamma  }  + \prth{  \int_\cdot^{+\infty} h'' \overline{\varphi}_\gamma  } \prth{ D_\gamma + \frac{ F_\gamma}{f_\gamma } B_\gamma  }
\end{align*}

Thus, we get
\begin{align*}
\abs{ D^3 \Le_\gamma\inv \widehat{h}_\gamma  }    & \leq    \norm{ h'' }_\infty \Bigg( 1 + \crochet{ \int_{-\infty}^\cdot   F_\gamma  + \int_\cdot^{+\infty}   \overline{F}_\gamma  }\abs{ \rho_\gamma + X D_\gamma - \frac{\psi_\gamma}{f_\gamma } B_\gamma   }  \\
              & \hspace{+2.2cm} + \prth{ \int_{-\infty}^\cdot   \varphi_\gamma } \abs{  - D_\gamma + \frac{\overline{F}_\gamma}{f_\gamma } B_\gamma  }   +   \prth{  \int_\cdot^{+\infty}   \overline{\varphi}_\gamma  } \abs{ D_\gamma + \frac{ F_\gamma}{f_\gamma } B_\gamma  } \Bigg) \\
              & = \norm{ h'' }_\infty \Bigg( 1 +  \crochet{ \varphi_\gamma  +  \overline{\varphi}_\gamma  } \abs{ \rho_\gamma + X D_\gamma - \frac{\psi_\gamma}{f_\gamma } B_\gamma   } \\
              & \hspace{+2.2cm} + \chi_\gamma \abs{  - D_\gamma + \frac{\overline{F}_\gamma}{f_\gamma } B_\gamma  } + \overline{\chi}_\gamma \abs{ D_\gamma + \frac{ F_\gamma}{f_\gamma } B_\gamma  } \Bigg)
\end{align*}

Each term in the absolute values of the LHS of this former inequality are symmetric on $ \Rr $, using the parity of $ \rho_\gamma $, \eqref{Ineq:InequalityPsiGamma} and \eqref{Ineq:PositivityDerivativesOfG}, thus, this is enough to restrict ourselves to $ \Rr_+ $.

Let $ x > 0 $. Using the fact that $ 2 \rho_\gamma' + \rho_\gamma^2 \geq 0 $ and $ \rho_\gamma(x) = x (a + b x^2) $, we have, using \eqref{Ineq:InequalityPsiGamma} and \eqref{Ineq:PositivityDerivativesOfG}
\begin{align*}
\frac{\abs{ D^3 \Le_\gamma\inv \widehat{h}_\gamma  } }{ \norm{ h'' }_\infty  } - 1 & \leq     \crochet{ \varphi_\gamma  +  \overline{\varphi}_\gamma  }\prth{ \frac{\psi_\gamma}{f_\gamma } B_\gamma - (\rho_\gamma + X D_\gamma )  }  \\
              & \hspace{+2cm} +     \chi_\gamma   \prth{   D_\gamma - \frac{\overline{F}_\gamma}{f_\gamma } B_\gamma  }   + \overline{\chi}_\gamma  \prth{ D_\gamma + \frac{ F_\gamma}{f_\gamma } B_\gamma }  \\
              & =: \crochet{ \varphi_\gamma  +  \overline{\varphi}_\gamma  } I_1 + \chi_\gamma I_2 +  \overline{\chi}_\gamma  I_3
\end{align*}

We now estimate each of the previous terms.

$ $

\noindent\textbf{\underline{$ (1) $ Estimation of $ \prth{ \varphi_\gamma  +  \overline{\varphi}_\gamma  } I_1 $ :}} \linebreak
%
Using \eqref{Ineq:QGaussienneXRho} and the obvious positivity of $ \rho_\gamma'' + 3 \rho_\gamma \rho_\gamma' + \rho_\gamma^3 $ on $ \Rr_+ $, we have 
\begin{align*}
I_1 & \leq   \frac{ X }{ \rho_\gamma } (\rho_\gamma'' + 3 \rho_\gamma \rho_\gamma' + \rho_\gamma^3 ) - (\rho_\gamma + X(2\rho_\gamma' + \rho_\gamma^2) )   \\
       & = 2 b X \frac{ b X^4 + a X^2 + 3}{ b X^2 + a } 
\end{align*}

Using \eqref{Ineq:InequalityPhiPlusPhiBar}, we have $ 0 \leq  \varphi_\gamma  +  \overline{\varphi}_\gamma \leq 2 f_\gamma / \widetilde{\rho}_\gamma  $, which implies
\begin{align*}
\prth{ \varphi_\gamma  +  \overline{\varphi}_\gamma  } I_1  & \leq  4 b f_\gamma X \frac{ b X^4 + a X^2 + 3}{ (b X^2 + a)^2 }  
\end{align*}

If $ x \in \crochet{ 0, 2 \gamma} $, we have $ a x^2 \leq 4 a \gamma^2 = 4 $, thus
\begin{align*}
\prth{ \varphi_\gamma(x)  +  \overline{\varphi}_\gamma(x)  } I_1(x) \leq   4  \, b \gamma \, f_\gamma(x)  \frac{ b x^4 + 7}{ (b x^2 + a)^2 }
\end{align*}

Using $ f_\gamma' = - \rho_\gamma f_\gamma $, we have 
\begin{align*}
\frac{d}{dx} \prth{  f_\gamma(x)  \frac{ b x^4 + 7}{ (b x^2 + a)^2 } } = -4 f_\gamma(x) \, \frac{ x \left(  b^3 x^8 + 2 a b^2 x^6 + b \left( a^2  + 7 b  \right) x^4 + 10 a b x^2 + 7 a^2 + 28 b \right) }{ \left( b x^2 + a \right)^3 } \leq 0
\end{align*}
thus, the maximum of this last function on $ \crochet{0, 2\gamma} $ is obtained in $ x = 0 $, as it is decreasing. We thus deduce that
\begin{align*}
\prth{ \varphi_\gamma(x)  +  \overline{\varphi}_\gamma(x)  } I_1(x) \leq   4  \, b \gamma \, f_\gamma(0)  \frac{   7}{ a^2 } = 28 \frac{\gamma }{ z_\gamma } \frac{C}{\gamma^8 } \gamma^4 \leq \frac{28 C}{ \sqrt{2\pi} } \frac{1}{\gamma^4 } + \frac{3C }{4 \gamma^8 } \leq \prth{ \frac{28  }{ \sqrt{2\pi}  } + \frac{3}{4 } } \frac{C}{ \gamma^4 } 
\end{align*}
for $ \gamma \geq 1 $.

If $ x \in [2 \gamma, +\infty [ $, define
\begin{align*}
M_\gamma(x) := 	6 b f_\gamma(x) x \frac{ b x^4 + a x^2 + 3}{ (b x^2 + a)^2 }
\end{align*}

Using $ f_\gamma' = - \rho_\gamma f_\gamma $, we have 
\begin{align*}
M_\gamma'(x) = -6 b f_\gamma(x) \, \frac{ P_\gamma(x^2)  }{  ( b x^2 + a )^3 }  
\end{align*}
with
\begin{align*}
P_\gamma(x) := b^3 x^5 + 3 a b^2 x^4 + 3 a^2 b x^3 + a^3 x^2 + 3 b x - 3 a 
\end{align*}

As $ P_\gamma' > 0 $ on $ \Rr_+^* $ with $ P_\gamma(0) = -3a < 0  $ and $ \lim_{+\infty } P_\gamma = +\infty $, there exists a unique $ x_\gamma^* \in \Rr_+ $ such that $ P_\gamma'(x_\gamma^*) = 0 $ ; hence, $ P_\gamma $ is increasing on $ \crochet{0, x_\gamma^*} $ and decreasing on $ [x_\gamma^* , +\infty[ $. As 
\begin{align*}
P_\gamma'(2\gamma) = 80 \frac{ C^3 }{ \gamma^{20} }   + 96 \frac{ C^2 }{ \gamma^{15} } + 36 \frac{C }{ \gamma^{10} }  + 3 \frac{ C }{ \gamma^8 } + \frac{4}{\gamma^5 } \geq 0
\end{align*}
and $ P_\gamma' $ is bijective from $ \Rr_+ $ to $ [3b, +\infty[ $, we deduce that $ x_\gamma^* \leq 2 \gamma $ and in particular, that $ M_\gamma $ is decreasing on $ [2 \gamma, +\infty [ $. We thus have, for all $ x \geq 2 \gamma $
\begin{align*}
M_\gamma(x) & \leq M_\gamma(2\gamma) = 	6 b f_\gamma(2\gamma) 2\gamma \frac{ b (2\gamma)^4 + a (2\gamma)^2 + 3}{ (b (2\gamma)^2 + a)^2 } \leq \frac{6}{z_\gamma} \frac{ 4C( 7 \gamma^4 + 16 C) }{ \gamma^5 ( \gamma^4 + 4 C) } \leq \frac{168 C}{ z_\gamma \gamma^5} \\
             & \leq \frac{ C}{ \gamma^6 } \prth{ \frac{ 168 }{ \sqrt{2\pi } } + \frac{3}{4} }
\end{align*}
for $ \gamma \geq 1 $. This last bound is smaller than the bound on $ \crochet{0, 2\gamma} $ for $ \gamma \geq 2C $.

Finally, using $\frac{28  }{ \sqrt{2\pi}  } + \frac{3}{4 } \approx 11,\!920 \leq 12  $, we have on $ \Rr_+ $ and for $ \gamma \geq 2C $
\begin{align}\label{Ineq:I1}
\prth{ \varphi_\gamma  +  \overline{\varphi}_\gamma  } I_1  \leq   \prth{ \frac{28  }{ \sqrt{2\pi}  } + \frac{3}{4 } } \frac{C}{ \gamma^4 } \leq \frac{12 C}{ \gamma^4 }
\end{align}

\noindent\textbf{\underline{$ (2) $ Estimation of $ \chi_\gamma I_2 + \overline{\chi}_\gamma I_3 $ :}}\linebreak
%
We have 
\begin{align*}
\chi_\gamma(x) + \overline{\chi}_\gamma(x) = \Esp{ \frac{( H_\gamma - x )^2 }{ 2 } \Unens{ H_\gamma - x \leq 0 } } + \Esp{ \frac{( H_\gamma - x )^2 }{ 2 } \Unens{ H_\gamma - x \geq 0 } } = \frac{ x^2 + \Esp{ H_\gamma^2 }}{ 2 } =: V_\gamma(x)
\end{align*}
using definition \eqref{Def:GeneralDefinitions}.

Then, 
\begin{align*}
\chi_\gamma I_2 + \overline{\chi}_\gamma I_3 =   V_\gamma I_2 + \overline{ \chi}_\gamma \prth{ I_2 - I_1 } =  V_\gamma   I_2 + \overline{ \chi}_\gamma \frac{ B_\gamma }{ f_\gamma }
\end{align*}

Using \eqref{Ineq:BoundOnVgammaI2}, we get
\begin{align*}
V_\gamma I_2 \leq 1 + \frac{18 C }{ 10 } \leq 1 + 2C
\end{align*}

Using \eqref{Ineq:BoundOnChiBarBoverf}, we get
\begin{align*}
\overline{ \chi}_\gamma \frac{ B_\gamma }{ f_\gamma } \leq 1
\end{align*}

Finally, combining these last two inequalities, we get
\begin{align}\label{Ineq:I2-3}
\chi_\gamma I_2 + \overline{\chi}_\gamma I_3 \leq 2 + 2C
\end{align}

\noindent\textbf{\underline{$ (3) $ Conclusion :}} \linebreak
%
Using \eqref{Ineq:I1} and \eqref{Ineq:I2-3}, we get the desired bound. 
\end{proof}


\section*{Acknowledgements}

The author expresses his thanks to A. Nikeghbali for guidance and introduction to the topic, and to G. Borot, P. Dey and P.-L. M\'eliot for interesting discussions and remarks related to earlier drafts of this work.
The first version of this paper was written while the author was a guest at the University of Jinan (China) ; many thanks are due to this institution and in particular to S. Peng and A. Nikeghbali. 
Particular thanks are given to J. Najnudel for several discussions and corrections of earlier drafts, to O. H\'enard for careful discussion of the results and to A.D. Barbour for encouragements and helpful criticism that led to \ref{Thm:MASSmodG2}.

The author was supported by EPSRC grant EP/L012154/1 and the Schweizerischer Nationalfonds PDFMP2 134897/1.


\bibliographystyle{amsplain}


\end{document}